\documentclass{amsart}
\usepackage[margin=1.8cm]{geometry}
\usepackage{amssymb}
\usepackage{relsize}
\usepackage{url}
\numberwithin{equation}{section}

\newcommand{\di}{d_i}

\newcommand{\ls}{{(L, \mathcal S)}}

\newcommand{\ali}{\al_{ _i}}

\newcommand{\ovl}{\overline}

\newcommand\numeq[1]%
  {\stackrel{\scriptscriptstyle(\mkern-1.5mu#1\mkern-1.5mu)}{=}}

\newcounter{relctr} 
\everydisplay\expandafter{\the\everydisplay\setcounter{relctr}{0}} 

\AtBeginDocument{} 

\newcommand{\io}{{i_1}}
\newcommand{\itw}{{i_2}}
\newcommand{\jo}{{j_1}}
\newcommand{\jtw}{{j_2}}

\usepackage{tikz-cd}
\usepackage{latexsym,amssymb,amsmath}
\usepackage{hyperref}
\usepackage{comment}
\usepackage{mathrsfs}
\usepackage{amsmath,amscd}
\usepackage[all]{xy}
\usepackage{empheq}

\usepackage{amsmath}

\usepackage{empheq}
\usepackage{xcolor}
\definecolor{lightgreen}{HTML}{90EE90}

\newcommand{\irrcc}{\irr(\cc)}

\newcommand{\lag}{\langle}
\newcommand{\rag}{\rangle}

 \newcommand{\vsk}{\vskip 0.15cm \noindent}

 \newcommand{\fpcc}{\fp(\cc)}

\newtheorem{theorem}{Theorem}[section]
\newtheorem{lemma}[theorem]{Lemma}
\newtheorem{proposition}[theorem]{Proposition}
\newtheorem{corollary}[theorem]{Corollary}
\newtheorem{conjecture}[theorem]{Conjecture}

\newtheorem{question}[theorem]{Question}

\newcommand{\vect}{\mtr{Vec}}

\newcommand{\sent}{\mapsto}

\newcommand{\mui}{\mu_i}\newcommand{\muj}{\mu_j}

\newcommand\C{\mathcal{C}}
\DeclareMathOperator{\id}{id}

\DeclareMathOperator{\ev}{ev}

\DeclareMathOperator{\FPdim}{\mathsf{FPdim}}

\excludecomment{verlong}
\includecomment{vershort}
\excludecomment{noncompile}

\usepackage{tikz}
\usetikzlibrary{matrix}
\usepackage{combelow}
\newcommand{\ccb}{{\mathcal B}}

\usepackage{amsmath,amsthm,amsfonts,amssymb}
\usepackage{amsmath, amsthm, amssymb, amscd, amsfonts}
\usepackage[all]{xy}
\usepackage{amssymb, amsthm, amsmath}
\usepackage{amsfonts}
\usepackage{color}

\newcommand{\ra}{\rightarrow}

\newcommand{\ot}{\otimes}

\newcommand{\xra}{\xrightarrow}
\newcommand{\mtc}{\mathcal}
\newcommand{\cs}{\mtc S}

\newcommand{\lam}{\lambda}

\newcommand{\Lam}{\Lambda}

\newcommand{\al}{\alpha}
\newcommand{\eps}{\epsilon}

\newcommand{\ul}{\underline}

\newcommand{\lh}{\leftharpoonup}
\newcommand{\whb}{{\widehat{(H, \mtc B)}}}
\newcommand{\hb}{{(H, \mtc B)}}
\newtheorem{thm}[theorem]{Theorem}
\newtheorem{rem}[theorem]{Remark}
\newtheorem{defn}[theorem]{Definition}
\newtheorem{lem}[theorem]{Lemma}
\newtheorem{conj}[theorem]{Conjecture}
\theoremstyle{plain}

\newcommand{\ch}{\chi}
\newcommand{\mtr}{\mathrm}

\newcommand{\ncm}{\newcommand}
\ncm{\np}{\newpage}
\ncm{\ebl}{\end{thebibliography}}
\ncm{\bbl}{\begin{thebibliography}}
\ncm{\chd}{_{ _{\ch}}}
\ncm{\ald}{_{ _{\al}}}
\newcommand{\blam}{\Lam}
\ncm{\cP}{\mathcal{P}}
\ncm{\ei}{e_i}
\ncm{\eij}{e_{i,\;j}}
\ncm{\bt}{\begin{thm}}
\ncm{\bdef}{\begin{defn}}
\ncm{\edf}{\end{defn}}
\ncm{\et}{\end{thm}}
\ncm{\bc}{\begin{corollary}}
\ncm{\bl}{\begin{lem}}
\ncm{\el}{\end{lem}}
\ncm{\bpf}{\begin{proof}}
\ncm{\epf}{\end{proof}}
\ncm{\ec}{\end{corollary}}
\ncm{\ord}{\mtr{ord}}
\ncm{\er}{\end{rem}}
\ncm{\br}{\begin{rem}}
\ncm{\bn}{\begin}

\ncm{\bp}{\begin{proposition}}
\ncm{\ep}{\end{proposition}}
\ncm{\bd}{
\begin{document}}
\ncm{\ed}{\end{document}}
\ncm{\beq}{\begin{equation}}
\ncm{\beqn}{\begin{equation*}}
\ncm{\eeq}{\end{equation}}
\ncm{\eeqn}{\end{equation*}}
\ncm{\bea}{\begin{eqnarray}}
\ncm{\eea}{\end{eqnarray}}
\ncm{\beanon}{\begin{eqnarray*}}
\ncm{\eeanon}{\end{eqnarray*}}\ncm{\ek}{\eps|_K}\ncm{\diez}{\#}
\ncm{\bwt}{\bowtie}
\ncm{\cC}{\mtc{C}}\ncm{\cc}{\mtc{C}}
\ncm{\cX}{\mtc{X}}
\ncm{\wt}{\widetilde}
\ncm{\sg}{\sigma}
\ncm{\Rep}{\mathrm{Rep}}
\ncm{\Aut}{\mathrm{Aut}}
\DeclareMathOperator{\Irr}{Irr}
\ncm{\X}{\mathcal{X}}
\ncm{\cA}{\mathcal{A}}
\ncm{\HKer}{\mtr{HKer}}
\ncm{\LKER}{\mtr{LKer}}
\ncm{\aad}{\mtr{ad}}
\newcommand{\mbf}{\mathbb F}
\ncm{\Dr}{\mtr{D}}
\ncm{\cD}{{\mathcal{D}}}\ncm{\cd}{{\mathcal{D}}}\ncm{\ce}{{\mathcal{E}}}
\ncm{\G}{\mathcal{G}}
\ncm{\Dc}{\mtc{D}}
\ncm{\E}{\mtc{E}}
\ncm{\fp}{\mtr{FPdim}}
\ncm{\Vc}{\mtr{Vec}}
\ncm{\cK}{\mtc{K}}
\ncm{\cM}{\mtc{M}}
\ncm{\cE}{\mtc{E}}
\ncm{\cS}{\mtc{S}}

\newcommand{{\ipr}}{i'}
\newcommand{\tomega}{\widetilde{\omega}}

\DeclareMathOperator{\End}{End}
\ncm{\cop}{\mtr{cop}}
\ncm{\op}{\mtr{op}}
\ncm{\chr}{character }\ncm{\ck}{\mtc{K}}
\ncm{\bw}{\bwt}
\ncm{\hker}{\mtr{HKer}}
\ncm{\bx}{\boxtimes}
\ncm{\blue}{\textcolor[rgb]{.00, .00, 1.00}}
\ncm{\bluer}{\textcolor[rgb]{.30, .30, .70}}
\ncm{\red}{\textcolor[rgb]{1.00, .00, .00}}
\ncm{\green}{\textcolor[rgb]{.50, 0.20, .90}}
\ncm{\bne}{\begin{enumerate}}
\ncm{\ene}{\end{enumerate}}
\ncm{\lker}{\mtr{LKer}}
\ncm{\md}{\medbreak}
\ncm{\rep}{\Rep}\ncm{\ind}{\mtr{ind}}
\ncm{\mdn}{\md\noindent}
\ncm{\dd}{$}
\ncm{\up}{^}
\newcommand{\tcs}{\text}
\newcommand{\mbb}{\mathbb B}
\newcommand{\vs}{\mathbb V}
\newcommand{\sth}{suppose that\;}
\newcommand\rad{\operatorname{rad}}
\newcommand{\itm}{\item}
\newcommand{\dbd}{$$}
\newcommand{\mol}{\mtr{mod}}
 \newcommand{\ro}{\rho}
\newcommand{\irr}{\mathrm{Irr}}
\newcommand{\mbc}{\mathbb C}
\newcommand{\mbs}{\mathbb S}
\newcommand{\mbz}{\mathbb Z}
\newcommand{\ct}{\mtc T}
\newcommand{\sm}{\setminus}
\newcommand{\epl}{^{+}}
\newcommand{\sbsq}{\subseteq}
\newcommand{\sbs}{\subset}
\newcommand{\cco}{\mtr{co}}
\newcommand{\cz}{\mathcal{Z}}
\newcommand{\dual}{^{*}}
\newcommand{\Gm}{\Gamma}
\ncm{\cY}{\mtc{Y}}
\newcommand\ZZ{{\mathbb Z}} 
\newcommand{\bab}{\color{DarkOrchid}{}}
\newcommand{\eab}{\normalcolor{}}
\newcommand{\subs}{\subsection}
\newcommand{\cv}{\mtc{V}}
  \newcommand{\grn}{\green}
\newcommand{\dt}{\delta}

\newcommand{\ccf}{\mathrm{ {CF}(\cc)}}
\newcommand{\cce}{\mathrm{ {CE}(\cc)}}
\newcommand{\cecc}{\mathrm{ {CE}(\cc)}}
\newcommand{\cecd}{\mathrm{ {CE}(\cd)}}
\newcommand{\kk}{\Bbbk}
\newcommand{\otL}{\ot_{L}}
\newcommand{\otl}{\ot_{L}}
\newcommand{\unpsi}{1_{\psi}}
\newcommand{\epsi}{e_{\psi}}
\newcommand{\ephi}{e_{\phi}}
\newcommand{\ech}{e_{\ch}}
\newcommand{\nleftcid}{\text{left normal  coideal subalgebra}}
\newcommand{\dimL}{\dim_{\kk}L}
\newcommand{\cl}{\mtc L}
\newcommand{\mj}{\mtc J}
\newcommand{\tl}{\tilde L}
\newcommand{\tL}{\tilde L}
\newcommand{\tpsi}{\tilde(\psi)}
\newcommand{\tmx}{\tilde{\mtc X}}
\newcommand{\zlh}{\mathrm{ZL}}
\newcommand{\ba}{\mathrm A}
\newcommand{\bv}{\mathrm V}
\newcommand{\zhopf}{\mtc{Z}_{\mtr{Hopf}}}
\newcommand{\lstar}{L^{*}}
\newcommand{\ldstar}{L^{**}}
\newcommand{\mstar}{M^{*}}
\newcommand{\mdstar}{M^{**}}
\newcommand{\lkera}{\lker_{A}}
\newcommand{\mdprime}{M''}
\newcommand{\ldprime}{L''}
\newcommand{\cm}{\mtc M}
\newcommand{\ccm}{\mathcal M}
\newcommand{\cn}{\mathcal N}
\newcommand{\ccn}{\mathcal N}
\newcommand{\rx}{\mtr{Rex}}
\newcommand{\cca}{\ca}
\newcommand{\ih}{\underline{\mtr{Hom}}}
\newcommand{\cih}{\underline{\mtr{coHom}}}
\newcommand{\hm}{\mtr{ {Hom}}}
\newcommand{\cov}{\mtr{coev}}
\newcommand{\rora}{\rho^{\mtr{ra}}}
\newcommand{\rola}{\rho^{\mtr{la}}}
\newcommand{\cx}{\mtc X}
 \newcommand{\cZ}{\cz}
 \newcommand{\ca}{\cA}
 \newcommand{\stat}{\noindent}
 \newcommand{\bfa}{{\bf A}}
 \newcommand{\unu}{\mathbf{1}}
 \newcommand{\barzu}{{\bar {  Z}(\unu)}}
 
\newcommand{\idx}{\id_X}
\newcommand{\lprime}{L'}
\newcommand{\mprime}{M'}
\newcommand{\nat}{ \mtr{{  Nat}}}
\newcommand{\ft}{\mtc F_\lam}
\newcommand{\rhau}{\rightharpoonup}
\newcommand{\lhau}{\leftharpoonup}
\newcommand{\cf}{\mathrm{ {CF}}}

\newcommand{\cfc}{\mathrm{{CF}}(\cc)}
\newcommand{\csu}{\overline{\mathfrak{  C}}}
\newcommand{\cfcc}{{\mathrm{CF}(\cc)}}
\newcommand{\catfcc}{\mathrm{ {CF}}(\cc)}
\newcommand{\cfcd}{\mathrm{CF}(\cd)}
\newcommand{\cfd}{\mathrm{CF}(\cd)}
\newcommand{\czcc}{{\cz(\cc)}}
\newcommand{\czcd}{{\cz(\cd)}}
\newcommand{\czt}{{\cz(\cz(\cc))}}
\newcommand{\enx}{\mtr{  End}}
\newcommand{\runu}{R(\unu)}

\newcommand{\bdfn}{\bn{defn}}
\newcommand{\edfn}{\end{defn}}
\newcommand{\deltax}{\delta_X}
\newcommand{\deltav}{\delta_V}
\newcommand{\repcca}{\rep_\cc(A)}
\newcommand{\xotay}{X \ot_A Y}
\newcommand{\xoty}{X \ot Y}
\newcommand{\votw}{V \ot W}
\newcommand{\votaw}{V \ot_A W}
\newcommand{\dimax}{\dim_AX}
\newcommand{\dimccx}{\dim_\cc(X)}
\newcommand{\dimcca}{\dim_\cc(A)}
\newcommand{\dimccv}{\dim_\cc(V)}
\newcommand{\dima}{\dim_A}
\newcommand{\biga}{A}
\newcommand{\comp}{\mathbb C}
\newcommand{\tehtaa}{\theta_A}
\newcommand{\tetaa}{\theta_A}
\newcommand{\ida}{\id_A}
\newcommand{\hma}{\hm_A}
\newcommand{\hmcc}{\hm_\cc}
\newcommand{\fv}{F(V)}
\newcommand{\fw}{F(W)}
\newcommand{\ota}{\ot_A}
\newcommand{\repza}{\rep_\cc^0(A)}
\newcommand{\epsa}{\eps_A}
\newcommand{\bndefn}{\bn{defn}}
\newcommand{\edefn}{\end{defn}}
\newcommand{\bdefn}{\bn{defn}}

\newcommand{\vld}{V^{*}}
\newcommand{\vldd}{V^{**}}
\newcommand{\xld}{X^{*}}
\newcommand{\xldd}{X^{**}}
\newcommand{\yld}{Y^{*}}
\newcommand{\yldd}{Y^{**}}
\newcommand{\aldu}{A^{*}}
\newcommand{\aldd}{A^{**}}

\newcommand{\ia}{\mtr{i}_A}
\newcommand{\aota}{A\ot A}

\newcommand{\idv}{\id_V}

\newcommand{\ld}{^*}
\newcommand{\repg}{\rep(G)}

\newcommand{\thetav}{\theta_V}

\newcommand{\tta}{\theta_A}

\newcommand{\muv}{\mu_V}
\newcommand{\muw}{\mu_W}

\newcommand{\dimcc}{\dim(\cc)}
\newcommand{\chii}{\chi_i}
\newcommand{\chistar}{\ch_{i^*}}
\newcommand{\chj}{\ch_j}
\newcommand{\chm}{\ch_m}
\newcommand{\chn}{\ch_n}
\newcommand{\dimvi}{\dim(V_i)}
\newcommand{\mtcd}{Q}
\newcommand{\mtca}{\mtc A}
\newcommand{\lamcd}{\lam_\cd}
\newcommand{\fpdimcd}{\fp(\cd)}
\newcommand{\laml}{\lam_L}
\newcommand{\apm}{A//M}
\newcommand{\apl}{A//L}
\newcommand{\repapm}{\rep(\apm)}
\newcommand{\repapl}{\rep(\apl)}
\newcommand{\dimvj}{\dim(V_j)}
\newcommand{\dvi}{\dim(V_i)}
\newcommand{\dvj}{\dim(V_j)}
\newcommand{\sumjtom}{\sum_{j \in \mathcal I}}
\newcommand{\sumitom}{\sum_{i \in \mathcal I}}
\newcommand{\sij}{s_{ij}}
\newcommand{\sji}{s_{ji}}
\newcommand{\dxj}{d_j}
\newcommand{\dxi}{\di }
\newcommand{\dimka}{\dim_{\kk}(A)}
\newcommand{\dimk}{\dim_{\kk}}
\newcommand{\blaml}{\blam_L}
\newcommand{\sumjtor}{\sum_{j=1}^r}
\newcommand{\dimkl}{\dim_{\kk}(L)}
\newcommand{\mtcjl}{\mtc J_L}
\newcommand{\vota}{ V\ot A}
\newcommand{\vi}{V_i}
\newcommand{\vj}{V_j}
\newcommand{\dimcd}{\dim(\cd)}

\newcommand{\alij}{{\al_{ _{ij}}}}
\newcommand{\alji}{{\al_{ _{ji}}}}
\newcommand{\rcc}{r_\cc}
\newcommand{\rcd}{r_\cd}
\newcommand{\clsx}{[X]}
\newcommand{\clsy}{[Y]}
\newcommand{\clsz}{[Z]}
\newcommand{\rcdp}{r_{\cd'}}
\newcommand{\sumjtorp}{\sum_{j=1}^{r'}}
\newcommand{\aljm}{{\al_{ _{jm}}}}
\newcommand{\aljn}{{\al_{ _{jn}}}}
\newcommand{\sjm}{s_{jm}}
\newcommand{\smj}{s_{mj}}
\newcommand{\snj}{s_{nj}}

\newcommand{\betaij}{\beta_{ _{ij}}}
\newcommand{\betaji}{\beta_{ _{ji}}}
\newcommand{\gammaij}{\gamma_{ _{ij}}}
\newcommand{\gammaji}{\gamma_{ _{ji}}}
 \newcommand{\ip}{i'}
\newcommand{\sumjtoprp}{\sum_{j=1}^{r'}}
\newcommand{\sumjtopr}{\sum_{j=1}^{r}}
 \newcommand{\teh}{\tilde{h}}
\newcommand{\cdp}{{\cd'}}\newcommand{\xphii}{X_{\phi(i)}}
\newcommand{\inv}{^{-1}}

\newcommand{\fq}{\mtr f_{ Q}}
\newcommand{\tr}{\mtr{tr}}
\newcommand{\rtwone}{R_{21}R}

\newcommand{\ccad}{{\cc_{\mtr{ad}}}}
\newcommand{\ccpt}{{\cc_{\mtr{pt}}}}
\newcommand{\qtr}{quasi-triangular\;}
\newcommand{\trq}{\tr_q}

\newcommand{\repal}{\mtr{Rep}(A//L)}
\newcommand{\lkeravi}{\lker_A(V_i)}
\newcommand{\lkeravj}{\lker_A(V_j)}
\newcommand{\cross}[1][1pt]{\ooalign{%
 \rule[1ex]{1ex}{#1}\cr
 \hss\rule{#1}{.7em}\hss\cr}}
\newcommand{\blml}{\blam_L} 
\newcommand{\phir}{\phi_R}
\newcommand{\kda}{{  \Phi(A)}}

\newcommand{\mtcil}{\mtc{I}_L}

\newcommand{\un}{\unu}
\newcommand{\tfl}{\mtc{T}}
\newcommand{\barzm}{\barz(M)}
\newcommand{\barzn}{\barz(N)}
\newcommand{\ccr}{\mtc R^{\cc}}
\newcommand{\ulc}{\ul{\cc}}

\newcommand{\pimx}{\pi_{M;\;X}}
\newcommand{\pinx}{\pi_{N;\;X}}
\newcommand{\acc}{{\mathrm A_\cc}}
\newcommand{\epsu}{\eps_\unu}

\newcommand{\ob}{\mtr{Obj}}
\newcommand{\obc}{\mtr{Obj(\cc)}}
\newcommand{\ccop}{\cc^{\mtr{op}}}
\newcommand{\mtf}{\mtc F}
\newcommand{\mtfi}{\mtc F^{-1}_\lam}
\newcommand{\elcd}{\ell_\cd}
\newcommand{\mcid}{\mtc I_\cd}
\newcommand{\mcidp}{\mtc I_{\cd'}}
\newcommand{\wtildelcd}{\widetilde{\elcd}}
\newcommand{\wtildelcdp}{\widetilde{\ell_{\cd'}}}
\newcommand{\cpt}{\cc_{\mtr{pt}}}
\newcommand{\barzr}{\barz_\cd}
\newcommand{\barzv}{\barz(V)}
\newcommand{\acd}{\mathrm A_\cd}
\newcommand{\czrcd}{\cz_\cc(\cd)}
\newcommand{\sml}{\Small}
\newcommand{\bs}{{\Small }}
\newcommand{\yd}{Yetter-Drinfeld}

\newcommand{\sumitor}{\sum_{i=1}^r}
\newcommand{\cdop}{\cd^{\mtr{op}}}
\newcommand{\ccrev}{\cc^{\mtr{rev}}}
\newcommand{\barz}{{\bar{\mathrm Z}}}
\newcommand{\etl}{etale\;}
\newcommand{\czca}{\cz(\ca)}

\newcommand{\tetx}{\text}
\newcommand{\widehta}{\widehat}
\newcommand{\wdhat}{\widehat}
\newcommand{\wht}{\widehat}
\newcommand{\cofa}{{\mathbb C[\mtc B]}}
\newcommand{\wdt}{\widehat}
\newcommand{\dl}{{^\#}}
\newcommand{\comx}{\mathbb C}

\newcommand{\sgj}{{\sg(j)}}

\newcommand{\mujo}{\mu_\jo}
\newcommand{\mujtw}{\mu_\jtw}
\newcommand{\adz}{a^{\#}}
\newcommand{\bdz}{b^{\#}}

\newcommand{\spr}{S^\perp}
\newcommand{\cofs}{\comp [S]}
\newcommand{\spz}{S^{\perp_z}}

\newcommand{\omz}{\omega_z}
\newcommand{\zg}{\mathrm{Z}(S)}
\newcommand{\aling}{{\al \in g}}

\newcommand{\blkg}{\mtr{Bl}(g)}
\newcommand{\clsg}{\mtr{Cl}(g)}
\newcommand{\mtadinv}{\mtc G^{{-1}}}
\newcommand{\muk}{\mu_{k}}
\newcommand{\mta}{\mtc F}
\newcommand{\cofad}{\comp[\wdht A]}
\newcommand{\wtau}{\wdht{\tau}}
\newcommand{\mtainv}{{\mta}^{-1}}
\newcommand{\wdht}{\widehat}
\newcommand{\augm}{\mtr{aug}}
\newcommand{\mua}{\wdht {\wdht a}}
\newcommand{\aps}{A//S}
\newcommand{\cssa}{\cc(S, A)}
\newcommand{\aug}{\mtr{aug}}
\newcommand{\rss}{{\big|_S}}
\newcommand{\gprp}{g^\perp}
\newcommand{\alins}{{s \in S}}

\newcommand{\sz}{s^{D}}
\newcommand{\wmu}{\widehta{\mu}}
\newcommand{\wmui}{\widehta{\mu}_i}
\newcommand{\wmuj}{\widehta{\mu}_j}
\newcommand{\wch}{\widehta{\ch}}

\newcommand{\wzd}{\widehat{d}}
\newcommand{\wpm}{\widehat{P}}
\newcommand{\wps}{\widehat{p}}

\newcommand{\gal}{\mtr{Gal}}
\newcommand{\galkq}{\gal(\mathbb K/\mathbb Q)}
\newcommand{\sgh}{\sg_{ _{H}}}
\newcommand{\sggi}{{\sg(i)}}
\newcommand{\sge}{\sg_{_{\widehat R}}}
\newcommand{\unue}{{\unu_{\cecc}}}

\newcommand{\mtcf}{\mtc {F}}

\newcommand{\wsgf}{\widehat{{\sg}_{ _F}}}
\newcommand{\sghstar}{{{\sg}_{ _{H^*}}}}
\newcommand{\we}{\widehta{E}}
\newcommand{\sumktom}{\sum_{k=1}^m}

\newcommand{\wf}{\widehat{F}}

\newcommand{\hsgj}{\widehat{\sg}(j)}
\newcommand{\whsgi}{\widehta{\sg}(i)}

\newcommand{\wpp}{\widehat{p}}
\newcommand{\tauj}{{\tau(j)}}
\newcommand{\dimcctauj}{\dim(\cc^\tauj)}
\newcommand{\etas}{{\eta(s)}}
\newcommand{\mcc}{m_H}

\newcommand{\wal}{\widehta{\al}}
\newcommand{\wj}{\widehat{\mtc J}}
\newcommand{\galc}{\mtr{Gal}_{\cc}}
\newcommand{\galz}{\mtr{Gal}_{\czcc}}
\newcommand{\wjr}{\widehat{J}_{R}}

\newcommand{\dimcck}{\dim(\cc^k)}

\newcommand{\wgrcc}{\widehat{\mtr{Gr}(\cc)}}
\newcommand{\nchi}{{\frac{\ch_i}{\di }}} \newcommand{\nchj}{{\frac{\ch_j}{\dxj}}}
\newcommand{\wni}{{\widehat{n}_i}}

\newcommand{\sgte}{\widetilde{\sg_E}}

\newcommand{\mtad}{\mtc G}
\newcommand{\whj}{\widehta{h}_j}
\newcommand{\jdl}{{j\dl}}
\newcommand{\wcfcc}{\widehat{\cfcc}}
\newcommand{\mutauj}{\mu_{\tau(j)}}
\newcommand{\tauk}{\tau(k)}
\newcommand{\muzm}{{\mu_1^{-}}}
\newcommand{\sqrtog}{\sqrt{|G|}}
\newcommand{\muz}{\mu_1}
\newcommand{\njtw}{n_\jtw}
\newcommand{\njo}{n_\jo}
\newcommand{\fjo}{F_\jo}
\newcommand{\fjtw}{F_\jtw}
\newcommand{\wta}{\widehat{A}}

\newcommand{\dol}{{^{\circ}}}
\newcommand{\bdl}{{b\dl}}
\newcommand{\jdol}{{j\dol}}
\newcommand{\fj}{F_j}

\newcommand{\cwta}{\comp[\wta]}

\newcommand{\hx}{\widehta{x}}
\newcommand{\hy}{\widehta{y}}

\newcommand{\cal}{\mtc A_{\al}}
\newcommand{\xuu}{x_{uu}}
\newcommand{\wxuu}{\widehat{\xuu}}
\newcommand{\xvv}{x_{vv}}
\newcommand{\xuv}{x_{uv}}
\newcommand{\xmn}{x_{m,n}}
\newcommand{\buvmn}{B^{u,v}_{m,n}}
\newcommand{\blm}{\blam}
\newcommand{\dimccr}{\dim(\cc^r)}
\newcommand{\adl}{a\dl}
\newcommand{\sumltom}{\sum_{l=1}^m}

\newcommand{\mbq}{\mathbb Q}
\newcommand{\mbqs}{\mathbb Q(S)}
\newcommand{\mbk}{\mathbb K}
\newcommand{\mz}{\mathbb Z}

\newcommand{\wsgj}{\widehat{\sigma}(j)}
\newcommand{\wsgi}{\widehat{\sigma}(i)}
\newcommand{\wg}{\widehat{g}}
\newcommand{\wtf}{\widehat{F}}
\newcommand{\galqspq}{\mtr{Gal}(\mathbb Q(S)/\mathbb Q)}
\newcommand{\cctauj}{\cc^{\tau(j)}}
\newcommand{\cctauk}{\cc^{\tau(k)}}
\newcommand{\wtfj}{\widetilde{F_j}}
\newcommand{\wfj}{\widetilde{F_j}}
\newcommand{\wtmuj}{\widetilde{\mu_j}}
\newcommand{\wmtcfj}{{\widetilde{\mtc F}_j}}
\newcommand{\mtfr}{\mtr{F_a}}
\newcommand{\wdr}{R_\comp^*}

\newcommand{\fgph}{{F_{G/H}}}
\newcommand{\wcfj}{\wmtcfj}

\newcommand{\nxi}{{\frac{x_i}{\di }}}
\newcommand{\fpr}{{\fp(R)}}
\newcommand{\nxs}{{\frac{x_s}{d_s}}}

 \newcommand{\mtfme}{\mtc F}
\newcommand{\chic}{\ch_i^{\circ}}
\newcommand{\chjc}{\ch_j^{\circ}}
\newcommand{\mtfsh}{{\mtc F_\lam}}
\newcommand{\mupq}{{\mu_{pq}}}

 \newcommand{\tlam}{{\widetilde{\lam}}}
 \newcommand{\chid}{{\ch_i^{\circ}}}
\newcommand{\rc}{{R_\comp}}
\newcommand{\rgo}{{\mathbb R_{\geq 0}}}
\newcommand{\sumrorc}{{{\sum\limits_{\ro \in \rc}}}}
\newcommand{\aliro}{{\ro(x_i)}}

\newcommand{\barjd}{{\bar{\mtc J_\cd}}}
\newcommand{\lbarcj}{{\frac{C_j}{{\dim(\mathcal C^j)}}}}
\newcommand{\omtcb}{{\overline{\mathcal B}}}
\newcommand{\whr}{{\widehat{R}}}
\newcommand{\nxj}{\frac{x_j}{\dxj}}
\newcommand{\nxk}{\frac{x_k}{d_k}}
\newcommand{\onkij}{{\overline{N^k_{ij}}}}
\newcommand{\sgk}{{\sigma(k)}}
\newcommand{\sgl}{{\sigma(l)}}
\newcommand{\fqi}{{\fq^{-1}}}
\newcommand{\wdb}{{\widehat{\mtc B}}}
\newcommand{\mtcb}{{\mtc B}}

\newcommand{\nif}{{h_i}}
\newcommand{\rb}{(R, \mtc B)}
 \newcommand{\nxip}{{\frac{x_{i'}}{d_{i'}}}}
\newcommand{\mujp}{{\mu_{j'}}}

\newcommand{\etai}{{\eta(i)}}
\newcommand{\wsg}{{\widehat{\sg}}}
\newcommand{\wsgh}{{\wsg_{ _{H}}}}
\newcommand{\wtaui}{{\widehat{\tau}(i)}}
\newcommand{\wsghstar}{{\wsg_{H^*}}}
\newcommand{\wtauj}{{\wtau(j)}}
\newcommand{\weta}{{\widehat{\eta}}}
\newcommand{\detai}{{d_{ _{\eta(i)}}}}

\newcommand{\hbz}{{(H, \mtc B, \mu_1)}}
\newcommand{\whbz}{{\widehat{\hbz}}}
\newcommand{\tsgh}{{{\widetilde{\sgh}}}}
\newcommand{\ghb}{{G\hb}}
\newcommand{\taujo}{{\tau(\jo)}}
\newcommand{\taujtw}{{\tau(\jtw)}}
\newcommand{\mutauk}{{\mu_{\tauk}}}

\newcommand{\hetai}{{h_{ _\etai}}}
\newcommand{\xetai}{{x_{ _\etai}}}
\newcommand{\wn}{{\widehat{n}}}
\newcommand{\wh}{{\widehat{h}}}
\newcommand{\distar}{{d_{i^*}}}
\newcommand{\dwtaui}{{d_{ _{\wtaui}}}}
\newcommand{\sumiptom}{{\sum_{\ip=1}^m}}
\newcommand{\alitaugj}{{\al_{ _{i\tau_g(j)}}}}
\newcommand{\dtaugj}{{d_{ _{\tau_g(j)}}}}
\newcommand{\mip}{{M(\ip)}}
 \newcommand{\sumttom}{{\sum_{t=1}^m}}
 \newcommand{\muxi}{{\mu_{ _{[X_i]}}}}
\newcommand{\muxip}{{\mu_{ _{[X_{\ip}]}}}}
\newcommand{\ncj}{{\frac{C_j}{\dim(\cc^j)}}}
\newcommand{\jp}{{j'}}
\newcommand{\minv}{{M^{-1}}}
\newcommand{\tfq}{\widehat{\mtr{f}}_Q} 
\newcommand{\catcecc}{{\mtr{CE}(\cc)}}
\newcommand{\wcatfcc}{{\widehat{\catfcc}}}
\newcommand{\mforall}{{\;\;\text{for all}\;\;}}
\newcommand{\what}{\widehat}
\newcommand{\wfz}{{\widehat{F}_1}}
\newcommand{\hbfr}{{(H, \mtc B, \fp)}}
\newcommand{\sgn}{{\mtr{sgn}}}
\newcommand{\wir}{{\widehat{R}}}
\newcommand{\gcc}{{G(\cc)}}
\newcommand{\jccpt}{{{\mtc I}_{ _{\ccpt}}}}
\newcommand{\jccad}{{{\mtc I}_{ _{\ccad}}}}
\newcommand{\fpccad}{{\fp(\ccad)}}
\newcommand{\fpccpt}{{\fp(\ccpt)}}
\newcommand{\kc}{{K(\cc)}}
\newcommand{\wkc}{{\widehat{\kc}}}

\newcommand{\hs}{{(L,\mtc S)}}
\newcommand{\rbad}{{H_{ _{ad}}}}
\newcommand{\hbad}{{\hb_{ad}}}
\newcommand{\jhbad}{{{\mtc I}_{\hbad}}}
\newcommand{\htt}{{(K, \mtc T)}}
\newcommand{\jhtt}{{\mtc I_{ _{\htt}}}}
\newcommand{\jhs}{{\mtc I_{ _{\hs}}}}
\newcommand{\lamhs}{{\lam_{ _{\hs}}}}
\newcommand{\lamhtt}{{\lam_{ _{\htt}}}}

\newcommand{\coo}{{co}}
\newcommand{\wdhad}{{(\wdh)_{ad}}}
\newcommand{\had}{{H_{ad}}}
\newcommand{\wdh}{\widehat{H}}
\newcommand{\whbad}{{\whb_{ _{ad}}}}
\newcommand{\kerhb}{{\ker_{ _{\hb}}}}
\newcommand{\gwdh}{{G(\wdh)}}
\newcommand{\nxl}{\frac{x_l}{d_l}}

\newcommand{\proditom}{{\prod_{i=1}^m}}

\newcommand{\wdhn}{{\wdh^{(n)}}}
\newcommand{\hn}{{H_{(n)}}}
\newcommand{\nox}{{\frac{x}{\fp(x)}}}
\newcommand{\mujstar}{{\mu_{j^\#}}}
\newcommand{\sco}{{S^\coo}}
\newcommand{\rrad}{{I(1)}}
\newcommand{\istar}{{i^*}}
\newcommand{\mtcs}{{\mtc S}}

\newcommand{\qghb}{{\ghb}}
\newcommand{\wqghb}{{{G\whb}}}
\newcommand{\nxp}{{\frac{x_p}{d_p}}}
\newcommand{\nxm}{{\frac{x_m}{d_m}}}

\newcommand{\whp}{{\widehat{P}}}
\newcommand{\prodjtom}{{\prod_{j=1}^m}}
\newtheorem{statement}[theorem]{Statement}
\newcommand{\wghb}{{G\whb}}
\newcommand{\gwhb}{{G\whb}}
\newcommand{\invisible}[1]{\iffalse #1 \fi}

\title{Burnside type results for fusion categories}

\author{Sebastian Burciu}
\address{Inst.\ of Math.\ ``Simion Stoilow" of the Romanian Academy P.O. Box 1-764, RO-014700, Bucharest, Romania}
\email{sebastian.burciu@imar.ro}

\author{Sebastien Palcoux}
\address{S. Palcoux, Beijing Institute of Mathematical Sciences and Applications, Huairou District, Beijing, China}
\email{sebastien.palcoux@gmail.com}
\urladdr{https://sites.google.com/view/sebastienpalcoux}
\date{\today}
\subjclass{18M20; 20N20; 20C15; 18N25; 16T20; 16T30}
\bd
\thanks{The first author is supported by a grant of the Ministry of Research, Innovation and Digitization, CNCS/CCCDI - UEFISCDI, project number PN-III-P4-ID-PCE-2020-0878, within PNCDI III. The second author is supported by BIMSA Start-up Research Fund,  Foreign Youth Talent Program from the Ministry of Sciences and Technology of China and National Natural Science Foundation of China (NSFC, Grant no. 12471031)}
\maketitle
\begin{abstract}
In this paper, we extend a classical vanishing result of Burnside from the character tables of finite groups to the character tables of commutative fusion rings, or more generally to a certain class of abelian normalizable hypergroups. We also treat the dual vanishing result. We show that any nilpotent unitary fusion categories satisfy both Burnside's property and its dual. Using Drinfeld's map, we obtain that the Grothendieck ring of any weakly-integral modular fusion category satisfies both properties. As applications, we prove new identities that hold in the Grothendieck ring of any weakly-integral fusion category satisfying the dual-Burnside's property, thus providing new categorification criteria. In particular we improve \cite[Theorem 4.5]{o-yu} as follows: A weakly integral modular fusion category of $\fp$ $md$ with $d$ square-free coprime with $m$ and $\fp(X)^2$ for every simple object $X$, has a pointed modular fusion subcategory of $\fp$ $d$. We also present new results on perfect modular fusion categories, including a Cauchy-type theorem.

\end{abstract}
\setcounter{tocdepth}{1}
\tableofcontents
\section{Introduction}\label{introd}
A classical result of Burnside in the character theory of finite groups states that any  irreducible non-linear character  of a finite group vanishes on at least one element of the group. This can be stated as follows: in the character table of a finite group, the row of every character of degree $\neq 1$ contains a zero entry.

More recently, the dual version of this result was also intensively studied in the literature. This dual version describes the group elements that vanish on at least one irreducible character, see \cite{inw} and the references therein.
Burnside's result was extended to every weakly integral fusion category with a commutative Grothendieck ring, initially in the context of modular categories in \cite[Appendix]{gnn}, and subsequently in the general framework in \cite{b-galois}.

The main goal of this paper is to develop an analogue of Burnside's result for a certain large class of fusion rings, or even more generally, to a certain class of abelian rational normalizable hypergroups. We also consider  the dual version of Burnside's result in this more general settings.


Recall  that in \cite{b-blms} the author introduced the notion of dual of a Grothendieck ring. This was achieved based on work of Harrison \cite{hdk} on dualizable probability groups and on \cite{zz}. More generally, this notion of dual of a fusion ring was extended to arbitrary fusion rings in \cite{b-palcoux}. It was noticed in \cite{b-blms} that the dual of the Grothendieck ring of a pivotal fusion category is isomorphic to the center of the category as defined in \cite{scalg}. 

In this paper we show that both Burnside's property and its dual are related with the ring structure of the dual of the involving Grothendieck/fusion ring. In general, the dual of a fusion ring  is no longer a fusion ring, but an \emph{abelian normalizable hypergroup}, see \cite{b-blms}. The advantage of the class of hypergroups is to be closed under duality, see \cite[Proposition 2.3]{b-blms}.

For a finite dimensional $\comp$-algebra $H$, fix a (finite) $\mathbb C$-linear basis $\mtc B$ with the unit $1\in \mtc B$ and write $ab=\sum_{c\in \mtcb} N^c_{ab}c$ for some scalars $N^c_{ab}\in \comp$. Note that $H=\cofa$, the complex linear span of $\mtc B$.
\bn{defn} \label{hyp:defn}
A finite dimensional algebra $H$ together  with a $\comp$-linear basis $\mathcal B$ is called a \emph{hypergroup} if there is an involution map $(-)^*:\mtcb\ra \mtcb$ such that   $N^1_{aa^*}> 0$ and $N^1_{ab}= 0$ if $b\neq a^* $. The elements $b\in \mtc B$ are called the \emph{standard elements} of $\hb$.
\end{defn}
We may refer to $\mtcb$ (or $H$) only as a hypergroup, when $H$ (or $\mtcb$) is implicitly understood.
A hypergroup $\mtc B$ is called:
\begin{itemize}
\item \emph{symmetric} if $N^1_{ab}=N^1_{ba}$,
\item \emph{normalized} if $\sum_{c\in \mtcb} N^c_{ab} =1$,
\item \emph{real} if $N_{ab}^c\in \mathbb R$,
\item \emph{rational} if $N_{ab}^c\in \mathbb Q$,
\item \emph{real non-negative} if $N_{ab}^c \geq 0$,
\item \emph{abelian} if $ab=ba$,
\end{itemize} 
for all $a,b,c \in \mtcb$. Throughout this paper we use the abbreviation (A)RN-hypergroup for an (abelian) real non-negative hypergroup.

For a normalized hypergroup, the linear map $\mu_1:H\ra \comp$ defined on the linear basis $\mtcb$ by $ \mu_1(a)= 1$ for all $a\in \mtcb$, is a morphism of algebras. It is called the \emph{augmentation} map of $(H, \mtc B)$. 
The underlying algebra $H$ of a hypergroup $(H, \mtc B)$ is a finite dimensional $*$-algebra, so is semisimple. For a normalized hypergroup $\hb$ the central primitive idempotent $F_1$ corresponding to the character $\mu_1$ is denoted by $\lam_\hb$ or simply $\lam_H$ if $\mtc B$ is implicitly understood.  A formula for this idempotent is described in Corollary \ref{fz:gen}.

A RN-hypergroup is sometimes called a \emph{table algebra} in the literature.  Recall that a \emph{fusion ring} $\mtcb$ is a hypergroup such that $N^a_{bc}\in \mathbb Z_{\geq 0}$ and $N^1_{aa^*}=N^1_{a^*a}=1$ for all $a,b,c \in \mtc B$. 
The first main result of this paper is the following generalization of Burnside's result to some weakly-integral fusion rings:
\bt\label{burnside-fr}
Let $\hb$ be a commutative fusion ring with a $h$-integral dual. 
Let $a$ be a standard element such that $\fp(a)>1$. Then there is  $\mu\in \wdb$ such that $\mu(a)=0$. 
\et  
A fusion ring whose dual is h-integral (see Definition \ref{h:integral}) is always weakly-integral ($\fp\hb\in \mathbb Z$). The converse is true for a weakly-integral fusion category (see Remark \ref{owh}). The set $\wdb$ is the set of characters $\muj:H\ra \comp$ (see \S \ref{afr} for more details). In fact we prove a more general version of Theorem \ref{burnside-fr}, for a certain class of abelian rational normalizable hypergroups, see Theorem \ref{burnside}.
\bn{defn}\label{grplike:def}
For any hypergroup $\hb$, an element $a\in \mtc B$ with 
\beq\label{grplk:id:eq}
a a^* = N_{a,a^*}^1 1  \text{ (or equivalently, }h_a a a^* = 1\text{, with }h_a:=1/N_{a,a^*}^1)
\eeq 
is called a \emph{grouplike element} of $\hb$. In other words, $N_{a,a^*}^c = \delta_{1,c} N_{a,a^*}^1$.
\end{defn}
By Lemma \ref{gh}, the set of grouplike elements form a group (with the multiplication structure) denoted $\ghb$.
\bn{defn}\label{v-property}
A hypergroup $\hb$ has \emph{Burnside's vanishing property} (or shortly, is Burnside) if for all $a\in  \mtc B$, the following are equivalent:  
\bne
\item There is some $\mu\in \wdb$ such that 
$\mu(a)=0$,
\item
the standard element $a$ is not a grouplike element.
\ene
\end{defn}

Theorem \ref{burnside-fr} states that a commutative fusion ring with h-integral dual is Burnside. Non-Burnside simple integral fusion rings can be found in \S \ref{sec:exacrit}.  It is easy to see that a hypergroup is Burnside if and only if the set of invertible standard elements coincides to the set of grouplike elements, see Proposition \ref{grouplike:set}.
\subsection{The dual hypergroup of a normalized hypergroup} \label{dual:nzd}
Let $\hb$ be an abelian  normalized  hypergroup. By the abelian assumption, $\wdb$ forms a basis for the dual $H^*$, which then is $\comp[\wdb]$, and on which we can define a multiplication. On the basis $\wdb$, the multiplication $\mui\star \muj$ is defined by declaring
\beq\label{mwa}
[\mui\star \muj](a):=\mui(a)\muj(a), \;\text{for all}\;a\in \mtc B
\eeq
and then extending linearly $\mui\star \muj$ on the entire $\cofa$. In this way, we obtain an algebra structure on $H^*$. 

It follows that there are some non-zero scalars ${\widehat p}_k(i,j)\in \comp$ such that 
\beq\label{hatp:eq}
\mui \star \muj=\sumktom{\wdht p}_k(i,j)\muk.
\eeq 
 By \cite[Proposition 2.3]{b-blms}, for any abelian normalized hypergroup $(H, \mtcb)$, the pair $(H^*, \wdb)$ is also an abelian normalized hypergroup. Moreover, the involution on $\wdb$ is given by $\muj\sent \mu_{j^\#}$, where $\mu_{j^\#} \in \wdb$ is defined by $\mu_{j^\#}(a):=\muj(a^*)$ for all $a\in \mtc B$. We denote by $\widehat{\hb}:=(H^*, \wdb)$ the dual hypergroup of the abelian normalized hypergroup $\hb$. It is also well known that for a normalized hypergroup $\hb$, then $\muj(a^*)=\overline{\muj(a)}$ for all $a\in \mtc B$. We also use the notation $\wh_j:={\wdht p}_1(j,j^\#)^{-1}$, called the \emph{order} of $\muj$ in the dual hypergroup $\whb$, see Lemma \ref{dual:order}.
 
For any RN-hypergroup $\hb$, we denote by $\hb_{ad}$ its adjoint sub-hypergroup, see \S \ref{adjoint:hyp}. As for fusion rings, this is defined as the sub-hypergroup of $\hb$ generated by the element $I(1):=\sumitom h_ix_ix_{i^*}$, where $h_i:=1/N^1_{i,i^*}$.
 
For any RN-hypergroup $\hb$, the central primitive idempotent corresponding to the character $\fp:H\ra\comp$ is denoted by $\lam_H$. For $H=K(\cc)$, the Grothendieck ring of a fusion category $\cc$, this primitive idempotent corresponds to the regular element.
\bn{defn}\label{def:hbz}
 A hypergroup $\hb$ is called \emph{normalizable} if we obtain a normalized hypergroup after rescaling it. (see \S \ref{nzble} for more details.)
\end{defn}
Based on Lemma \ref{nbz:cond} we denote such a hypergroup by $\hbz$ where $\mu_{1}\in \wdb$.

Inspired by the notion of dualizable probability groups from \cite{hdk}, we propose the following:
\bn{defn}\label{def:dualzb}
An ARN-hypergroup is called \emph{dualizable} if its dual is also ARN.
\end{defn}
\subsection{Main results of the paper}
\bt\label{hwrn}
Let $(H, B, \mu_1)$ be an abelian normalizable hypergroup such that $G\hb$ is a finite group. Then $(H, \ccb)$ is Burnside if and only if 
$$
(\prod_{j=1}^m\mu_j)^2=\frac{|G\hb|}{n\hbz}\big(\sum_{j\in \mtc I_{G\hb}}\wh_j\mu_j\big)
$$ 
\et
See Definition \ref{def:order} for $n(\_)$). We characterize Burnside dualizable ARN-hypergroups as follows:
\bt\label{hbz:hyp}
Let $\hbz$ be an abelian normalizable hypergroup such that $\whbz$ is RN. Then $(H, \ccb)$ is Burnside if and only if the following identity holds on $\whb$:
\beq\label{muj:gen:hyp}
\big(\prod_{j=1}^m\muj\big)^2=\frac{1}{n(\widehat{(H,{\mtc B})}_{ad})}\big(\sum_{\muj \in {\wdb}_{ad}}\wh_j\mu_j\big).
\eeq 
\et
Note that the RHS of the equation above corresponds to the integral $\lam_{ _{\whb_{ad}}}$ of the adjoint hypergroup $ \whb_{ad}$.
\br
If the abelian normalizable hypergroup $ \hb $ is either RN or dual RN, then according to Lemmas \ref{gh} and \ref{dual:rn}, the group $ G\hb $ must be finite. As a result, we can apply Theorem \ref{hwrn}. Specifically, by comparing this with Theorem \ref{hbz:hyp}, we can infer that in the dual RN case, $ n(\widehat{(H,{\mtc B})}_{ad})|G\hb| = n\hb $. By duality, we also have $ n((H,{\mtc B})_{ad})|G\whb| = n\hb $ in the RN case. These equalities can also be obtained from Harrison's Equations \eqref{quot:gh} and \eqref{w:qout:gh}, although in the more restrictive dualizable ARN case.
\er
%

For a fusion category $\cc$, its Grothendieck ring is denoted by $K(\cc)$, and  the set of isomorphism class representatives of simple objects of $\cc$ by $\irr(\cc):=\{X_1,\dots,X_m\}$. Let $d_i:=\fp(X_i)$  be the Frobenius-Perron dimension of $X_{i}$, $i\in {\mtc I}  = \{1,\dots, m\}$, and  $x_i:=[X_i]$ be the class of the simple object $X_i$ in the Grothendieck ring $K(\cc)$.

For any weakly-integral premodular category $\mathcal{C}$, $\wkc$ is both RN and rational, as demonstrated in \cite[Theorem 1.2]{b-blms}. Furthermore, for any unitary fusion category $\mathcal{C}$, $\wkc$ is also RN, see \cite{lpw,eno-nec}, and \emph{primary $3$-criterion} in \cite{hlpw}.
A normalizable hypergroup $\hbz$ is called \emph{dual-Burnside} if its dual $\whbz$ is Burnside (see examples from group theory in \S \ref{sec:exacrit}). The dual version of Theorem \ref{hwrn}, is the following:
\bt\label{hwrndual}
Let $(H, B, \mu_1)$ be an abelian normalizable hypergroup such that $G\whbz$ is a finite group. Then $(H, \ccb)$ is dual-Burnside if and only if 
$$
(\prod_{i=1}^m\nxi)^2=\frac{|G\whbz|}{n\hbz}\big(\sum_{i\in \mtc I_{G\whbz}}h_id_{i^*}x_i\big)
$$ 
\et
On the other hand, the dual version of Theorem \ref{hbz:hyp} is:
\bt\label{dual:hbz:hyp}
The following are equivalent for an ARN-hypergroup $\hb$.
\bne
\item 
The hypergroup $\hb$ is dual-Burnside.
\item
The following identity holds:
\beq\label{xi:gen:hyp}
\big(\prod_{i=1}^m\nxi\big)^2=\frac{1}{n(\hbad)}\big(\sum_{x_i\in \mtc B_{ad}} h_id_{i^*}x_i\big)
\eeq 
\ene 
\et

The dual-Burnside fusion categories can be characterized by:
\bt\label{dual:v:prop}
Let $\cc$ be a fusion category with a commutative Grothendieck ring.  Then $\kc$ is dual-Burnside if and only if the following holds:
\beq\label{dual:muj:gen:squared}
\big(\prod_{i=1}^m\nxi\big)^2=\frac{1}{\fp(\ccad)}\big(\sum_{x_i\in \ccad}d_ix_i\big).
\eeq
\et

About the dual-Burnside hypergroups, we prove the following:
\bt\label{db:integral:hypgs}
Let $(H,{\mtc B})$ be an abelian rational RN-hypergroup. If $(H,{\mtc B})$ is dual-Burnside then $\fp(H,{\mtc B}) \in \mathbb Q$.
\et
As a corollary, any commutative dual-Burnside fusion ring is weakly-integral.
Similar to the fusion ring settings, \cite{NG}, the universal grading group can be defined for any ARN-hypergroup. Using this grading, we can similarly define the concept of nilpotent ARN-hypergroup, and prove the following:
\bt\label{iff:nilpotent}
Let $\hb$ be a dualizable ARN-hypergroup. Then $\hb$ is nilpotent if and only if its dual $\whb$  is also nilpotent. Moreover, in this case they have the same nilpotency class.
\et
For a nilpotent ARN-hypergroup, we prove the following:
\bt\label{nilpotent:burnside}
A nilpotent dualizable ARN-hypergroup is both Burnside and dual-Burnside.
\et
Let $\cc$ be a modular fusion category. Then Corollary \ref{mtc:case} states that $\kc$ is Burnside if and only if it is dual-Burnside. The proof involves the Drinfeld map of $\cc$. Then Theorem \ref{dual:v:prop} implies that:
\bt\label{mtc:even:order}
In any modular fusion category $\cc$ with  $\kc$  Burnside:
\beq\label{square:star:ad}
\big(\prod_{i=1}^m\frac{x_i}{d_i}\big)^2=\frac{1}{\fp(\ccad)}\big(\sum_{x_i\in \irr(\ccad)}d_ix_i\big)
\eeq
\et
Note that the above theorem holds for weakly-integral modular categories since their Grothendieck rings are Burnside by \cite[Appendix]{gnn} or \cite[Theorem 2]{b-galois}.
\bc\label{mtc:odd:order}
In any modular fusion category $\cc$ with  $\kc$  Burnside and with the group of invertible objects $G(\cc)$ of odd order,
\beq\label{star:ad}
\prod_{i=1}^m\frac{x_i}{d_i}
=\frac{1}{\dim(\ccad)}\big(\sum_{x_i\in \irr(\ccad)}d_ix_i\big)\eeq
\ec
For any non-negative integer $n$, let $\mtc V(n)$ be the set of prime divisors of $n$. For a weakly-integral fusion category we also let $\mtc V(\cc):=\mtc V(\fpcc)$.
From Theorem \ref{dual:v:prop}, we can also derive the following result:
\bt\label{dual:burns:div}
Let $\cc$ be a fusion category such that $\kc$ is commutative and dual-Burnside. Then
\beq\label{ccad:di-sq}
\frac{(\proditom d_i)^2}{\fp(\ccad)}\in \mathbb Z.
\eeq
Moreover, if additionally $\cc$ is nilpotent  then  
\beq\label{nilp:ccad}
\mtc V(\ccad)=\bigcup_{i=1}^m \mtc V(d_i^{2}).
\eeq
\et
Note that since any nilpotent fusion category $\cc$ is weakly-integral then  $d_i^2\in \mathbb Z$ for all $i$, see \cite[Theorem 3.10]{NG}. Note that in the case of an integral fusion category Equation \eqref{nilp:ccad} can be written as 
\beq\label{nilp:ccad:int}
\mtc V(\ccad)=\bigcup_{i=1}^m \mtc V(d_i).
\eeq

Using the above identities, we prove the following results about the $\mathrm{FPdim}$ of simple objects:
\bt\label{first:div}
For any modular fusion category $\cc$ such that $\kc$ is dual-Burnside,
\beq\label{div:fpcc}
\frac{(\proditom d_i)^2}{\fp(\ccad)}\in \mathbb Z.
\eeq
For any weakly-integral modular fusion category $\cc$:
\beq\label{prime:set}
\mtc{V}(\cc)=\mtc V(\ccpt)\cup\big(\bigcup_{i=1}^m \mtc V(d_i^{2})\big)
\eeq
\et
Again, recall $d_i^2\in \mathbb Z$ for all $i$  by \cite[Theorem 3.10]{NG}. Moreover, if $\cc$ is an integral modular fusion category then the above equation becomes:
\beq\label{prime:set:int}
\mtc{V}(\cc)=\mtc V(\ccpt)\cup\big(\bigcup_{i=1}^m \mtc V(d_i)\big)
\eeq
In particular,
\br[Cauchy-type theorem] 
For every perfect integral modular fusion category,
\beq\label{prime:set:int:perfect}
\mtc{V}(\cc)=\bigcup_{i=1}^m \mtc V(d_i),
\eeq
thus, for all prime $p$ dividing $\fp(\cc)$, there is a simple object $X$ in $\cc$ such that $p$ divides $\fp(X)$; which can be interpreted as a Cauchy-type theorem. Consequently, $p^2$ divides $\fp(\cc)$ since $\fp(X)^2$ divides $\fp(\cc)$ by \cite[Proposition 8.14.6]{EGNO15}.
\er
Concrete examples where these results apply can be found in \S \ref{mtc:app}.
{Next theorem improves \cite[Theorem 4.5]{o-yu}.
\bt\label{gen:o-yu}
Let $\cc$ be a weakly-integral modular fusion category of $\fp = dm$ with $d$ square-free coprime with $m$ and $d_i^2$ for all $i$. Then $d\mid\fp(\ccpt)$ and $\cc$ admits a decomposition
$$
\cc\simeq \cd\boxtimes \cd'
$$
 where $\cd$ is a pointed modular fusion category of dimension $d$ and $\cd'$ a modular fusion category of dimension $m$.
\et
}
\bc\label{ccpt}
Let $\cc$ be any integral modular fusion category  of $\fp =dm$ with $d$ square-free coprime with $m$. 
{Then $d\mid\fp(\ccpt)$ and $\cc$ admits a decomposition
$$
\cc\simeq \cd\boxtimes \cd'
$$
 where $\cd$ is a pointed modular fusion category of dimension $d$ and $\cd'$ a modular fusion category of dimension $m$.
 }
\ec

In particular, the $\fp$ of a perfect integral modular fusion category does not have a powerless prime factor (Corollary \ref{p:square}). In particular, a perfect modular fusion category of even $\fp$ has its $\fp$ divisible by $4$.

\br
By the classification of non-pointed integral modular fusion categories of small rank in \cite{abpp} and by Corollary \ref{ccpt}, the $\fp$ of a non-pointed integral modular fusion category does not have a powerless prime factor for every rank less than $11 \times 3 = 33$, and every rank less than $17 \times 7 = 119$ in the odd-dimensional case.
\er
%

Integral modular fusion categories were recently intensively studied in the literature. The following conjecture is formulated in \cite{j-plav-odd}:
\begin{conjecture}\label{odd-j81}
There is no non-trivial perfect modular fusion category of odd $\fp$.
\end{conjecture}

This conjecture was checked for ranks less than $25$ in \cite{abpp} or \cite{CGP2023}. Note that the above result also applies to non-pointed simple integral modular fusion categories (since they are in particular perfect). The importance of the existence of such categories was described in \cite[\S 5]{lpr2}, where the following two \textbf{open} statements are proved to be equivalent, whereas Proposition \ref{main:link} states that they imply Conjecture \ref{odd-j81}.
\begin{statement}\label{5.4} Every simple integral fusion category is weakly group-theoretical.
\end{statement}
\begin{statement} \label{5.5} Every simple integral modular fusion category is pointed.
\end{statement}

We finally prove the following result:

\bt \label{thm:perfmodburn}
A perfect modular fusion category is (dual-)Burnside if and only if it is integral.
\et

Note that the integral modular fusion category $\mathcal{Z}(\rep(G))$ is perfect if and only if the finite group $G$ is perfect and centerless (more generally, see \S \ref{sub:perfmod}).

Shortly, the organization of this paper is the following. In \S \ref{afr} we recall the basics on hypergroups and abelian rational normalizable  hypergroups.  In \S \ref{gal} we develop few basic facts about Galois symmetries of rational hypergroups.  In \S \ref{v} hypergroups with Burnside property are studied. In the same section, we prove Theorem \ref{burnside-fr} and its more general version Theorem \ref{burnside} for hypergroups. In this section, Burnside's property for dual hypergroups is also studied. In \S \ref{kernels} we extend the notion of kernels of objects from the settings of fusion categories to the settings of ARN-hypergroups. In \S \ref{sec:adjsub} and \S \ref{universe}, we construct the universal grading group for ARN-hypergroups. A connection between this universal grading group and the group of grouplike elements of the dual is proven in Theorem \ref{grpk:eq:hbad}. Theorems \ref{hbz:hyp} and \ref{dual:hbz:hyp} are proven. In \S \ref{lucs} we define nilpotent ARN-hypergroup and prove Theorems \ref{iff:nilpotent} and \ref{nilpotent:burnside}. In \S \ref{afc} we prove Theorem \ref{dual:v:prop}. In \S \ref{pmc}, the applications to modular fusion categories are given (e.g. the proofs from Theorem \ref{mtc:even:order} to Corollary \ref{ccpt}), and we prove Proposition \ref{main:link} showing that any of the two statements from \cite[\S 5]{lpr2} implies Conjecture \ref{odd-j81}.

Next, \S \ref{sec:AppExtra} delves into applications and supplementary results: \S \ref{sub:perfmod} offers a characterization for a Drinfeld center to be perfect; \S \ref{sub:BurnInt} establishes some results concerning integrality and Burnside properties, proving Theorems \ref{db:integral:hypgs} and \ref{thm:perfmodburn}, particularly highlighting their equivalence in the perfect modular case; \S \ref{sub:nearmod} characterizes the near-group fusion categories which admit a modular structure.

Finally, \S \ref{sec:exacrit} presents concrete examples and counterexamples: \S \ref{sub:burn} introduces non-Burnside simple integral fusion rings; \S \ref{sub:dualburn} classifies certain dual-Burnside $\Rep(G)$; \S \ref{mtc:app} provides certain families of fusion rings lacking modular categorification by Theorems \ref{first:div} or \ref{ccpt}.


Throughout of this paper, all vector spaces and linear categories  are considered over the ground field  $\comp$ of complex numbers.

\subsection*{Acknowledgements} We thank Vicen\c tiu Pa\c sol and Andrew Schopieray for fruitful discussions on Galois groups of fusion rings and hypergroups. 
\section{Abelian rational normalizable hypergroups}\label{afr}
Let $\hb$ be a hypergroup as in Definition \ref{hyp:defn}.
A \emph{morphism of hypergroups} $\psi: (H, \mtcb)\ra (H', \mtcb')$  is an algebra morphism $\psi:H \ra H'$ such that $\psi(\mtcb)=\mtcb'$ and $\psi(b^{*})=\psi(b)^{*}$ for all $b\in \mtc B$. The morphism $\psi$ is called an \emph{isomorphism of hypergroups} if it is bijective.

Let $\hb$ be a normalized abelian hypergroup and $\whb$ its dual as defined in the introduction. There is a formula for the dual fusion coefficients ${\wdht p}_k(i,j)$ given in \cite[Proposition 2.1]{b-blms} by:
\beq\label{hat:pk}
\widehat{p}_k(i,j)=\frac{1}{n_k}\bigg(\sum_{a\in \mtcb}h_a\mu_i(a)\mu_j( a)\mu_k(a^*)\bigg)
\eeq
where the scalars $n_k$ are defined below (\ref{tau-form}), and $h_a:=1/N^1_{a,a^*}$ (then $h_1=1$). 
\br
Let $\hb$ be an abelian normalized hypergroup. We usually let $\mtcb:=\{x_{i}\}_{i \in \mtc I}$, with $\mtc I:=\{1, \dots m\}$ and $x_1 = 1$, and say that $\hb$ has rank $m$. We can simplify the notation $N_{x_i x_j}^{x_k}$ by $N_{i,j}^k$, and $h_{x_i}$ by $h_i$.
\er
\br \label{idmpt:hyp}
Let $\{x_i^{\circ}\}_{i \in \mtc I}$ 
be the linear dual basis of the linear basis $\mtcb$ of $H$. Therefore $x_i^{\circ} (x_{j}) =\delta_{i,j}$, for all $i,j\in \mtc I$. It is easy to see that in this case $\widehat F_i:=x_i^\circ\in H^*$ is a primitive central idempotent of $H^*$ corresponding to the character
$$\begin{array}{cccc}
             \widetilde{\omega}_i=\ev_{x_i}: & H^* & \to & \comp \\ 
                             & f & \mapsto &  f(x_i).
\end{array}$$
\er
\subsection{Function $\tau$ and associative non-degenerate bilinear form} \label{sub:FunTau}
Let $(H, \mtc B)$ be any abelian hypergroup. The finite dimensional $*$-algebra $H$ (so semismple) is commutative, therefore  $H\simeq \bigoplus_{i=1}^m \comp F_{i}$, where $(F_i)$ are the central primitive idempotents of $H$. We can define $\tau\in H^*$ with $\tau(x_i)=\delta_{i,1}$, where as above, $x_1=1$. Note that $\tau=x_1^\circ$ with the above notations. Moreover, in this case
$$\begin{array}{cccc}
             (\;,\;)_\tau : & H\times H & \to & \comp \\ 
                            & (a,b) & \mapsto & \tau(ab) 
\end{array}$$
is an associative symmetric non-degenerate bilinear form. Thus, $\tau(F_i)$ must be nonzero for all $i$. The nonzero scalars $n_i:=\tau(F_i)^{-1}$ are called the \emph{formal codegrees} of $(H, \mathcal B)$, see \cite{O3}. We can write
\beq\label{tau-form}
\tau =  \sum_{i=1}^m \tau(F_i) \mu_i = \sum_{i=1}^m\frac{1}{n_i}\mu_i.
\eeq
Note that
$\sumitom \frac{1}{n_i}=\tau(1)=1$.

From the definition of a hypergroup it follows that $\{h_ix_{i}\}_{i \in \mtc I}$ and $\{x_{i^{*}}\}_{i \in \mtc I}$
form a pair of dual bases for $(\;,\;)_\tau$.
Note that in this case $\{n_jF_j\}_{j \in \mtc I}$ and $\{F_j\}_{j \in \mtc I}$ form also a pair of dual basis for the bilinear form $(\;,\;)_\tau$.   The existence of the two pairs of dual bases implies that:
\beq\label{db}
\sumitom \nif x_i\ot  x_i^*=\sumjtom n_jF_j\ot F_j
\eeq
By applying $\id \ot \muj$, for the character $\muj$ corresponding to $F_j$, it follows that
\beq\label{fj}
F_j=\frac{1}{n_{j}}\sumitom \nif\mu_j(x_{i^*})x_{i}\eeq 
Applying $\muj\ot \muk$ to the above, we obtain the first orthogonality relation
\beq\label{first:gen:orth}
\sumitom \nif\mu_j(x_i)\mu_k(x_{i^*})=\delta_{j,k}n_j.
\eeq

Since $\muj(x_{i^*})=\overline{\muj(x_i)}$, for $j=k$, it follows from here that $n_j>0$. By a classical argument (see, for example, the proof of \cite[Theorem 2.4]{lpr2}), as for finite groups orthogonality, we obtain the second orthogonality relation:
\beq\label{second:orth}
\sumjtom\frac{1}{n_j}\muj(x_i)\muj(x_{l^*})=\delta_{i, l}h_i^{-1}.
\eeq
\subsection{Normalizable abelian hypergroups and their duals} \label{nzble}
Note that hypergroups can be rescaled as follows.  Let   $y_i:= \frac{x_i}{\ali}$, for some non-zero complex numbers $\ali $ with $\al_1=1$, $\al_{ _{i^*}}=\ovl{\ali}$. Let $\mtc B':=\{y_i\}_{i \in \mtc I}$.  Then, $(H, \mtc B')$ is also a hypergroup, with coefficients rescaled as $\frac{\alpha_k}{\alpha_i \alpha_j} N_{i,j}^k$; and then $h_i:=(N_{i,i^*}^1)^{-1}$ rescaled as $|\alpha_i|^2 h_i$. Consider $\tau'$ defined by $\tau'(y_i):= \delta_{i,1}$, then $\tau'(x_i) = \alpha_i  \delta_{i,1} = \delta_{i,1}$, because $\alpha_1=1$. Therefore $\tau'=\tau$ and  the formal codegrees $(n_j)$ are invariant by rescaling.

Recall the notion of a normalizable hypergroup from Definition \ref{def:hbz}. Note that any rescaled normalizable hypergroup is also normalizable.
\bl\label{nbz:cond}
An hypergroup $\hb$ is normalizable if and only if there is $\mu_1\in \wdb$ such that $\mu_1(x_i)$ nonzero, for all $i\in \mtc I$.
\el
\bpf
The rescaling $(H, \{\frac{x_i}{\ali}\})$ is normalized if and only if $\sum_{i \in \mtc I} \frac{\alpha_k}{\alpha_i \alpha_j} N_{i,j}^k = 1$, if and only if  $ \alpha_i \alpha_j = \sum_{i \in \mtc I} N_{i,j}^k \alpha_k$, if and only if $\mu_1: x_i \mapsto \alpha_i$ is an algebra morphism (i.e. an element of $\wdb$) with $\mu_1(x_i)$ nonzero, for all $i\in \mtc I$.
\epf
The normalizable hypergroup $\hb$ together with $\mu_1$ as in Lemma \ref{nbz:cond} is denoted $\hbz$. The normalized hypergroup $(H, \{\frac{x_i}{\mu_1(x_i)}\})$ is denoted $\ovl{\hbz}$.
\bn{defn} \label{dual:hypergroup}
Recall from \S \ref{dual:nzd} that any abelian normalized hypergroup $H$ admits a dual hypergroup $H^*$. In the case of $\ovl{\hbz}$, we denote this dual by ${\whbz}$.
\end{defn}

Let $\hbz$ be an abelian normalizable hypergroup.
From Equation \eqref{mwa}, it follows that the multiplication on the dual $\widehat{\hbz}$  can be written as 
\beq\label{mwa:nzd}
[f\star g]\left(\frac{x_i}{\mu_1(x_i)}\right)=f\left(\frac{x_i}{\mu_1(x_i)}\right)g\left(\frac{x_i}{\mu_1(x_i)}\right)
\eeq
for all $f,g \in H^*$. 
Then it is easy to see that  the algebra unit of $\whbz$ is $\mu_1$. Moreover, the involution on $\wdb$,  as  given  in \S \ref{dual:nzd}, becomes:
\beq \label{eq:dualinv}
\mu_{j^\#}(\frac{x_i}{\mu_1(x_i)})=\muj(\frac{x_{i^*}}{\ovl{\mu_{1}(x_i})}).
\eeq
%
%
%
%
\bn{defn} \label{def:order}
The \emph{order} of a normalizable hypergroup $\hbz$ is defined by
$$
n\hbz:=\sumitom h_i|\mu_1(x_i)|^2.
$$ 
\end{defn}
In the abelian case, the following holds by Equation \eqref{first:gen:orth}:
\beq\label{n0:nrzd}
n\hbz = n_1
\eeq 

\bp \label{resc0}
Let $\hbz$ be a normalizable hypergroup. Then the order $n\hbz$ is invariant under rescaling.
\ep
\bpf 
Let   $x'_i:= \frac{x_i}{\ali }$, for some non-zero complex numbers $\ali $ with $\al_1=1$, $\al_{ _{i^*}}=\ovl{\ali }$. Consider $(H, \mtc{B}', \mu_1)$ the rescaled hypergroup with $\mtc B':=\{x'_i\}_{i \in \mtc I}$. Recall that the order of $x'_i$ is $ h'_i=|\ali |^2h_i>0.$
Therefore
$$n(H, \mtc{B}', \mu_1)=\sumitom h'_i\mid\mu_1(x'_i)\mid ^2=\sumitom  h_i \mid\mu_1(x_i)\mid ^2=n\hbz. \qedhere $$ 
\epf

We shortly write $n(H)$ when $\mtc B$ and $\mu_1$ are implicitly understood. In particular, for a normalized hypergroup we may choose $\mu_1(x_i)=1$ for all $1\leq i \leq m$ and therefore $n(H)=\sumitom h_i$. For the rest of this paper we use the notation $d_i:=\mu_1(x_i)$ for any normalizable hypergroup $\hbz$.

We say that a normalizable hypergroup $\hbz$ is \emph{weakly-rational} if its order is a rational number, i.e. $n\hbz\in \mathbb Q$. Also, we say that $\hbz$ is \emph{weakly-integral} if its order is an integer, i.e. $n\hbz\in \mathbb Z$. Note that a fusion ring is weakly-rational if and only if it is weakly-integral since $\fp(R)$ is an algebraic integer in this case.
\br
Note that the formula from Equation \eqref{hat:pk} can be written for abelian normalizable hypergroups as follows (using the fact that $h_l$ rescales as $|d_l|^2h_l$):
 \beq\label{hat:pk:fus-alg}
\widehat{p}_k(i,j)=\frac{1}{n_k}\bigg(\sum_{x_{l}\in \mtcb}\frac{h_l}{d_l}\mu_i(x_l)\mu_j(x_l)\mu_k(x_{l^*})\bigg)
\eeq
\er 
\bl\label{dual:order} 
Let $\hbz$ be an abelian normalizable hypergroup. Then with the above notations, $n\hbz=n\whbz$. 
\el
\bpf
By applying (\ref{hat:pk:fus-alg}) to $i=j^\#$ and $k=1$, and then (\ref{first:gen:orth}), we obtain (see also \cite[Equation (2.10)]{b-blms}),
\beq\label{whj:eq}
\widehat{p}_1(i, i^\#)=\frac{n_i}{n\hbz}
\eeq
Let $\wh_i:=\frac{1}{\widehat{p}_1(i,i^\#)}$ be the \emph{order} of $\mu_i$ in $\whbz$. It follows that $n\whbz=\sumitom \wh_i=n\hbz(\sumitom \frac{1}{n_i})=n\hbz$.
\epf
The next proposition holds for any normalized hypergroup, not necessarily abelian. 
\bp
For any normalized hypergroup $\hb$,
\beq\label{f0}
F_1=\frac{1}{n\hb}\big(\sumitom h_{i^*}x_i\big)
\eeq
where $F_1$ is the idempotent corresponding to the augmentation map $\mu_1:H\ra \comp,\;a\sent 1$ for all $a\in \mtc B$.
\ep
%
\bpf
Suppose that
$F_1=\sumitom \ali x_i$ for some $\ali \in \comp$. Recall that $x_lF_i = \mu_i(x_l)F_i$. Then, $x_lF_1=F_1=\sumitom \ali  x_lx_i$. It follows that $\tau(F_1)=\al_1=\tau( x_lF_1)=\al_{l^*}h_{l}^{-1}.$ Thus $\al_{l^*} = h_l \al_1$, so $\al_{l} = h_{l^*} \al_1$,  $F_1=\sumitom \al_ix_i={\al_1}(\sumitom {h_{i^*}}x_i)$ and $1=\mu_1(F_1)={\al_1}(\sumitom {h_i})$. This implies that  $\al_1=\frac{1}{n(H)}$ and therefore $\al_i=\frac{h_{i^*}}{n(H)}$ for all $i\in \mtc I$. 
\epf

\bc\label{fz:gen}
In any normalizable hypergroup $\hbz$,
\beq\label{f0:gen}
F_1=\frac{1}{n\hbz}\big(\sumitom h_{i^*}d_{i^*}x_i\big)
\eeq
\ec

\br\label{idmpt:fus:alg}
Let $\hbz$ be an abelian (rational) normalizable hypergroup.  In order to find the central primitive idempotents of $H^*$, we can apply Remark \ref{idmpt:hyp} after normalization. If  $x_i^{\circ}\in H^*$ are defined as above by $\lag x_i^{\circ}, x_{i'}\rag =\delta_{i,i'}$,
then  $\wtf_i:=\di x_i^\circ\in H^*$ are the primitive central idempotents of $H^*$ corresponding to the character
$$\begin{array}{cccc}
             \widetilde{\omega}_i=\ev_{\overline{x}_i}: & H^* & \to & \comp \\ 
                             & \muj & \mapsto &  \muj(\overline{x}_i),
\end{array}$$
where $\overline{x}_i = \frac{x_i}{d_i}$. By identifying $H^{**}$ with $H$, we get that $\ev_{\overline{x}_i} = \overline{x}_i$. 
\er
Recall that $\overline{x}_1 = x_1$.

\bp\label{dual:ni}
Let $\hbz$ be an abelian normalizable hypergroup. Then $\whbz$ is an abelian normalized hypergroup whose formal codegrees are
\beq\label{dual:codeg}
\widehta{n}_i=\frac{n\hbz}{h_i|d_i|^2}.
\eeq
\ep
\bpf 
As already mentioned, \cite[Proposition 2.3]{b-blms} shows that $\whbz$ is an abelian normalized hypergroup.

Let $\overline{x}_i = \frac{x_i}{d_i}$ be the normalization of $x_i$. By Remark \ref{idmpt:fus:alg}, $\overline{x}_i:H^*\ra\comp$ are the characters of $\whbz$. By definition of $\wtau\in (H^*)^*=H$, $\wtau(\mu_j)=\delta_{j,1}$ as $\mu_1$ is the unit of $H^*$.  Thus $\wtau=F_1$. Equation \eqref{tau-form} for $H^*$ gives that $\wtau=\sumitom \frac{1}{\widehta{n}_i}\overline{x}_i$. So, by Equation \eqref{f0:gen},
$$\frac{1}{n\hbz}\big(\sumitom h_{i^*}d_{i^*}x_i\big) = F_1 = \wtau = \sumitom \frac{1}{\widehta{n}_i}\frac{x_i}{d_i}.$$ The result follows by $h_{i^*} = (N_{i^*,i}^1)^{-1} = (N_{i,i^*}^1)^{-1} = h_i$, in the abelian case.
\epf

\bc\label{wfz}
Let $\hbz$ be an abelian normalizable hypergroup. Then the primitive central idempotent of the linear character $x_{1}\in \widehat{\widehat{\mtc B}}=\mtc B$ is given by
\beq
\widehat{F_{1}}=\sumjtom\frac{\mu_j}{n_{j}}\in H^{*}.
\eeq
\ec
\bpf
Proposition \ref{dual:ni} shows in its proof  that $\wtau=F_1$. By duality, since  $\widehat{\widehat{\hbz}}=\overline{\hbz}$, it follows that $\tau=\what{F}_1$. Equation \eqref{tau-form} finishes the proof.
\epf

A  \emph{morphism of normalizable hypergroups }$\psi: (H, \mtcb, \mu_1)\ra (H', \mtcb', \mu_1')$  is  a  morphism of hypergroups with the property that  $\mu'_{1}\circ \psi=\mu_1$.
\subsection{Frobenius-Perron theory for a RN-hypergroups}\label{fp:hy} 
Frobenius-Perron theory can be defined for RN-hypergroups, in the same manner as for fusion rings, see \cite[\S 3]{EGNO15}.  We denote by $\fp(x_i)$ the Frobenius-Perron eigenvalue of the left multiplication operator by $x_i$ on $H$, i.e. the ($\ell^2$) matrix norm of $N_i = (N_{i,j}^k)_{j,k}$.  Recall that a hypergroup is called \emph{symmetric} when $N_{a,b}^1 = N_{b,a}^1$, for all $a,b \in \mathcal{B}$, but $N_{a,b}^1 = \delta_{a^*,b} h_a^{-1}$, so symmetric means that $h_a = h_{a^*}$, for all $a \in \mathcal{B}$

\begin{lemma} \label{lem:absym}
An abelian hypergroup is symmetric.
\end{lemma}
\begin{proof}
By definition and abelian assumption, $h_i^{-1} := N_{i,i^*}^1 = N_{i^*,i}^1 = h_{i^*}^{-1}$.
\end{proof}

\begin{lemma} \label{lem:hihjhk}
If $N_{i,j}^k$ is nonzero then 
\beq \label{eq:hihjhk} \frac{h_{i^*}h_{j^*}h_{k}}{h_{i}h_{j}h_{k^*}} = 1. \eeq
\end{lemma}
\begin{proof}
The associativity of $H$ reformulates as $\sum_s N_{i,j}^s N_{s,k}^t = \sum_s N_{j,k}^s N_{i,s}^t$, for all $i,j,k,t$. If $t=1$, we get that $\sum_s N_{i,j}^s N_{s,k}^1 = \sum_s N_{j,k}^s N_{i,s}^1$. But recall that $N_{a,b}^1 = \delta_{a^*,b} h_a^{-1}$. It follows that $N_{i,j}^{k^*} h_{k^*}^{-1} = N_{j,k}^{i^*} h_{i}^{-1}$, in other words, $N_{i,j}^{k}  = \frac{h_{k}}{h_{i}} N_{j,k^{*}}^{i^*}$.
By applying this last equality, we get that $N_{j,k^{*}}^{i^*} = \frac{h_{i^*}}{h_{j}} N_{k^{*},i}^{j^*}$, and $N_{k^{*},i}^{j^*} =  \frac{h_{j^*}}{h_{k^*}}N_{i,j}^{k}$. Thus 
\beq \label{eq:preFR}
N_{i,j}^{k} = \frac{h_{k}}{h_{i}} N_{j,k^{*}}^{i^*} = \frac{h_{i^*} h_{k}}{h_{i}h_{j}} N_{k^{*},i}^{j^*} =   \frac{h_{i^*}h_{j^*}h_{k}}{h_{i}h_{j}h_{k^*}} N_{i,j}^{k}
\eeq
The result follows.
\end{proof}

\bt  \label{thm:Sym}
A (finite-dimensional) hypergroup is symmetric, i.e. $h_{i^*} = h_i$, for all $i \in \mathcal{I}$. In particular, the function $\tau$ from \S \ref{sub:FunTau} is symmetric.
\et
\bpf Introduce the notation \(r_i := \frac{h_{i^*}}{h_i}\). Lemma \ref{lem:hihjhk} reformulates as: $$N_{i,j}^k \neq 0 \Rightarrow  r_i r_j = r_k.$$ We aim to demonstrate that \(r_i\) is a root of unity for every \(i \in \mathcal{I}\). First, select \(k_2\) such that \(N_{i,i}^{k_2}\) is nonzero, which implies \(r_i^2 = r_{k_2}\). Next, choose \(k_3\) such that \(N_{i,k_2}^{k_3}\) is nonzero, leading to \(r_i r_{k_2} = r_{k_3}\) and thus \(r_i^3 = r_{k_3}\). Continuing this process, we find that the sequence \(\{r_i^n\}\) for \(n \in \mathbb{N}\) matches the set \(\{r_j\}\) for \(j \in \mathcal{J}\), with \(\mathcal{J}\) being a subset of \(\mathcal{I}\). Due to the finite dimensionality of the hypergroup, \(\mathcal{I}\) is finite, which implies that \(\mathcal{J}\) is also finite. Consequently, there exist integers \(m\) and \(n\) with \(m > n\) such that \(r_i^m = r_i^n\). Thus, \(r_i^{s} = 1\) where \(s = m-n > 0\), confirming that \(r_i\) is a root of unity. Additionally, the axioms of the hypergroup state that \(N_{a,a^*}^1\) is positive for all \(a \in \mathcal{B}\); hence \(h_i\) is positive for all \(i \in \mathcal{I}\). Therefore, \(r_i = \frac{h_{i^*}}{h_i}\) is a positive root of unity, which must be \(1\).
\epf

In the rest of the paper, we assume that the involution on $\mathcal{B}$ extends into an anti-involution on $H$, i.e. $N_{i,j}^k = N_{j^*,i^*}^{k^*}$, for all $i,j,k$.

\bp[Frobenius Reciprocity]
The following equalities hold:
\beq \label{eq:FR}
h_k^{-1}N_{i,j}^{k} = h_j^{-1} N_{i^*,k}^{j} = h_i^{-1} N_{j,k^*}^{i^*}  = h_k^{-1}N_{j^*,i^*}^{k^*} = h_j^{-1} N_{k^*,i}^{j^*} = h_i^{-1} N_{k,j^*}^{i}.
\eeq
\ep
\bpf
By (\ref{eq:preFR}), $N_{i,j}^{k} = \frac{h_{i^*} h_{k}}{h_{i}h_{j}}N_{k^*,i}^{j^*}$, then by anti-involution, $N_{k^*,i}^{j^*} = N_{i^*,k}^{j}$, and by (\ref{eq:hihjhk}), $\frac{h_{i^*} h_{k}}{h_{i}h_{j}} = \frac{h_{k^*}}{h_{j^*}}$. But $h_{k^*} = h_k$ and $h_{j^*} = h_j$ by Theorem \ref{thm:Sym}. The first equality follows. The rest is similar.
\epf

\bc \label{cor:di=di*}
For every $x_i \in \mathcal{B}$ then $\fp(x_i) = \fp(x_{i^*})$.
\ec
\bpf
Let $N_i$ be the matrix $(N_{i,j}^k)_{j,k}$. By the first equality of (\ref{eq:FR}), $N_{i,j}^{k} = \frac{h_{k}}{h_{j}} N_{i^*,k}^{j} $, so $N_i = D N_{i^*}^T D^{-1}$, with $D = {\rm diag}(h_{i})$ and $()^T$ the matrix transpose. Thus $$ \fp(x_{i^*})= \Vert N_{i^*} \Vert = \Vert N_{i^*}^T \Vert =  \Vert D^{-1} N_{i} D  \Vert =  \Vert N_i \Vert =  \fp(x_{i}). $$
The second last equality hold because $D^{-1} N_{i} D$ have the same eigenvalues than $N_i$, because if $N_iv=\lambda v$ then $D^{-1} N_{i} Dw = \lambda w$, with $w = D^{-1}v$.
\epf

%


Since we are dealing with the finite-dimensional case, by Theorem \ref{thm:Sym}, every hypergroup is symmetric. Therefore, we can omit this assumption in the rest of the paper, as it is automatically satisfied.
\br\label{unique:fp}
Note that any hypergroup $\hb$ is \emph{transitive} in the sense of \cite[Definition 3.3.1]{EGNO15}, i.e. $\forall i,j \in \mathcal{I}$ there are $k_1, k_2 \in \mathcal{I}$ such that $N_{i,k_1}^j$ and $N_{k_2,i}^{j}$ are nonzero.  Indeed, take $k_1, k_2$ such that $N_{i^*,j}^{k_1}$ and $N_{j,i^*}^{k_2}$ are nonzero, the result follows by Frobenius reciprocity (\ref{eq:FR}).
%

Observe that an analogue of  \cite[Proposition 3.3.6]{EGNO15}  holds for \invisible{symmetric}RN-hypergroups as it requires Frobenius-Perron \cite[Theorem 3.2.1]{EGNO15} involving a RN-matrix (i.e. with real nonnegative entries). In particular,  $\fp:H\ra \comp$ is the unique algebra morphism which takes positive values on $\mtc B$.
\er
\bn{defn}
The Frobenius-Perron dimension of a RN-hypergroup $\hb$ is
$$
\fp\hb:=\sumitom h_i\fp(x_i)^2.
$$ 
\end{defn}

By Lemma \ref{nbz:cond} and Remark \ref{unique:fp}, any RN-hypergroup is normalizable via the linear character $\mu_1=\fp$. By Definition \ref{def:order} and positivity of $\fp$,
$$\fp\hb=n\hbfr.$$
Moreover, in the normalized case, by uniqueness in Remark \ref{unique:fp}, $\fp$ must be the augmentation map, i.e. $\fp(x_i) = 1$ for all $i \in \mathcal{I}$.


Recall from Proposition \ref{resc0} that for a normalizable hypergroup $\hbz$, the order $n\hbz$ is invariant under rescaling.

\bp \label{resc}
If $\hb$ is a RN-hypergroup, then the Frobenius-Perron dimension of $\hb$ is invariant under rescaling with real positive numbers.
\ep
\bpf 
If $\hb$ is a RN-hypergroup  and  $\ali =\al_{ _{i^*}}>0$, for all $i \in \mathcal{I}$, then $(H', \mtc B') = (H,\{\frac{x_i}{\alpha_i} \})$ is RN and 
\begin{eqnarray*}
\fp(H', \mtc B')&=&\sumitom h'_i \fp(\frac{x_i}{\alpha_i})^2=\\&=& \sumitom | \ali |^2h_i\frac{\fp(x_i)^2}{\ali^2 }\\&=&\fp\hb. \qedhere
\end{eqnarray*}
\epf
\subsection{Fourier transform}
Let $\hbz$ be an abelian normalizable hypergroup. Define
\beq\label{mtf-eqq}
\mtf :\hb \ra \whb, \; x_i\sent \frac{n\hbz}{h_{i^*}}x_{i^*}^\circ.
\eeq
It is clear that $ \mtf $ is a linear isomorphism. Additionally, as noted in Remark \ref{idmpt:fus:alg}, we have $ \wtf_i := \di x_i^\circ $. Therefore:
\beq\label{mtf-eqq-2}
\mtf(x_i)=\frac{n\hbz}{d_{i^*} h_{i^*}}\wtf_{i^*}.
\eeq
For the sake of brevity, we define $ |H| := n \hbz $.
\bl\label{mtf-id}
Let $ \hbz $ be an abelian normalizable hypergroup. With the above notations, for all $x, y\in H$:
$$
\langle \mtf(y), x\rangle=|H|\tau(xy).
$$
\el
\bpf
It suffices to verify the above identity for $ (x,y) = (x_i, x_j) $, for all $ i, j$. We have:
$$
\langle\mtf(x_j), x_i\rangle \numeq{\ref{mtf-eqq}}\langle \frac{|H|}{h_{j^*}}x_{j^*}^\circ, x_i\rangle=\frac{|H|}{h_{i}}\delta_{i,j^*}=|H|\tau(x_ix_j).
$$
since $ h_i = h_{i^*} $ (see Lemma \ref{lem:absym}), and $ \delta_{i,j^*} h_i^{-1} = \tau(x_i x_j) $, because $ h_i^{-1} = N_{i,i^*}^1 $ and $ \tau = x_1^\circ $.
\epf
\bp
Using the above notations, we have that 
\beq\label{mtf-muk}
\mtf(F_k)=\frac{|H|}{n_k}\mu_k.
\eeq	
\ep
\bpf
By taking $ (x,y) = (F_j,F_k) $ in Lemma \ref{mtf-id}, we obtain:
$$
\langle\mtf(F_k), F_j\rangle=|H|\tau(F_jF_k)=\delta_{j,k}|H|\tau(F_j).
$$
From Equation \eqref{fj}, we have:
$$
\tau(F_j)=\frac{1}{n_j}\big(\sumitom h_i\muj(x_{i^*})\tau(x_i)\big)=\frac{1}{n_j}.
$$
Thus
$$
\langle \mtf(F_k), F_j \rangle=\delta_{j,k}\frac{|H|}{n_j},
$$
which shows that $\mtf(F_k)=\frac{|H|}{n_k}\mu_k$, since $\mu_i(F_j) = \delta_{i,j}$.
\epf

\subsection{Grouplike elements in hypergroups} \label{grp:lk}
In  the rest of this paper, all abelian RN-hypergroups are normalized by $\mu_1=\fp$, and therefore their duals are also considered with respect to $\mu_{1}=\fp$.  Recall that $d_i:=\mu_{1}(x_i) = \fp(x_i)$.
\bl\label{grouplike:nn} 
For any \invisible{symmetric}RN-hypergroup $\hb$, $h_id_i^2\geq 1$. Moreover, $h_id_i^2=1$ if and only $h_ix_ix_{i^*}=x_1.$ 
\el
\bpf
 The following equality holds
 \beq\label{adj}
 x_ix_{i^*}=\frac{1}{h_i}x_1+\sum_{k=2}^m N^k_{ii^*}x_k.
 \eeq
 
Passing to $\fp$, $h_id_i^2=1+h_i\big(\sum_{k=2}^m N^k_{ii^*}d_k\big)\geq 1.$ If $h_id_i^2=1$ then $\sum_{k=2}^m N^k_{ii^*}d_k=0$, and therefore $N^k_{ii^*}=0$ for all $k\neq 1$. Thus $h_ix_ix_{i^*}=x_1.$ Conversely, if $h_ix_ix_{i^*}=x_1$, applying $\fp$, it follows that $h_id_i^2=1$.
\epf
\bc\label{grplike:fr}
In any fusion ring $\hb$, $x_ix_{i^*}=x_1$ if and only if $d_i=1$.
\ec
Recall the set of grouplike elements $\ghb$ from Definition \ref{grplike:def}.
\br  \label{rk:grlk:id}
For a normalizable hypergroup $\hbz$, we establish that $x_i \in G\hb$ if and only if:
\begin{equation}\label{grlk:id}
\frac{x_i x_{i^*}}{d_id_{i^*}} = x_1.
\end{equation}
Applying $\mu_1$ to Equation (\ref{grplk:id:eq}), we find that $h_i d_i d_{i^*} = 1$, from which Equation \eqref{grlk:id} follows. The converse is similar.
\er
\bl\label{grp:id}
Suppose that $\hbz$ is an abelian normalizable hypergroup. Then $x_i\in G\hb$ if and only if 
\beq\label{xi:abs:value}
|\muj(\nxi)|=1\;\;\text{for all}\;\; j\in \mtc I.
\eeq
\el
\bpf
If $x_i\in G\hb$ then applying $\muj$ to Equation \eqref{grlk:id} results in $|\muj(\nxi)|^2=1$, since $\muj(x_{i^*})=\overline{\muj(x_i)}$.  Conversely, if $\muj(\nxi\frac{x_{i^*}}{d_{i^*}})=1$ for all $\muj$, then $\nxi\frac{x_{i^*}}{d_{i^*}}=x_1$, since $\hb$ is abelian, and therefore $x_i\in G\hb$ by Remark \ref{rk:grlk:id}.
\epf
Dually we have the following:
\bl\label{dual:grp:id}
Let $\hbz$ be an abelian normalizable hypergroup. Then $\muj\in G\whbz$ if and only if 
\beq\label{muj:abs:value}
|\muj(\nxi)|=1\;\; \text{for all} \;\; i\in \mtc I.
\eeq
\el
\bpf
The dual version of Equation \eqref{grlk:id} implies that $\muj\in G\whbz$ if and only if $\muj \mu_{j^\#}=\mu_1.$ Evaluating at $\nxi$ results in the desired identity.
\epf
\bc\label{frmlcdg:grplike}
In any abelian normalizable hypergroup $\hbz$, for any $\muj \in G\whbz$, it holds that $n_j=n(H, B, \mu_1)$.
\ec
\bpf
By Lemma \ref{dual:grp:id}, we find that $\mid\muj(x_i)|=|d_i|$ for all $i\in \mtc I$. 
From Equation \eqref{first:gen:orth}, we deduce that $n_j=\sumitom h_i|\muj(x_i)|^2=\sumitom h_i|d_i|^2=n(H, B, \mu_1)$.
\epf
\bl\label{l:fp:ineq}
Suppose that $\hbz$ is an abelian normalizable hypergroup such that one of the following holds: 
\bne
\item
$\hb$ is RN and $\mu_{1}=\fp$,
\item $\whbz$ is RN. 
\ene
Then 
\beq\label{fp:ineq}
|\muj(x_{i})|\leq |\mu_1(x_{i})| = |d_{i}|\;\;\text{, for all}\; i,j\in \mtc I.
\eeq
\el
\newcommand{\fpxi}{{\fp(x_{i})}}
\bpf
Suppose first that $\hb$ is RN and $\mu_{1}=\fp$. Thus $d_i=\fpxi>0$. But $x_{i}=\sumjtom \muj(x_{i})F_{j}$, and therefore $x_iF_j = \muj(x_{i})F_{j}$ meaning that $\muj(x_{i})$ is  an eigenvalue of $L_{x_{i}}$. It follows that $|\muj(x_i)|\leq d_i$, by definition of $\fp$.

Now, suppose that $\whbz$ is RN. It follows that $\fp(\muj)=1$,  since by Remark \ref{unique:fp}, in this case, $\muj\sent 1$ is the only algebra morphism on $\whbz$ which takes positive values on $\wdb$. Then, the equation
\beq\label{muj:wtei}
\muj=\sumitom \muj(\nxi)\wtf_{i}
\eeq
implies as above that $\big|\muj(\nxi)\big|\leq 1$.
\epf
\bl\label{max:formal:codeg}
Under the same assumptions as Lemma \ref{l:fp:ineq}, we have $n_j\leq n\hbz$, for all $j \in \mtc I$. Moreover, $n_j=n\hbz$ if and only if $\muj$ is a grouplike element.
\el
\bpf
Lemma \ref{l:fp:ineq} and Equation \eqref{first:gen:orth} give:
$$n_j=\sumitom h_i|\muj(x_i)|^2\leq \sumitom h_i|d_i|^2=n(H, B, \mu_1),
$$
and the equality holds if and only if $|\muj(x_i)|^2=|d_i|^2$, if and only if $\muj\in G\whbz$, by Lemma \ref{dual:grp:id}.
\epf

\bl \label{gh}
In any \invisible{symmetric}RN-hypergroup $\hb$, the set $\{\nxi\}$ of normalized grouplike-elements forms a finite group. Moreover, the inverse of $\nxi$ is given by $\frac{x_{i^*}}{d_{i^*}}$ for all $x_i\in \ghb$.
\el
\bpf
If $x_ix_{i^*}=\frac{1}{h_i}x_1$ and $x_jx_{j^*}=\frac{1}{h_j}x_1$ then
$
(x_ix_j)(x_ix_j)^*=\frac{1}{h_ih_j}x_1. $
Thus $(\sumktom N^k_{ij}x_k)(\sumktom N^k_{ij}x_{k^*})=\frac{1}{h_ih_j}x_1,$
so that there is a unique $k \in \mathcal{I}$ such that $N_{i,j}^k$ is nonzero (if it were not the case, then by the axioms of hypergroups coupled with the RN assumption, the left-hand side of the aforementioned equation would decompose into several components). Consequently, we have $(N^k_{ij}x_k)(N^k_{ij}x_{k^*})=\frac{1}{h_ih_j}x_1$, and furthermore, $x_i x_j = N^k_{ij}x_k$. However, applying $\fp$ yields $N^k_{ij} = (d_i d_j)/d_k$, so $\frac{x_i}{d_i} \frac{x_j}{d_j} = \frac{x_k}{d_k}$, and by the definition of a grouplike element, we have $h_i^{-1} = d_i^2$. Synthesizing all these equalities, we deduce that $((d_i d_j)/d_k)^2 x_k x_{k^*} = (d_i d_j)^2 x_1$, leading to the conclusion that $x_k x_{k^*} = d_k^2 x_1 = h_k^{-1}x_1$. There remain to prove that if $x_i \in \ghb$ then so is $x_{i^*}$. By Lemma \ref{grouplike:nn}, $x_{i^*}$ is grouplike if and only if $h_{i^*} d_{i^*}^2 = 1$. But $h_{i^*} = h_i$ by Theorem \ref{thm:Sym}, and $d_{i^*} = d_i$ by Corollary \ref{cor:di=di*}, so $h_{i^*} d_{i^*}^2 = h_{i} d_{i}^2 = 1$, because $x_i$ is grouplike. Finally, by above, there is $k$ such that $\frac{x_i}{d_i} \frac{x_{i^*}}{d_{i^*}} = \frac{x_k}{d_k}$, and by the hypergroup axioms, $k$ must be $1$.
\epf
We denote $\overline{G\hb}:=\{\nxi\;|\; x_i\in G\hb\}$ the above group of normalized grouplike elements. Clearly $\overline{G\hb}=G\overline{\hb}$, the grouplike elements of the normalized hypergroup $\overline{\hb}$.

By duality we obtain the following:
\bl\label{gh:dual}
Let $\hbz$ be an abelian normalizable hypergroup such that its dual $\whbz$ is RN. Then $\wghb$ is a group.
\el
Recall that $\overline{\mtc B}=\{\nxk\;|\; x_k\in \mtcb\}$.
\bl\label{perm} 
Let $\hb$ be a \invisible{symmetric}RN-hypergroup. Let $x_i \in \ghb$, then $\nxi\nxj, \nxj\nxi \in \overline{\mtc B}$, for any $x_j\in \mtc B$.
\el
\bpf
Suppose that 
$$
x_ix_j=\sum_{k\in \mtc A} N^k_{ij}x_k,
$$
where $\mtc A$ is  a set such that $N^k_{ij}>0$. Multiplying the above equality by $x_{i^*}$, we obtain
$$
\frac{1}{h_i}x_j=\sum_{k\in \mtc A} N^k_{ij}x_{i^*}x_k,
$$ 
since $h_{i^*}x_{i^*}x_i=x_1$ as $x_{i^*}$ is also grouplike by Lemma \ref{gh}, and $h_{i^*} = h_i$ by Theorem \ref{thm:Sym}.  By RN assumption, $x_{i^*}x_k=\al_{ik}x_j$, for some scalar $\al_{ik} > 0$ and all $k\in \mtc A$. Thus $\al_{ik}=\frac{d_id_k}{d_j}$. Multiplying the second last equality by $x_i$, we get that $h_i^{-1} x_k = \al_{ik} x_i x_j$, which reformulates as $\nxi\nxj = \nxk \in \overline{\mtc B}$. Idem for $\nxj\nxi$ using right multiplications.
\epf
\bl\label{const:scalars}
Suppose that $\hb$ is an abelian RN hypergroup and $\muj\in \widehat{\mtc B}$ a character such that $\muj(x_{m})=d_m\omega_{m}$ and $\muj(x_{n})=d_n\omega_{n}$ for some scalars $|\omega_{m}|=|\omega_{n}|=1$. Then $\muj(x_{p})=\omega_{m}\omega_{n}d_{p}$ for all constituents $x_{p}$ of $x_{m}x_{n}$.
\el
\bpf
Apply the absolute value triangle inequality and Lemma \ref{l:fp:ineq}.
\epf
\bl\label{dual:rn}
If $\hb$ is an abelian RN hypergroup then $\wqghb$ is a group.
\el
\bpf
Suppose that $\muj, \muk\in \wqghb$. By Lemma \ref{dual:grp:id}, $\muj(x_m)=\omega_md_m$ and $\muk(x_m)=\eta_md_m$, for some roots of unity $\omega_m$ and $\eta_{m}$. Now, 
$$(\muj\star \muk)(x_mx_n) = \sum_p N^p_{m,n}d_p(\muj\star \muk)(\nxp),$$
but by Equation (\ref{mwa:nzd}), $$ (\muj\star \muk)(\nxp) =  \muj(\nxp) \muk(\nxp),$$
and by Lemma \ref{const:scalars}, $\muj(x_p) = \omega_m \omega_n d_p$ and  $\muk(x_p) = \eta_m \eta_n d_p$, when $N^p_{m,n}$ is nonzero. Thus $$(\muj\star \muk)(x_mx_n) =  \sum_p N^p_{m,n}d_p \omega_m \omega_n  \eta_m \eta_n = \omega_m \omega_n  \eta_m \eta_n d_m d_n$$

On the other hand, 
$$(\muj\star \muk)(x_m) = d_m(\muj\star \muk)(\nxm) =  d_m \muj(\nxm) \muk(\nxm) = d_m \omega_m \eta_m.$$
Idem, $(\muj\star \muk)(x_n) = d_n \omega_n \eta_n$. So 
$$ (\muj\star \muk)(x_m) (\muj\star \muk)(x_n) = d_m \omega_m \eta_m d_n \omega_n \eta_n = (\muj\star \muk)(x_mx_n),$$
by above, which shows that $(\muj\star \muk)\in \wdb$. 

Moreover, $(\muj\star \muk)(\nxp)=\muj(\nxp)\muk(\nxp)=\omega_p\eta_p$, which by Lemma \ref{dual:grp:id}, implies that $\muj\star \muk\in G\whb$.
\epf
\bp\label{perm:dual}
Let $\hb$ be an abelian RN hypergroup and $\mu\in \wqghb$. Then, $\mu\star \muk\in \widehta{B}$, for any $\muk\in \widehat{B}$.
\ep
\bpf
Since $\mu\in \wqghb$ and $\hb$ is RN, by Lemma \ref{dual:grp:id}, $\mu(x_i)=\omega_id_i$ for some root of unity $\omega_i$ and for all $i\in \mtc I$. 
By Lemma \ref{const:scalars}, $\mu(x_l)=\omega_i\omega_jd_l$ for any constituent $x_l$ of $x_ix_j$, therefore
$$
(\mu\star\muk)(x_ix_j)=\sumltom N^l_{ij} (\mu\star\muk)(x_l)=\omega_i\omega_j \sumltom N^l_{ij} \muk(x_l)=\omega_i\omega_j \muk(x_ix_j).
$$
On the other hand 
$$
(\mu\star\muk)(x_i)(\mu\star \muk)(x_j)=\omega_i\muk(x_i)\omega_j\muk(x_j)=\omega_i\omega_j \muk(x_ix_j).
$$
Thus $(\mu\star\muk)(x_ix_j)=(\mu\star\muk)(x_i)(\mu\star\muk)(x_j)$ which shows that $\mu\star\muk\in \wdb$.
\epf
By duality, Lemma \ref{dual:rn} and Proposition \ref{perm:dual} reformulate as follows:
\bp\label{whb:rn}
If $\hbz$ is an abelian normalizable  hypergroup such that $\whbz$ is RN, then $\overline{G\hb}$ is a group. Moreover, the left multiplication $L_{\nxi}$, for any $x_i\in \ghb$, permutes the normalized basis $\overline{\mtc B}$.
\ep
Lemmas \ref{gh}, \ref{gh:dual}, \ref{perm}, \ref{dual:rn} and Propositions \ref{perm:dual},  \ref{whb:rn} imply the following:
\bc\label{all:in:one}
Suppose that $\hbz$ is an abelian normalizable hypergroup such that one of the following holds: 
\bne
\item
$\hb$ is RN and $\mu_{1}=\fp$ 
\item $\whbz$ is RN. 
\ene
Then, the following statements hold:
\bne
\item 
Both $\overline{G\hb}$ and $G\whb$ are finite groups.  
\item 
The left multiplication operator $L_\nxi$, with $x_i\in G\hb$, permutes the normalized basis $\overline{\mtc B}$. 
\item 
The left multiplication operator $L_{\mui}$, with $\mui\in G\whb$, permutes the basis ${\wdb}$.
\ene
\ec
\subsection{Some results on $P$ and $\whp$.}
They denote $\prod_{i \in \mathcal{I}} \nxi$ and $\prod_{j \in \mathcal{I}}\muj$, respectively.
\bl\label{ev:grplike}
Let $\hbz$ be an abelian normalizable hypergroup, and let $x_i\in G\hb$ be a grouplike element. Then 
\beq
(\prod_{j \in \mathcal{I}}\muj)(\nxi)=\pm 1.
\eeq
\el
\bpf
By Equation (\ref{mwa:nzd}), defining the multiplication in the dual, $$(\prod_{j \in \mathcal{I}}\muj)(\nxi)=\prod_{j \in \mathcal{I}}\muj(\nxi).$$
By Equation (\ref{eq:dualinv}), defining the involution in the dual, and then Lemma \ref{grp:id}, $$\muj(\nxi)\mu_{j^\#}(\nxi)=|\muj(\nxi)|^2=1.$$ 
We group together the factors $\muj$ and $\mu_{j\#}$ in the product above. Now, if $\muj=\mu_{j^\#}$, then $\muj(\nxi)\in \mathbb R$, and therefore, $\muj(\nxi)=\pm 1$ by Lemma \ref{grp:id}.
\epf
\bn{defn}\label{sgn:def}
Let $\hbz$ be  an abelian normalizable hypergroup, and let $x_i\in \ghb$ be a grouplike element. Then $(\prod_{j \in \mathcal{I}}\muj)(\nxi)$ will be denoted $\sgn(x_i)$.
\end{defn}
The dual version of the above result is the following:
\bl\label{dl:ev:grplike}
Let $\hbz$ be an abelian normalizable hypergroup, and let $\muj\in G\whb$ be a grouplike element. Then 
\beq
\muj(\prod_{i \in \mathcal{I}}\nxi)=\pm 1.
\eeq
\el
By duality, we can also define:
\bn{defn}\label{dual:sgn:def}
Let $\hbz$ be  an abelian normalizable hypergroup, and let $\muj\in G\whb$ be a grouplike element. Then $\muj(\prod_{i \in \mathcal{I}}\nxi)$ will be denoted $\sgn(\muj)$.
\end{defn}
\bp
Following the assumptions of Corollary \ref{all:in:one}, let $x_i\in \ghb$ and $\muj\in G\whbz$ be grouplike elements. Then $$\sgn(x_i) =  \det(L_\nxi) \text{ and } \sgn(\muj) = \det(L_{\muj})$$ are the signature of the permutations of $\overline{\mtc B}$ by $L_{\nxi}$, and of $\widehat{\mtc B}$ by $L_{\muj}$, respectively.
\ep
\bpf
By Corollary \ref{all:in:one}, the left multiplication operator $L_{\nxi}$ permutes $\overline{\mtc B}$, but the signature of this permutation is  the determinant $\det(L_\nxi)$. Now $\nxi=\sum_{j \in \mathcal{I}} \muj(\nxi)F_j$, so $\det(L_\nxi) = \prod_{j \in \mathcal{I}}\muj(\nxi)$. Idem for $\sgn(\mu_j)$.
\epf
\bn{defn} \label{perfect}
A \invisible{symmetric}hypergroup $\hb$ is called \emph{perfect} if it has no nontrivial grouplike elements.
\end{defn}
\section{Galois symmetries for hypergroups}\label{gal}
Let $(H, \mtc B)$ be an abelian rational hypergroup (i.e. $N^k_{ij}\in \mathbb Q$).  For any standard element $x_i \in \mathcal B$, we can write $x_i=\sum_{j \in \mathcal{I}} \alij F_j$, with $\alij=\muj(x_i) \in \mathbb C$. Let $\mbk=\mathbb Q(\alij)$ be the field obtained by adjoining all  $\alij$ to $\mbq$. 
\br\label{alijjp}
Note that if $\alij=\al_{ _{ij'}}$ for all $i$, then $\mu_j(x_i)=\mu_{j'}(x_i)$, thus $\mu_j=\mu_{j'}$ and therefore $j=j'$.
\er
\subsection{Permutation of characters $\mu_j$}\label{tauperm} Let $\hb$ be any abelian rational hypergroup.  For any character $\muj:H\ra \comp$ and $\sg \in \gal( \overline{\mathbb Q}/\mathbb Q)$, define  $\sg.\muj\in H^*$ as the linear function on $H$ such that $[\sg.\muj](x_i)=\sg(\muj(x_i))=\sg(\alij)$, for all $x_i \in {\mtc B}$.
\bl \label{def:sghstar}
For any abelian rational hypergroup $\hb$, the function $\sg.\muj$ is an algebra map. Thus $\gal( \overline{\mathbb Q}/\mathbb Q)$ acts on the set $\wdb$.
\el
\bpf
Suppose that $x_{i_1}x_{i_2}=\sum_{k \in \mathcal{I}} N^k_{\io,\itw}x_k$. Since $N^k_{\io,\itw}\in \mathbb Q$, 
$$[\sg.\muj]({x_{i_1}}{x_{i_2}})=[\sg.\muj](\sum_{k \in \mathcal{I}} N^k_{\io,\itw}x_k)=\sum_{k \in \mathcal{I}} N^k_{\io,\itw}[\sg.\muj](x_k)=\sum_{k \in \mathcal{I}} N^k_{\io,\itw}\sg(\muj(x_k)).$$
On the other hand,
\begin{eqnarray*}
[\sg.\muj]({x_{i_1}})[\sg.\muj]({x_{i_2}}) &=& \sg(\muj({x_{i_1}}))\sg(\muj({x_{i_2}}))=\sg(\mu_j({x_{i_1}})\muj({x_{i_2}}))\\ &= & \sg(\muj({x_{i_1}}{x_{i_2}}))=\sg(\muj(\sum_{k \in \mathcal{I}} N^k_{\io,\itw}x_k))\\ &= &\sum_{k \in \mathcal{I}} N^k_{\io,\itw}\sg(\muj(x_k)).
\end{eqnarray*}
Thus $[\sg.\muj]({x_{i_1}}{x_{i_2}})=[\sg.\muj]({x_{i_1}})[\sg.\muj]({x_{i_2}})$, so that $\sg.\muj$ is an algebra map. Now,
$$ [\sg.(\sg'.\muj)](x_i) = \sg([\sg'.\muj](x_i)) = \sg(\sg'(\muj(x_i))) = (\sg \sg')(\muj(x_i)) = [\sg \sg' . \muj](x_i)$$
 It follows that $\gal( \overline{\mathbb Q}/\mathbb Q)$ acts on the set $\wdb$. 
\epf
Thus, there is a permutation $\tau_\sg$ of ${\mtc I}$ such that $\sg.\muj=\mu_{\tau_\sg(j)}$, and therefore
\beq\label{tau} 
\sg(\alij)=\al_{ _{i\tau_\sg(j)}}\;\; \text{for all}\;i,j \in {\mtc I}.
\eeq
\bc\label{Galois theory}
Let $\hb$ be an abelian rational hypergroup. Then $\mbq\subseteq \mbk$ is a Galois extension.
\ec
\bpf
Equation \eqref{tau} shows that $\sg(\mbk)\subseteq \mbk$, for  all $\sg \in \gal( \overline{\mathbb Q}/\mathbb Q)$. Therefore, this is a normal extension and consequently a Galois extension (as ${\mathrm char}(\mathbb{Q}) = 0$).
\epf
\br\label{injectivity} 
The map $\ro:\galkq\ra \mathbb S_{\mathcal I}, \sg\sent \tau_{\sg}$ is injective. Indeed, if $\tau_\sg=\tau_{\sg'}$, then $ \sg(\alij)=\sg'(\alij)$, i.e. $ \sg^{-1}\circ \sg'(\alij)=\alij$, for all $i,j$. Thus $\sg=\sg'$. 
\er
\bp \label{sg:nk:eq:ntauk}
Let $(H, \mathcal B)$ be an abelian rational hypergroup. For any $\sg \in \galkq$, let $\tau = \tau_{\sg}$, with the above notations,
\beq\label{sg:nj}
\sg( n_j)={n_{\tau(j)}}.
\eeq
\ep
\bpf 
Applying $\sg\in \galkq$ to the orthogonality relation \eqref{first:gen:orth}, we obtain:
$$
\sum_{i \in \mathcal I} h_i \mu_{\tau(j)}(x_i) \mu_{\tau(k)}(x_{i^*}) = \delta_{j,k} \sg(n_j),
$$
On the other hand, by the same orthogonality relation:
$$
\sum_{i \in \mathcal I} h_i \mu_{\tau(j)}(x_i) \mu_{\tau(k)}(x_{i^*})=\delta_{\tau(j), \tau(k)}{n_{\tau(j)}}.
$$
Therefore, with $j=k$, we get that $\sg(n_j)={n_{\tau(j)}}$.
\epf
\bn{defn}\label{h:integral}
A hypergroup $\hb$ is \emph{$h$-integral} if $h_i :=1/N^1_{i,i^*} \in \mathbb Z$, for all $i \in \mathcal I$.
\end{defn}
\br\label{owh}
Let $\hbz$ be the Grothendieck ring of any weakly-integral fusion category $\cc$. By \cite[Theorem 2.13]{O3} and Equation \eqref{whj:eq}, $\widehat{h}_i = \frac{n\hbz}{n_i}$ is an integer, for all $i \in \mathcal I$. Thus, the dual $\wkc$ is $h$-integral. See also Remark \ref{hj}.
\er
\bl\label{sgccj}
Suppose that $\hbz$ is an abelian  normalizable hypergroup with a $h$-integral dual. Then it is weakly-integral. Moreover, if it is rational, then $\widehat{h}_k=\widehat{h}_{\tau(k)}$, for all $k \in \mathcal I$.
\el
\bpf
By Lemma \ref{dual:order} and Definition \ref{def:order}, $$n\hbz=n\whbz=\sum_{i \in \mathcal I} \wh_i\in \mathbb Z.$$ If $\hbz$ is rational, then Equation \eqref{sg:nj} can be written as
$\sg \big( \frac{n\hbz}{\widehat{h}_{i}} \big)=\frac{n\hbz}{\widehat{h}_{\tau(i)}}.$ In particular,   
$
\widehat{h}_k=\widehat{h}_{\tau(i)}
$, since $\whbz$ is $h$-integral.
\epf
\section{Burnside's vanishing property for hypergroups}\label{v}
Let $\hb$ be any  hypergroup. An element $x_i\in \mtc B$ is called a \emph{vanishing element} if there is $\muj\in \wdb$ such that $\muj(x_i)= 0$. Otherwise, $x_i$ is called \emph{a non-vanishing element}.  

Let $\mtc B_0$ be the set of all vanishing elements of $\mtc B$, and $\mtc B_1$ be the set of all non-vanishing elements of $\mtc B$. Thus $\mtc B=\mtc B_1\sqcup \mtc B_0$. By the proof of Lemma \ref{grp:id}, if $\hb$ is a RN-hypergroup then $\ghb\subseteq \mtc B_1$. In the abelian case, since $x_i=\sum_{j \in \mathcal I} \muj(x_i)F_j$, it follows that $x_i$ is non-vanishing if and only if it is invertible in $H$ (i.e. $\det(L_{x_i})$ nonzero).  Thus, in this case, $\mtc B_1$ coincides with the set of all standard elements that are invertible in $H$.

In any abelian normalizable hypergroup $\hbz$:
\beq\label{wp}
\whp:=\prod_{j \in \mathcal I} \muj=\sum_{i \in \mathcal I} \whp(\nxi)\wtf_i.
\eeq
\bp\label{gen:muj:id}
Let $\hbz$ be any abelian normalizable hypergroup. Then:
\beq\label{gen}
\prod_{j \in \mathcal I} \mu_j=\sum_{x_i\in \mtc B_1}\det(L_{\nxi})\wtf_i
\eeq
where $\det(L_{\nxi})$ is the determinant of the left multiplication operator by $\nxi$ on $H$.
\ep
\bpf
Recall that $\wtf_i (\nxj)=\delta_{i,j}$, and $\frac{x_j}{d_j}=\sum_{i \in \mathcal I} \mui(\frac{x_j}{d_j})F_i$. Thus, for all $j \in \mathcal I$,
$$\left[\sum_{i\in \mtc I}\det(L_{\nxi})\wtf_i \right](\nxj) 
= \det(L_{\nxj}) = \prod_{i \in \mathcal I} \mu_i(\nxj) = \left[\prod_{i \in \mathcal I} \mu_i \right](\nxj)$$
It follows that $$\prod_{i \in \mathcal I} \mu_i = \sum_{i\in \mtc I}\det(L_{\nxi})\wtf_i = \sum_{x_i\in \mtc B_1}\det(L_{\nxi})\wtf_i.$$ 
Indeed, $x_i \in {\mtc B_1}$ if and only if $\det(L_{\nxi})$ is nonzero, justifying the last equality.%
%
%
\epf
Dually, applying the above result to $\whbz$, we obtain a decomposition of the basis $\wdb=\wdb_0\sqcup \wdb_1$. Proposition \ref{gen:muj:id} implies the following:
\bp\label{dl:gen:v:id}
Let $\hbz$ be any  abelian normalizable hypergroup. Then:
\beq\label{gen:xi}
\prod_{i \in \mathcal I}\nxi=\sum_{\mu_j\in \wdb_1}\det(L_{\muj})F_j,
\eeq
where $\det(L_{\muj})$ is the determinant of the left multiplication operator by $\muj$ on $H^*$.
\ep
Recall the definition of an abelian Burnside hypergroup from Definition \ref{v-property}. Note that the decomposition $\mtc B=\mtc B_0\sqcup \mtc B_1$ from above implies the following:
\bp\label{grouplike:set}
An abelian hypergroup $\hb$ is  Burnside  if and only if we have $\ghb=\mtc B_1$.
 \ep
\bc \label{bsd:iff:pt:id}
An abelian normalizable hypergroup $\hbz$ is Burnside if and only if the following equality holds:
\beq\label{gen:gh}
\prod_{j \in \mathcal I}\mu_j=\sum_{x_i\in \ghb}\sgn(x_i)\wtf_i,
\eeq
where the notation $\sgn(x_i)$ was introduced in Definition \ref{sgn:def}.
\ec
\bpf
Immediate from Propositions \ref{gen:muj:id} and \ref{grouplike:set}.
\epf
\bc\label{B:muj:sq:id}
An abelian normalizable hypergroup $\hbz$ is Burnside if and only if
\beq\label{gen:gh:squared}
(\prod_{j \in \mathcal I}\mu_j)^2=\sum_{x_i\in \ghb}\wtf_i.
\eeq
\ec
\bc\label{odd-order}
An abelian normalizable hypergroup $\hbz$, such that $\ghb$ has odd order, is Burnside if and only if 
\beq\label{muj:gen:odd}
\prod_{j \in \mathcal I}\muj=\sum_{x_i \in G(H, B)}\wtf_i.
\eeq
\ec
\bpf
Every $\nxi\in \ghb$ has odd order, so $\sgn(x_i)=1$.
\epf
\bn{defn} \label{def:d-B}
An abelian normalizable  hypergroup $\hbz$ is called \emph{dual-Burnside} if $\whbz$ is Burnside; in other words, if for all $\muj\in  \wdb$, the following are equivalent:
\bne
\item 
For all $x_i\in \mtcb$ then $\muj(x_i) \neq 0$,  
\item
$\muj\in G\widehat{\hb}$, i.e $\muj$ is a grouplike element.
\ene
This can be shortened to $G\widehat{\hb} = \wdb_1$.
\end{defn}
Here are the dual of Corollaries \ref{bsd:iff:pt:id}, \ref{B:muj:sq:id} and \ref{odd-order}:
\bc\label{dl:brsd:id}
An abelian normalizable hypergroup $\hbz$ is dual-Burnside if and only if:
\beq\label{gen:db:xi}
\prod_{i \in \mathcal I}\nxi=\sum_{\mu_j\in G\whb}\sgn(\muj)F_j.
\eeq
\ec
\bc\label{d-b:sq:idm}
An abelian normalizable hypergroup $\hbz$ is dual-Burnside if and only if :
\beq\label{d-b:xi:sq:idm}
(\prod_{i \in \mathcal I}\nxi)^{2}=\sum_{\mu_j\in G\whb}F_j.
\eeq
\ec
\bc\label{dl:odd-order}
An abelian normalizable hypergroup $\hbz$, with $\wqghb$ of odd order, is dual-Burnside if and only if 
\beq\label{d-b:odd:idm}
\prod_{i \in \mathcal I}\nxi=\sum_{\mu_j\in \wqghb}F_j.
\eeq
\ec
\subsection{Inequalities for $P$ and $\whp$}
Recall that $P := \prod_{i \in \mathcal I} \frac{x_i}{d_i}$ and $\whp := \prod_{j\in \mathcal I} \muj$.
\bp  \label{eignv:whp}
Let $\hbz$ be an abelian normalizable hypergroup satisfying the hypothesis of Lemma \ref{l:fp:ineq}. Then $|\whp(\nxi)|\leq 1$. Moreover, the equality holds if and only if $x_i\in G\hb$.
\ep
\bpf
First, suppose that $\hb$ is RN. Then $|\muj(\nxi)|\leq 1$ by Frobenius-Perron theory. This implies that $|\whp(\nxi)|=\prod_{j \in \mathcal I} |\muj(\nxi)|\leq 1$. Moreover, the equality holds if and only if $|\nxi(\muj)|=1$ for all $j\in \mtc I$, if and only if $x_i$ is grouplike, by Lemma \ref{grp:id}.

Next, suppose that $\whbz$ is RN. Since
$\muj=\sum_{i \in \mathcal I} \muj(\nxi)\wtf_i$, then $\{\muj(\nxi)\}_{i \in \mathcal I}$ are the eigenvalues of $L_{\mu_j}$. By Frobenius-Perron theory, $|\muj(\nxi)| \le \fp(\muj) = 1$, as $\whbz$ is normalized. Thus, $|\whp(\nxi)|\leq 1$, again. The rest is as before.
\epf
\bp\label{whp-sq:idmpt} 
Let $\hbz$ be an abelian normalizable hypergroup. Then, $\whp^2$ is an idempotent if and only if $\whp(\nxi)=\pm 1$, for all $x_i \in \mathcal B_1$.
\ep
\bpf
By Equation \eqref{wp}, $\whp^4=\whp^2$ if and only if  $\whp(\nxi)^4=\whp(\nxi)^2$ for all $i\in \mtc I$. Therefore, $\whp^2$ is idempotent if and only if $\whp(\nxi)=0$ or $\whp(\nxi)^2=1$, for all $i\in \mtc I$. Since $\whp(\nxi)=\prod_{j \in \mathcal I} \muj(\nxi)$, the second case occurs if and only if $x_i \in \mathcal B_1$, if and only if $\whp(\nxi)=\pm 1$.  
\epf
\bc
Let $\hbz$ be an abelian normalizable hypergroup satisfying the hypothesis of Lemma \ref{l:fp:ineq}. Then, $\whp^2$ is an idempotent if and only if $\hbz$ is Burnside. 
\ec
\bpf
By Proposition \ref{whp-sq:idmpt}, $\whp^2$ is an idempotent if and only if $\whp(\nxi)=\pm 1$, for all $x_i \in \mathcal B_1$. But, by assumption and Proposition \ref{eignv:whp}, $|\whp(\nxi)| = 1$ if and only if $x_i \in G\hb$. Thus, $\whp^2$ idempotent implies that ${\mathcal B_1} \subset G\hb$. But $G\hb \subset \mathcal B_1$, so the equality holds, i.e.  $\hb$ is Burnside.

Conversely, if $\hb$ is Burnside, then $\whp^2=\sum_{x_i\in G\hb} \wtf_i$, by Corollary \ref{B:muj:sq:id}, which is an idempotent. 
\epf
Let $\hbz$ be an abelian normalizable hypergroup, and $P:=\prod_{i \in \mathcal I} \nxi$. Then
\beq\label{dual:wp}
P=\sum_{j \in \mathcal I} \muj(P)F_j.
\eeq  

Here are the dual version of the last three results: 

\bp  \label{eignv:p}
Let $\hbz$ be an abelian normalizable hypergroup satisfying the hypothesis of Lemma \ref{l:fp:ineq}. Then $|\muj(P)|\leq 1$. Moreover, the equality holds if and only if $\muj\in G\whb$.
\ep 

\bp \label{whp-sq:idmptdual}
Let $\hbz$ be an abelian normalizable hypergroup. Then, $P^2$ is an idempotent if and only if $\muj(P)=\pm 1$, for all $\muj \in  \widehat{\mathcal B}_1$.
\ep 

\bc \label{corP2DualB}
Let $\hbz$ be an abelian normalizable hypergroup satisfying the hypothesis of Lemma \ref{l:fp:ineq}. Then, $P^2$ is an idempotent if and only if  $\hbz$ is dual-Burnside.
\ec 

\subsection{The analogue of Burnside's theorem}
We will prove the following generalization of Theorem \ref{burnside-fr}, using the same approach as in \cite[Theorem 2]{b-galois}.\bt\label{burnside}
Let $\hbz$ be an abelian rational normalizable hypergroup with ${h}$-integral dual.  Let $x_i\in \mtc B$. If $h_i|d_i|^2>1$, and
\beq\label{alg}
h_i|\muj(x_i)|^2\in \mathbb A
\eeq 
for all $j\in \mtc I$.  Then $x_i \in \mathcal B_0$
\et
\bpf
For any $i \in \mathcal I$, denote $\mtc T_i:=\{j\in \mtc I\;|\; \mu_{j}(x_{i})= 0\}$ and $\cd_i  :=\mtc I\setminus (\mtc T_i\cup \{1\})$.
\vsk
We need to show that $\mtc T_i\neq \emptyset$, for any $x_i\in \mtc B$ satisfying the theorem's hypothesis.

The second orthogonality relation (\ref{second:orth}) and Equation (\ref{whj:eq}) implies
\beq\label{orthofes2-z}
\sum_{j \in \mathcal I} \whj|\muj(x_i)|^2=\frac{n\hbz}{h_i}.
\eeq

Since $\widehat{h}_1=1$, the above equation can be written as:
$$
\frac{n\hbz}{h_i}=|\di|^2+\sum_{j \in \cd_i }{\whj}|\mu_j(x_i)|^2$$
which gives that 
\beq\label{gnn1}
1=\frac{n\hbz}{h_i|\di|^2}-\sum_{j \in \cd_i  }\frac{\whj|\muj(x_i)|^2}{|\di|^2}.
\eeq

On the other hand,
$$
n\hbz=\sumjtom \whj=1+\sum_{j\in \mtc T_i}\whj+\sum_{j\in \cd_i  }\whj.
$$
Therefore, Equation \eqref{gnn1} can be written as
\beq\label{gnn2}
1=\frac{1+\sum_{j\in \mtc T_i}\whj}{h_i|\di|^2}-\big(\sum_{j \in \cd_i  }\frac{\whj|\muj(x_i)|^2}{|\di|^2}-\sum_{j\in \cd_i  }\frac{\whj}{h_i|\di|^2}
\big).
\eeq
Thus, in order to finish the proof, it is enough to show that
\beq\label{gnn2}
\big(\sum_{j \in \cd_i  }\frac{\whj|\mu_j(x_i)|^2}{|\di|^2}-\sum_{j\in \cd_i  }\frac{\whj}{h_i|\di|^2}\big)\geq 0,
\eeq
since then, $\frac{1+\sum_{j\in \mtc T_i}\whj}{h_i|\di|^2}\geq 1$, i.e. $1+\sum\limits_{j\in \mtc T_i}\whj\geq h_i |\di|^2$; and $h_i|\di|^2>1$ implies $\mtc T_i\neq\emptyset$. 

The inequality \eqref{gnn2} can be written as
\beq\label{gnn3}
(\sum_{j \in \cd_i}\whj)^{-1}(\sum_{j \in \cd_i  }{\whj h_i|\mu_j(x_i)|^2})\geq 1.
\eeq
\noindent On the other hand, the weighted AM-GM inequality gives that 
\beq\label{gnn4}
(\sum_{j \in \cd_i}\whj)^{-1}(\sum_{j \in \cd_i  }{\whj h_i|\mu_j(x_i)|^2})
\geq 
\bigg( \prod_{j\in \cd_i  } \big(h_i |\muj(x_i)|^2\big)^{\whj} \bigg)^{(\sum_{j \in \cd_i}\whj)^{-1}}.
\eeq
By rational assumption, $h_i$ is rational and $\muj(x_i)$ is an algebraic number. Equation \eqref{tau} implies that the set $\cd_i$ is stable under the action of $\gal( \overline{\mathbb Q}/\mathbb Q)$. It follows that the product  
$$P_i:=\prod_{j\in \cd_i  }(h_i|\muj(x_i)|^2)^{\whj}$$ is fixed by the action of $\gal( \overline{\mathbb Q}/\mathbb Q)$, since ${\what{h}_{\tau(j)}}=\whj$ by Lemma \ref{sgccj}. Thus $P_i$ is a rational number. On the other hand, each factor of $P_{i}$ is an algebraic integer (since $\whj \in \mathbb Z_{_{>0}}$), and therefore, the entire product is an integer. Since it is also positive, it must be at least $1$.
\epf
Note that Theorem \ref{burnside-fr} follows from the above theorem, since in the case of a fusion ring,  $h_i=1$, and $\muj(x_i)\in \mathbb A$, for all $i,j\in {\mtc I}$.
\bp \label{arnnh:inv}
Let $\hb$ be a rational ARN-hypergroup with $h$-integral dual. Let $x_i\in \mtc B$ be an invertible element in $H$ satisfying also Equation \eqref{alg}, for all $j\in \mtc I$. Then $x_i\in \ghb$.
\ep
\bpf
Suppose that $x_i$ is not a grouplike element, i.e. $h_i|d_i|^2>1$, by Lemma \ref{grouplike:nn}. Then, all the assumptions of Theorem \ref{burnside} are satisfied, so that, $x_i \in \mathcal B_0$, contradiction with $x_i$ invertible.
\epf
\bc
Any rational ARN-hypergroup $\hb$ with $h$-integral dual, such that Equation \eqref{alg} holds for all $i,j \in \mathcal I$, is Burnside.
\ec
\bpf
By Proposition \ref{grouplike:set}, it is enough to show that ${\mtc B_1} \subseteq \ghb$. Let $x_{i}\in \mtc B_1$, thus $x_i$ is invertible, so by Proposition \ref{arnnh:inv}, $x_i\in \ghb$.
\epf
\bc \label{fr:inv}
Any commutative fusion ring  with $h$-integral dual is Burnside.
\ec
The following corollary recovers \cite[Theorem 2]{b-galois}.
\bc\label{hhj}
A weakly-integral fusion category with a commutative Grothendieck ring is Burnside.
\ec
\bpf
It follows from Corollary \ref{fr:inv} by Remark \ref{owh}.
\epf
We can write a kind of dual version of Theorem \ref{burnside}: 
\bt\label{dual:burnside}
Let $\hbz$ be an abelian normalizable $h$-integral hypergroup with a rational dual. 
Let $\muj\in \wdb$ such that $\whj>1$ and
\beq\label{duality:alg}
\whj \left| \muj(\frac{x_i}{d_i})\right| ^2\in \mathbb A,
\eeq
for all $i\in \mtc I$. Then, $\mathcal B_0$ is non-empty.
\et
\bpf 
The condition from Equation \eqref{alg} applied on $\whbz$ becomes exactly Equation \eqref{duality:alg}. So we can apply Theorem \ref{burnside} to $\whbz$. Thus $\muj \in \widehat{\mathcal{B}}_0$, meaning the existence of $x_i \in \mathcal B$ such that $\muj(\frac{x_i}{d_i})=0$, so $x_i \in \mathcal B_0$.
\epf
\section{Kernels of fusion rings, Brauer's theorem}\label{kernels}

Let $\hb$ be an ARN-hypergroup. Let $\hb_+$ be the set of all elements $x\in H$ such that $x=\sum_{b\in \mtcb}x_b b$ with $x_b\in \mathbb R_{\geq 0}$.
A standard element $b \in \mtc B$ is called a \emph{constituent of $x\in \hb_+$} if $x_b > 0$. Recall that $\wdb$ is the set of all algebra morphisms $\muj:H\ra \comp$. Let $\psi \in \wdb$. For all $x\in \hb_+$, 
$$ |\psi(x)| = |\psi(\sum_{b\in \mtcb}x_b b)| \le \sum_{b\in \mtcb}x_b |\psi(b)| \le  \sum_{b\in \mtcb}x_b \fp(b) = \fp(x), $$ 
by Frobenius-Perron theory.

\bn{defn}\label{sub:hyp:gp}
Let $\mtc S\subseteq \mtc B$ be a subset such that
\begin{itemize}
\item $\mtc S$ is closed under the involution of $\hb$,
\item $L:=\comp[\mtc S]$ is a subalgebra of $H$.
\end{itemize}
Then $\hs$ is called a \emph{sub-hypergroup} of $\hb$.
\end{defn}
The notion of kernel of object of a fusion category from \cite{bmonat} can be extended to arbitrary ARN-hypergroups.
\bn{defn} \label{def:ker}
Let $\hb$ be any ARN-hypergroup.
For any $\psi\in \wdb$, define 
$$\ker_{\whb}(\psi):=\{x \in \mtc B\;|\; \psi(x)=\fp(x)\}.$$
\end{defn}
\bl\label{ker-sub-hyp}
Let $\hb$ be any ARN-hypergroup and $\psi\in \wdb$. Then $\ker_{\whb}(\psi)$ is (the basis of) a sub-hypergroup of $\hb$.
\el
\bpf 
Let $x, y\in \ker_{\whb}(\psi)$, i.e. $\psi(x)=\fp(x)$ and $\psi(y)=\fp(y)$. Suppose that $xy=\sum_{z \in \mtc B'}N^z_{xy}z$, where $\mtc B'$ is the subset of $\mtc B$ consisting of all standard elements for which $N^z_{xy}>0$.

It follows that
$$\fp(x)\fp(y)=\psi(x)\psi(y)=\psi(xy)=\sum_{z \in \mtc B'}N^z_{xy}\psi(z).$$
Using the triangle inequality for the complex absolute-value,
\begin{eqnarray*}
\fp(x)\fp(y)&=&\left| \sum_{z \in \mtc B'}N^z_{xy}\psi(z)\right|\\&\leq &\sum_{z \in \mtc B'}N^z_{xy}|\psi(z)|\leq \sum_{z \in \mtc B'}N^z_{xy}\fp(z)\\&=&\fp(x)\fp(y).
\end{eqnarray*}
It is easy to deduce that $\ker_{\whb}(\psi)$ is a sub-hypergroup of $H$.
\epf
Dually, we can define the following:
\bn{defn}\label{defn:ker:h} 
Let $\hb$ be an  ARN-hypergroup. Let $x\in H_+$. Define
$$\kerhb(x):=\{\psi\in \wdb \;|\; \psi(x)=\fp(x)\}.$$

\end{defn}
The dual version of Lemma \ref{ker-sub-hyp} implies that $\kerhb(x)$ is (the basis of) a sub-hypergroup of $\whb$ if $\whb$ is a RN-hypergroup.
\subsection{Brauer's theorem for ARN-hypergroup} \label{sub:brau}
For any hypergroup $\hb$, we define the bilinear function $m$ on $H$ by $$m(\sumitom \al_ix_i, \sumitom \beta_ix_i)=\sumitom \frac{\al_i\beta_i}{h_i}.$$ In the abelian case, Equation \eqref{tau-form} implies
\beq\label{m:eq}
m(x, y)=\tau(xy^*)=\sumjtom \frac{1}{n_j}\muj(x)\muj(y^*),\;\; \text{for all}\;\;x,y \in H.
\eeq
\bl  \label{lem:KerSum}
Let $\hb$ be an ARN-hypergroup. If $x:=\sum_{s \in \mtc B}p_ss\in \hb_+$ then
$$
\kerhb(x)=\bigcap_{ _{\{s \in \mtc B|\;p_s>0\}}}\kerhb(s).
$$
\el
\bpf
Note that $\fp(x)=\sum_{\{s \in \mtc B\;|\;p_s>0\}}p_s\fp(s)$. Let $\mu_j \in \widehat{B}$ then $$|\mu_j(x)|\leq \sum_{\{s \in \mtc B|\;p_s>0\}}p_s|\mu_j(s)|\leq \sum_sp_s\fp(s)=\fp(x).$$ Thus, $\muj\in \kerhb(x)$ if and only if $\muj\in \bigcap_{ _{\{s \in \mtc B|\;p_s>0\}}}\kerhb(s)$.
\epf
%
The following theorem is a generalization of Brauer's theorem. Although the proof is the same as that in \cite{bmonat}, we include it here for the sake of completeness.

For $x\in \hb_+$, we denote by $\langle x \rangle$ the sub-hypergroup generated by $x$, i.e. whose basis is the set of constituents of all powers $x^n$, with  $n\geq 1$.
\bt \label{Brauer}
Let $\hb$ be an ARN-hypergroup, and $x \in \hb_+$. Then $\langle x \rangle=H$ if and only if $\kerhb(x)=\{\mu_1\}$.
\et
\bpf 
Partition the set $\mtc I= \bigsqcup_{l \in \mathcal P} \mtc I_l$ such that the values $\muj(x)$ are constant on each component $\mtc I_l$ of the partition, but distinct otherwise (i.e. for all $l,l' \in \mathcal P$, $j \in \mtc I_l$, $j' \in \mtc I_{l'}$,  $\muj(x) = \mu_{j'}(x)$ if and only if $l=l'$).
Let $s\in \mtc B$.  By Equation (\ref{m:eq}),
$$
{m(x^n,s)}{}=\sumjtom \frac{1}{n_j} \muj(x)^n  \muj(s^*)=\sum_{l\in \mtc P}\left(\sum_{j\in \mathcal I_l}\frac{1}{n_j}\muj(s^*)\right)\alpha_l^n=\sum_{l\in \mtc P}a_l\alpha_l^n
$$
with 
$$a_l:=\sum_{j\in \mathcal I_l}\frac{1}{n_j}\muj(s^*), \text{ and } \alpha_l := \muj(x) \text{ when } j \in \mathcal I_l.$$
Note that $s \not \in \langle x \rangle$ if and only if $m(x^n,s)=0$ for all $n$, if and only if 
$Mv=0$, with $M$ be the matrix $(\alpha_l^n)$ and $v$ be the vector $(a_l)$. But $M$ is invertible because its determinant is nonzero, because all $\alpha_l$ are distinct (by Vandermonde determinant formula). Thus, $s \in \langle x \rangle$ if and only if there is $l \in \mathcal P$ such that $a_l$ is nonzero.

If $\kerhb(x)=\{\mu_1\}$ then the component $\mathcal I_1$ containing $1$ is just $\{1\}$. Thus, $a_1 = \frac{\mu_1(s^*)}{n_1} = \frac{\fp(s)}{\fp\hb}$ is nonzero for all $s \in \mathcal B$. It follows by above that $\langle x \rangle = H$.

Finally, $\muj \in \ker_\hb(x)$ if and only if $x \in \ker_{\whb}(\muj)$, if and only if $\langle x \rangle \subset \ker_{\whb}(\muj)$ by Lemma \ref{ker-sub-hyp}. So if $\langle x \rangle=H$ then  $\ker_{\whb}(\muj) = H$, i.e $\muj=\mu_1$.
\epf

\subsection{The cointegral $\lam_\ls$ in terms of kernels}
Recall from \S \ref{introd} that for any ARN-hypergroup $\hb$, we denoted by $\lam_H$ the primitive central idempotent $F_1$ corresponding to $\mu_1=\fp$. Let $\ls$ be a sub-hypergroup of $\hb$, and let $\lam_\ls$ be the idempotent in $L \subset H$ corresponding to $\mu_1\big|_\ls$. By Equation (\ref{fj}), we get that:
\begin{equation} \label{eq:lamls}
\lam_\ls = \frac{1}{n\ls} \sum_{x_i \in \mathcal S} h_i d_i x_i. 
\end{equation}
In particular, $\lam_\hs = \lam_\htt$ if and only if $\hs = \htt$.
\bt \label{ker:lam:rel}
For any $x \in \hb_+$, let $\mathcal{I}_x$ be the index set of $\kerhb(x)$. Then
$$\lam_{\langle x \rangle}=\sum_{j\in \mathcal{I}_x}F_{j}$$
\et
\bpf
Suppose that $\langle x \rangle=\hs$.  Let $G_1, \dots , G_s$ be the primitive central idempotents of $\hs$, and $\psi_1, \dots , \psi_s$ be their associated characters. As above, let us take $G_1 =\lam_{\hs}$, and consequently, $\psi_1=\fp$.

We look at the restrictions of the characters $\mu_j$ at $\hs$. Then, there is a surjective function $f :\mathcal{I} \ra \{1, \dots, s \}$ such that $\muj \big|_{\hs}=\psi_{f ( j )}$, for all $j \in \mathcal I$.

With the above notation, we are reduced to show that 
$f^{-1}(\{1\}) = \mtc I_x$ . By Theorem \ref{Brauer}, since $\langle x \rangle=\hs$, then  $\ker_{ _\hs}(x) = \{\psi_1\}$. On the other hand, $$\mu_j (x) = \mu_j\big|_{ _\hs} (x) = \psi_{f ( j )}(x).$$ Thus, $\mu_j \in \kerhb(x)$ if and only if
$\psi_{f(j)}\in \ker_{ _\hs}(x)$, i.e. $f(j)=1$.
\epf
 
\section{The adjoint sub-hypergroup and grouplike elements} \label{sec:adjsub}

\subsection{Notion of the center of a standard element}
Let $\hb$ be an ARN-hypergroup. We also define the notion of \emph{center} of  a standard element $x_i\in \mtc B$:
$$
Z_{ _\hb}(x_i)=\{\psi \in \wdb\;|\; |\psi(x_i)|=d_i\}.
$$

Define also the dual notion of the center, 
$$
Z_{ _\whb}(\psi)=\{x_i\in \mtcb\;|\; |\psi(x_i)|=d_i\}.
$$
As in Lemma \ref{ker-sub-hyp}, it is easy to see that $Z_{ _\whb}(\psi)$ is  a sub-hypergroup of $\hb$. Idem for $Z_{ _\hb}(x_i)$, if $\whb$ is a RN-hypergroup.
\subsection{Adjoint sub-hypergroup}\label{adjoint:hyp}
\bp\label{ker:int}
Let $\hb$ be an ARN-hypergroup and $I(1):=\sumitom h_ix_ix_{i^*}\in H$. Then 
$$\kerhb(I(1))=\bigcap_{i \in \mathcal I}  Z_{ _\hb}(x_i)$$
\ep
\bpf
Note that $$\psi(I(1)) =  \sumitom h_i\psi(x_ix_{i^*})=\sumitom h_i|\psi(x_i)|^2\leq \sumitom h_id_i^2=\fp\hb.$$ So if $\psi \in \kerhb(I(1))$, then $\psi \in \bigcap_{i \in \mathcal I}  Z_{ _\hb}(x_i)$. The converse is immediate.
\epf
\bn{defn}\label{def:hbad}
Denote by $\hbad$ the sub-hypergroup $\langle I(1) \rangle$ of $H$ generated by $I(1)$. It is called the \emph{adjoint sub-hypergroup} of $H$.
\end{defn}
\subsection{On the support $\jhs$}
\bn{defn}\label{def:supp:j}
Let $\hs$ be a sub-hypergroup of an abelian  hypergroup $\hb$. Its \emph{support} is the subset
$\jhs \subset \mathcal I$ such that
\beq\label{j:def}
\lamhs=\sum_{j\in \jhs}F_j.
\eeq
\end{defn}
\noindent In particular, by Equation (\ref{eq:lamls}), $\jhs = \jhtt$ if and only if $\hs = \htt$. 

The next result generalizes \cite[Lemma 6.4]{ccc-march}.
\bp\label{nj:max}
Let $\hb$ be an ARN-hypergroup. With the above notations,
$$
\jhbad=\{j\in \mtc I\;|\; n_j=\fp\hb\}.
$$
\ep
\bpf
By Theorem \ref{ker:lam:rel} and Definition \ref{def:hbad}, $\jhbad = {\mathcal I}_{I(1)}$, the index set of $\kerhb(I(1))$. Thus, $$\jhbad = \{j\in \mtc I\;|\; \muj(I(1)) = \fp(I(1))\}.$$ But $\fp(I(1)) = \fp\hb$, and $\muj(I(1)) = n_j$ by Equation \eqref{first:gen:orth}.
\epf
\bt\label{grpk:eq:hbad}
Let $\hb$ be an  ARN-hypergroup. Then $\muj\in G\whb$ if and only $j \in \jhbad$. 
\et
\bpf 
It follows from Proposition \ref{nj:max} and Lemma \ref{max:formal:codeg}.
%
\epf
\bc\label{lam:ad}
Let $\hb$ be an ARN-hypergroup. Then
\beq\label{lam:ad:eq}
\lam_{ _{\hbad}}=\sum_{\muj\in G\whb}F_{j}\eeq
\ec
\bpf
It follows from Equation \eqref{j:def} and Theorem \ref{grpk:eq:hbad}.
\epf
Here are dual versions of Theorem \ref{grpk:eq:hbad} and Corollary \ref{lam:ad}:
\bt\label{dl:grp:eq:hbad}
Let $\hbz$ be an abelian normalizable hypergroup such that $\whbz$ is RN. Then $x_i\in \ghb$ if and only $i \in \mtc I_{{\whb}_{ad}}$. 
\et
\bc\label{dual:lam:ad} Let $\hbz$ be an abelian normalizable hypergroup such that $\whbz$ is RN. Then
\beq\label{dual:lam:ad:eq}
\lam_{ _{{\whb}_{ad}}}=\sum_{x_i\in \ghb}\wtf_{i}.
\eeq
\ec
\subsection{Proof of Theorems \ref{hwrn} and \ref{hwrndual}}
Recall Definition \ref{def:supp:j}, which defines the support $\mtc{I}_{\hs} \subset \mtc{I}$ of a sub-hypergroup $\hs$ within an abelian hypergroup $\hb$.
For the sake of brevity, we define $|L|:=n(L, \cs, \mu_1)$.
\bt\label{ccc-eq}
Let $\hs$ be a sub-hypergroup of an abelian normalizable hypergroup $\hbz$. With the above notations, we have:
\beq\label{ccc-eqq}
\sum\limits_{x_i\in \cs}\wtf_i=\frac{|L|}{|H|}\big(\sum\limits_{j\in \mtc{I}_{\hs}}\wh_j\mu_j\big).
\eeq
\et
\bpf
Applying $\mtf$ to Equation \eqref{j:def} we have that
\beq\label{mtf-1}
\mtf(\lamhs)=\sum\limits_{j\in \jhs}\mtf(F_j)\numeq{\ref{mtf-muk}}|H|\big(\sum\limits_{j\in \jhs}\frac{\mu_j}{n_j}\big).
\eeq
On the other hand, from Equation \eqref{f0:gen} we have 
$
\lamhs=\frac{1}{|L|}\big(\sum\limits_{x_i\in \cs} h_{i^*}d_{i^*}x_i\big)
$
and therefore:
\beq\label{mtf-2}
\mtf(\lamhs)=\frac{1}{|L|}\big(\sum\limits_{x_i\in \cs} h_{i^*}d_{i^*}\mtf(x_i)\big) \numeq{\ref{mtf-eqq-2}} \frac{1}{|L|}\big(\sum\limits_{x_i\in \cs} h_{i^*}d_{i^*}\frac{|H|}{d_{i^*} h_{i^*}}\wtf_{i^*}\big)=\frac{|H|}{|L|}\big(\sum\limits_{x_i\in \cs} \wtf_{i^*}\big)
\eeq
Comparing \eqref{mtf-1} and \eqref{mtf-2}, we obtain:
\beqn
\sum\limits_{x_i\in \cs}\wtf_i=|L|\big(\sum\limits_{j\in \mtc{I}_\hs}\frac{\mu_j}{n_j}\big).
\eeqn
According to Equation \eqref{whj:eq}, we have $\wh_j = \frac{|H|}{n_j}$, thus concluding the proof.
\epf
\emph{Proof of Theorem \ref{hwrn}}:
By applying Corollary \ref{B:muj:sq:id}, we get that 
$$
\big(\prod_{j \in \mtc{I}} \mu_j\big)^2 = \sum_{x_i \in \ghb} \wtf_i.
$$
On the other hand, by utilizing Theorem \ref{ccc-eq} with $\cs = G\hb$, we arrive at the desired conclusion.

\emph{Proof of Theorem \ref{hwrndual}}: To derive the result through duality from Theorem \ref{hwrn}, we need to consider the normalized basis $(\frac{x_i}{d_i})_{i \in \mathcal{I}}$. Recall that by definition, $\langle x_i x_{i^*}, x_1 \rangle^{-1} = h_i$. Therefore, $\langle \frac{x_i}{d_i} \frac{x_{i^*}}{d_{i^*}}, x_1 \rangle^{-1} = d_i d_{i^*} h_i$. Consequently, $d_i d_{i^*} h_i \frac{x_i}{d_i} = d_{i^*} h_i x_i$, which precisely matches the summand on the RHS.

\subsection{Proof of Theorems \ref{hbz:hyp} and \ref{dual:hbz:hyp}.}
\bpf
By Corollary \ref{B:muj:sq:id}, $\hb$ is Burnside if and only if  
$$
\big(\prod_{j \in \mathcal I}\muj\big)^2=\sum_{x_i \in G(H, B)}\wtf_i,
$$
if and only if Equation \eqref{muj:gen:hyp} holds, by Corollary \ref{dual:lam:ad}, (the dual version of) Equation \eqref{f0} and Lemma \ref{lem:absym}.
%
\epf
By duality, we obtain directly Theorem \ref{dual:hbz:hyp}.
\subsection{On the sub-hypergroup generated by $P$}

\bl \label{lem:lambdaLs}
For any sub-hypergroup  $\hs$ of an ARN-hypergroup $\hb$, then ${\mathcal I}_{\hs} = \bigcap_{s \in \mathcal S} {\mathcal I}_s$. Thus, $s \in \mathcal S$ if and only if $s\lam_{\hs} = d_s\lam_{\hs}$.
\el
\bpf
Note that $\hs = \langle x \rangle$, with $x = \sum_{s \in \mathcal S}s \in \hb_{+}$. Recall that ${\mathcal I}_x$ is the index set of $\kerhb(x)$. By Theorem \ref{ker:lam:rel}, ${\mathcal I_x} = {\mathcal I_{\langle x \rangle}}$, and by Lemma \ref{lem:KerSum}, ${\mathcal I_x} = \bigcap_{s \in \mathcal S} {\mathcal I_s}$. The first sentence follows.
Recall that $x = \sum_{j \in \mathcal I} \muj(x) F_j$, thus $x\lam_{\hs} = \sum_{j \in \mathcal I_{\hs}} \muj(x) F_j$. But the first sentence means that $j \in {\mathcal I}_{\hs}$ if and only if $\muj(s) = \fp(s)$, for all $s \in \mathcal S$. So if $s \in \mathcal S$ then $s\lam_{\hs} = d_s\lam_{\hs}$. Regarding the converse, we deduce by positivity in Equation (\ref{eq:lamls}) that $sS \subseteq S$. Specifically, since $s = sx_1 \in sS$, it follows that $s \in S$.
\epf

In particular, for an ARN-hypergroup with $\mu_1 = \fp$, then $\jhs = \{j \in {\mtc I} \text{ such that } \mu_j|_\hs = \fp\}$. Let $\hs, \htt$ be two sub-hypergroups of an ARN-hypergroup $\hb$. Let $\hs\vee\htt$ be the sub-hypergroup generated by $\hs$ and $\htt$.

\bp\label{vee}
With above notations, $\mtc I_{ _{\hs\vee\htt}}=\jhtt\cap \jhs$. In other words, $\lam_{\hs\vee\htt} = \lam_{\hs} \lam_{\htt}$.
\ep
\bpf  
Let ${\mathcal S} \vee {\mathcal T}$ denotes the basis of $\hs\vee\htt$. By Lemmas \ref{lem:lambdaLs} and \ref{ker-sub-hyp}, 
\begin{align*}
\jhs\cap \jhtt = \bigcap_{a \in {\mathcal S} \cup {\mathcal T}} {\mathcal I}_a &= \{ j \in {\mathcal I} \ | \ \muj \in \ker_{\hb}(a), \ \forall a \in {\mathcal S} \cup {\mathcal T} \} \\
&= \{ j \in {\mathcal I} \ | \ {\mathcal S} \cup {\mathcal T} \subset \ker_{\whb}(\muj)\} \\
&= \{ j \in {\mathcal I} \ | \ {\mathcal S} \vee {\mathcal T} \subset \ker_{\whb}(\muj)\} \\
&= \{ j \in {\mathcal I} \ | \ \muj \in \ker_{\hb}(a), \ \forall a \in {\mathcal S} \vee {\mathcal T} \} \\
&=  \bigcap_{a \in {\mathcal S} \vee {\mathcal T}} {\mathcal I}_a = \mtc I_{ _{\hs\vee\htt}}. \qedhere
\end{align*}
%
%
\epf

\bp \label{j:incl}
For any two sub-hypergroups $\hs, \htt$ of an ARN-hypergroup $\hb$, then $\hs\subseteq \htt$ if and only if $\jhtt\subseteq \jhs$.
\ep
\bpf
By Lemma \ref{lem:lambdaLs}, if $\hs\subseteq \htt$ then $\jhtt\subseteq \jhs$. By Proposition \ref{vee}, if $\jhtt\subseteq \jhs$, then $\lamhs \lamhtt=\lamhtt$. Pick $s \in \mathcal S$, then $$s\lamhtt = s\lamhs \lamhtt = d_s \lamhs \lamhtt = d_s \lamhtt, $$ thus $s \in \mathcal T$, by Lemma \ref{lem:lambdaLs}.
\epf
\bc
For any two elements   $P, Q\in \hb_+$ of an ARN-hypergroup $\hb$, then  $\langle P \rangle \subseteq \langle Q \rangle$ if and only if $\kerhb(Q)\subseteq \kerhb(P)$.
\ec
For any $\hbz$, recall that $P: =\prod_{i \in \mathcal I}\nxi$.
\bp\label{gen:p}
For any abelian RN-hypergroup $\hb$ then
$$\langle P^2 \rangle=\hbad.$$ 
\ep
\bpf
By Definition \ref{def:hbad}, $\hbad := \langle I(1) \rangle$, so we are reduced to show that $\mathcal I_{\langle P^2 \rangle} = \mathcal I_{\langle I(1) \rangle}$, so (by Theorem \ref{ker:lam:rel}) that $\mathcal I_{P^2} = \mathcal I_{I(1)}$ i.e. $\kerhb(P^2)=\kerhb(I(1))$.

Note that $\big|\psi(P)\big|=\prod_{i \in \mathcal I}\big|\psi(\nxi)\big|\leq 1$. Thus, if $\psi \in \kerhb(P)$, then $\big|\psi(\nxi)\big|=1$, i.e $\psi\in \kerhb(I(1))$, by Proposition \ref{ker:int}. Thus, $\kerhb(P) \subseteq \kerhb(I(1))$, and the same argument shows that $\kerhb(P^2) \subseteq \kerhb(I(1))$. Conversely, if $\psi\in \kerhb(I(1))$ then $\psi(\nxi \frac{x_{i^*}}{d_{i^*}})=\big|\psi(\nxi)\big|^2=1$. Thus, if $i\neq i^*$ then, grouping together $x_i$ and $x_{i^*}$ in $P$, we obtain $\psi(\nxi)\psi(\frac{x_{\istar}}{\distar})=\big|\psi(\nxi)^2\big|=1$. Thus,
$$
\psi(P)=\prod_{i={i^*}}\psi(\nxi).
$$
Note that if $i=i^*$, then $\psi(\nxi)=\pm1$, since in this case $\psi(x_i)=\psi(x_\istar)\in \mathbb R$. Thus, if $\psi\in \kerhb(I(1))$, then $\psi(P)=\pm 1$ and $\psi(P^2)=1$, i.e. $\psi \in \kerhb(P^2)$.
%
\epf
\bc
For any ARN-hypergroup $\hb$, then $\hbad\subseteq \langle P \rangle$.
\ec
\bc \label{vars:grp:rn}
Let $\hb$ be an ARN hypergroup. The following holds:
\bne
\item 
$\hbad=\comp$ if and only if $\hb$ is pointed,
\item 
$\hbad=\comp$ if and only if $\whb$ is pointed,
\item 
$\hbad=\hb$ if and only if $\whb$ is perfect.
\ene
\ec
\bpf First, recall that $\hbad = \langle I(1) \rangle$, so $\hbad=\comp$ if and only if $I(1) = n\hb x_1$, if and only if $\hb$ is pointed, as $I(1) = \sum_{i \in \mathcal I} h_i x_i x_{i^*}$. Next, by Equation (\ref{eq:lamls}), $\hbad=\comp$ if and only if $\lambda_{\hbad} =  \lambda_{\comp} = x_1 = \id = \sum_{i \in \mathcal I} F_i, $
if and only if $ G\whb = \widehat{\mathcal B}$, by Equation \ref{lam:ad:eq}, meaning that $\whb$ is pointed. Similary, $\hbad=\hb$ if and only if $\lambda_{\hbad} = F_1$, if and only if $G\whb = \{\mu_1\}$, meaning that $\whb$ is perfect.
\epf

\br \label{rem:uni}
By \S \ref{universe}, $\hbad=\hb$ if and only if the universal grading group is trivial, if and only if there is no non-trivial grading (see Corollary \ref{cor:uni}).
\er
\bc
A simple non-pointed ARN-hypergroup $\hb$ has a perfect dual.
\ec
\bpf
By simplicity $\hbad = \comp$ or $\hb$. By non-pointed assumption and Corollary \ref{vars:grp:rn}, $\hbad \neq \comp$, therefore $\hbad = \hb$, and so $\whb$ is perfect by Corollary  \ref{vars:grp:rn}.
\epf
\bc\label{vars:grp:dl-rn}
Let $\hbz$ be an abelian normalizable hypergroup with $\whbz$ RN. Then the following holds:
\bne
\item 
$\whb_{ad}=\comp$  if and only if $\whb$ is pointed,
\item 
$\whb_{ad}=\comp$ if and only if $\hb$ is pointed,
\item 
$\whb_{ad}=\whb$ if and only if $\hb$ is perfect. 
\ene
\ec
\section{On the universal grading of hypergroups}\label{universe}
The universal grading construction for  fusion rings from \cite[\S 2]{NG} works word for word in the ARN-hypergroup settings. For the completeness of the paper, we sketch its main steps below.
\subsection{Based modules are completely reducible}
\bn{defn}
Given an ARN-hypergroup $\hb$, we define a  based left $\hb$-module as a pair $(M, \mtc M)$ where $M$ is a $H$-module with $\mtc M$ a finite free $H$-basis for $M$, i.e. for all $x_i\in \mtc B$ and all $m_j\in \mtc M$, then $x_im_j=\sum_k t^k_{ij}m_k$; and moreover, for all $i,j,k$ then $t^k_{ij}$ is zero if and only if $t^{j}_{i^*k}$ is zero.
\end{defn}
The last condition above means that $(x_i m_j,m_k)$ vanishes if and only if $(m_j,x_{i^*} m_k)$ vanishes, for the inner product on $M$ such that for all $m_i, m_j\in \mtc M$  
\beq\label{based:module:ip}
(m_i, m_j)=\delta_{i,j}.
\eeq

It is easy to deduce that any based $\hb$-modules is completely reducible, as in \cite[Lemma 2.1]{O}. Indeed, if $N \subset M$ and $HN \subset N$ then $(N, N^{\perp}) = \{0\} = (HN, N^{\perp}) = (N, HN^{\perp})$. Thus, $HN^{\perp} \subset N^{\perp}$.

Let $\hb$ be a RN-hypergroup and $x,y \in \hb_+$. 
We say that $x$ contains $y$ (or $y$ is a constituent of $x$) if the difference $x-y$ is in $\hb_{+}$.
\bn{defn}
Let $\hb$ be a \invisible{symmetric}RN-hypergroup, $G$ be a finite group. Suppose that $\mtc B=\sqcup_{g\in G}\ccb_g$ is a partition of $\mtc B$, and let  $H_g:=\comp[\mtc B_g]$. The decomposition $\hb=\bigoplus_{g\in G}(H_g,\mtc B_g)$ is called a \emph{grading} of $\hb$ by $G$ if $H_gH_h\subseteq H_{gh}$, for all $g,h\in G$. The grading is called \emph{faithful} if $\mtc B_g\neq\emptyset$, for all $g\in G$.
\end{defn}

Please note that a faithful grading group of an abelian hypergroup must also be abelian.
Let $\hb=\bigoplus_{g\in G}(H_g,\mtc B_g)$ be a faithful grading of $\hb$,
and define $R_g:=\sum_{x_i\in \mtc B_g}h_id_ix_i$, for all $g\in G$.
The analogue of  \cite[Proposition 8.20]{eno-annals} holds in the settings of ARN-hypergroup:
\beq\label{fprg}
\fp(R_g)=\frac{\fp\hb}{|G|}.
\eeq
Indeed, consider $R:=\sum_{g \in G}R_g$, and note that  $\frac{1}{\fp\hb}R=F_1$, the idempotent corresponding to $\fp$ in Corollary \ref{fz:gen}. Thus, $RR_h=\fp(R_h) R$. Since  $H_gH_h\subseteq H_{gh}$, we obtain from the previous equality that $R_gR_h=\fp(R_h) R_{gh}$, for all $g,h\in G$. Passing to $\fp$, we get that $\fp(R_g) \fp(R_h) = \fp(R_h) \fp(R_{gh})$;  but $\fp(R_h) \neq 0$ by faithfulness, so $\fp(R_{gh})=\fp(R_g)$, for all $g,h \in G$.
\subsection{Constructing the universal grading}
\bp 
Let $\hb$ be a RN-hypergroup. For all $x \in H$, define $I(x):=\sum_{j\in \mtc I}h_j x_j x x_{j^*}$. With the above notations:
\bne
\item $I(x)$ is central in $H$. 
\item Every based left $H_{ad}$-submodule $M$ of $H$ is also a $H_{ad}$-subbimodule of $H$.
\item 
A based $H_{ad}$-subbimodule $M$ of $H$ is indecomposable if and only if it is indecomposable as left $H_{ad}$-module.
\ene
\ep
\bpf
\bne
\item For all $i \in \mtc I$, by Frobenius reciprocity (\ref{eq:FR}): $$ x_i I(x) = \sum_{j\in \mtc I}h_j x_ix_j x x_{j^*} = \sum_{j,k\in \mtc I}h_j N_{i,j}^k x_k x x_{j^*} = \sum_{j,k\in \mtc I}h_k N_{k^*,i}^{j^*} x_k x x_{j^*} $$ $$= \sum_{k\in \mtc I} h_k x_k x (\sum_{j\in \mtc I} N_{k^*,i}^{j^*} x_{j^*})  = \sum_{k\in \mtc I}h_k x_k x x_{k^*} x_i = I(x)x_i.$$
\item Recall that $H_{ad}:=\langle I(1) \rangle$. By (1), $m\rrad^n=\rrad^nm$, thus by RN-assumption, $m H_{ad} = H_{ad} m \subset M$, for all $m\in M$, meaning that $M$ is also a based right $H_{ad}$-submodule.
\item Clear. \qedhere
\ene
\epf
Let $H=\oplus_{a\in A} H_a$ be a decomposition (into indecomposable) of $H$ as a $H_{ad}$-bimodule. This also corresponds to a partition $\mtc B=\sqcup_{a\in A}\mtc B_a$. Let $1$ be an element of $A$ such that $H_1=H_{ad}$.
\subsubsection{Definition of $a^*$ for any $a\in A$} Since $H_a$ is an indecomposable $H_{ad}$-bimodule then
$(H_a)^*$ is also an indecomposable $H_{ad}$-bimodule.  Let $a^*$ be an element of $A$ such that $(H_a)^* = H_{a^*}$. Next lemma is as \cite[Lemma 3.4]{NG}.
\bl 
If $x_a, y_a\in \mtc B_a$ then $x_ay_a^*\in H_{ad}$.
\el
\bpf 
Note that $M := H_a \cap (H_{ad} y_a)$ is a $H_{ad}$-submodule of $H_a$ (containing $y_a$). But $H_a$ is indecomposable (so irreducible by complete reducibility). It follows that $M=H_a$, so $H_a \subset H_{ad} y_a$. Thus $x_a \in H_{ad} y_a$ and $x_a y_a^* \in H_{ad} y_ay_a^* \subset H_{ad}$.
\epf
\bt
For all $a,b \in A$, there is $c \in A$ such that $H_a H_b \subset H_c$. This induces a group structure on $A$ given by $a b := c$. Moreover, $a^*$ is the inverse of $a$, and  $1$ the unit of $A$.
\et
\bpf
As for the proof of \cite[Theorem 3.5]{NG}.
\epf
\bn{defn}
Above group $A$ is called the \emph{universal grading group} of $\hb$. It is denoted $U_H$ when the standard basis $\mtc B$ is implicitly understood.
\end{defn}
\bc \label{cor:uni}
Every ARN-hypergroup $\hb$ has a canonical faithful grading by the group
$U(H)$. Any other faithful grading of $\hb$ by a group $G$ is determined by a surjective
 group morphism $\pi:U_H \ra G$.
\ec
\bpf
Let $\hb=\oplus_{g \in G}(S_g, \mtc D_g)$ be another faithful grading. It follows that $xx^*\in S_1$, for any $x\in \mtc D_g$, but $xx^* \in H_{ad}$ therefore, $\mtc B_{ad}\subseteq \mtc D_1$. Moreover, each $S_g$ is a $\had$-subbimodule of $H$. Since $H_a$ is an indecomposable module, then $H_a\subseteq S_{\pi(a)}$, for some well-defined $\pi(a)\in G$. It follows that $\pi: U_H\ra G$ is a surjective group morphism. Indeed, for all $a,b \in A$, then $H_a \subset S_{\pi(a)}$ and $H_b \subset S_{\pi(b)}$, thus $H_aH_b \subset S_{\pi(a)} S_{\pi(b)} \subset S_{\pi(a)\pi(b)}$. But $H_aH_b \subset H_{ab}$, thus $H_aH_b \subset H_{ab} \cap S_{\pi(a)\pi(b)} \neq \emptyset$, so $H_{ab} \subset S_{\pi(a)\pi(b)}$, meaning that $\pi(ab) = \pi(a) \pi(b)$. Finally, the surjectivity comes from the faithfulness.
\epf

\subsection{On the dual universal grading group and invertible}
In this subsection, we use the notation $H$ for a hypergroup $\hb$, and $\wdh$ for its dual $\whb$. The goal of this subsection is to prepare for the proof (\S \ref{sub:RegComp}) of the following:
\bt\label{univ:grd}
Let $H$ be an ARN-hypergroup.
Then,
$$U(H)\simeq \gwdh$$
\et
Dually, we have the following:
\bt \label{dual:univ:grd}
Let $\hbz$ be an abelian normalizable hypergroup such that $\whbz$ is RN. Then
$$
U(\wdh)\simeq G(H).
$$
\et 
Recall from Lemma \ref{dual:grp:id} that for any abelian normalizable hypergroup $\hbz$ then $\muj \in \gwdh$ if and only if $|\omega_{ij}| = 1$, for all $i \in \mathcal I$, with $\omega_{ij} := \mu_j(\frac{x_i}{d_i})$.
\bn{defn} 
Let $\hbz$ be an abelian normalizable hypergroup. For any $x_i\in \mtc B$, we define the linear character $\omega_{i} := \nxi\big|_{\gwdh}$, the restriction to $\gwdh$ of the linear character $\nxi:\wdh\ra \comp$, i.e. $\omega_{i}(\muj) = \omega_{ij}$.
\end{defn}
Let $\hbz$ be an abelian normalizable hypergroup. 
Given $\psi\in \widehat{\gwdh}$, let
\beq\label{bpsi:def}
\mtc B_{\psi}:=\{x_i\in \mtc B\;\mid\;\omega_i=\psi\} =\{x_i\;|\; \muj(\nxi)=\psi(\muj),\;\text{for all}\;\muj\in \gwdh\}.
\eeq
\bp Let $\hb$ be an ARN-hypergroup. Then $(\mtc B_{\psi})_{\psi\in \widehat{\gwdh}}$ defines a grading on $\hb$.
\ep
\bpf 
For any $\muj \in \gwdh$, $\psi,\phi \in \widehat{\gwdh}$, $x_i \in B_{\psi}$ and $x_k \in B_{\phi}$, then
$$
\psi\phi(\muj)=\psi(\muj)\phi(\muj)=\muj(\nxi)\muj(\nxk)=\muj(\nxi\nxk).$$
By the triangle inequality:
$$
1=|\psi\phi(\muj)|=|\muj(\nxi\nxk)|\leq \frac{1}{d_id_k} \sum_{l \in \mathcal I}d_lN^l_{ik}|\muj(\nxl)|\leq \frac{1}{d_id_k} \sum_{l \in \mathcal I}d_lN^l_{ik}=1,
$$
which forces $\muj(\frac{x_l}{d_l})=\muj(\frac{x_i}{d_i})\muj(\nxk)$ if $N^l_{ik}$ is nonzero, meaning that if $x_i\in \mtc B_{\psi}$ and $x_k\in \mtc B_{\phi}$, then $x_ix_k$ has all constituents in $\mtc B_{\psi \phi}$.
\epf
\subsection{Computation of the regular components} \label{sub:RegComp}
For any grading component $\mtc B_\psi$, we define (as above):
$$R_{\psi}=\sum_{x_i\in \mtc B_\psi}h_id_ix_i.$$
\bl \label{lem:FjRpsi}
Let $\hb$ be an ARN-hypergroup. For all $\muj\in \gwdh$,
\beq\label{fj:grd}
F_j=\frac{1}{n(H)}\sum_{\psi\in \widehta{\gwdh}}\psi(\muj)R_{\psi^{-1}}.
\eeq
\el
\bpf
For all $\muj\in \gwdh$, by Equation (\ref{fj}):
\begin{eqnarray*}
F_j&=&\frac{1}{n_j}\sumitom h_i\muj(x_{i^*})x_i=\frac{1}{n_j}\sum_{\psi \in \widehat{\gwdh}}\sum_{x_i\in \mtc B_{\psi}}h_i\muj(x_{i^*})x_i=\\&=&\frac{1}{n_j}\sum_{\psi\in \widehta{\gwdh}}\psi^{-1}(\muj)R_{\psi}=\frac{1}{n_j}\sum_{\psi\in \widehta{\gwdh}}\psi(\muj)R_{\psi^{-1}}.
\end{eqnarray*}
Finally, if $\muj\in \gwdh$ then $n_j=n(H)$ by Corollary \ref{frmlcdg:grplike}.
\epf
\bc
Let $\hb$ be an ARN-hypergroup. 
The set $\mtc B_{\psi}$ is not empty.
\ec
\bpf
Lemma \ref{lem:FjRpsi} shows that the linear span of $\{F_j\;|\;\muj\in \gwdh \}$ lies inside the span of $\{\mtc R_{\psi}\;|\; \psi \in \widehat{\gwdh} \text{ with }\mtc B_{\psi}\neq\emptyset \}$. But $\widehat{\gwdh} \simeq \gwdh$ as finite abelian group, so have same order. Thus, an argument involving the dimension of these vector spaces implies that $\mtc B_\psi$ is not empty, for any $\psi \in \widehat{\gwdh}$.
\epf
By Equation \eqref{fj:grd} and the second orthogonality relation (\ref{second:orth}) applied on the (hyper)group $\gwdh$, it follows that:
\bl \label{lem:rpsi}
Let $\hb$ be an ARN-hypergroup. Then 
$$
R_{\psi}=\frac{n(H)}{|\gwdh|}\big(\sum_{\muj\in \gwdh} \psi(\muj)F_j\big).
$$
\el
%
%
{\bf Proof of Theorem \ref{univ:grd} and Theorem \ref{dual:univ:grd}}\\
By Lemma \ref{lem:rpsi}, and Corollary \ref{lam:ad}, if $\psi=1$, then 
$$R_1=\frac{n(H)}{|\gwdh|}\big(\sum_{\muj\in \gwdh} F_j\big)= \frac{n(H)}{|\gwdh|} \lam_{\had}.
$$
Proposition \ref{j:incl} proves that $\mtc B_1=\mtc B_{ad}$. Therefore  the above grading $$\mtc B=\bigsqcup_{\psi \in \widehat{G(\wdh)}}\mtcb_\psi$$ coincides with the universal grading of $\hb$. Thus $U_H=\gwdh$.  By duality, $U_{\wdh}=G(H)$ if $\whbz$ is RN.
\subsection{On the perp of the adjoint sub-hypergroup} Recall from Definition \ref{def:dualzb} that an ARN-hypergroup is called dualizable if $\whb$ is also (A)RN.

For any sub-hypergroup $\hs$ of a dualizable ARN-hypergroup $\hb$, define 
$$
\mtc S^\perp:=\{\muj\;|\;\muj(s)=\fp(s)\;\text{for all}\; s\in \mtc S\}.
$$
With the above notations, $\mtc S^\perp=\bigcap_{s \in \mtc S}\ker_\hb(s)$.  Moreover, \cite[Proposition 2.11]{hdk} implies that $\mtc S^\perp$ is a sub-hypergroup of $\whb$, and $(\mtc S^{\perp})^\perp=\mtc S$.
\bc  \label{cor:adperp}
Let $H$ be a dualizable ARN-hypergroup. With the above notations,
$$
\gwdh^\perp=H_{ad} \text{ and } G(H)^\perp=\wdhad.
$$
Thus, by $(\mtc S^{\perp})^\perp=\mtc S$,
$$
\gwdh=(H_{ad})^\perp \text{ and } G(H)=\wdhad^\perp.
$$
\ec
\bpf
Consider the trivial character $\psi_1=x_1\big|_\gwdh$. It follows that
\begin{eqnarray*}
\had=B_{\psi_1}&=&\{x_i\;|\; \muj(\nxi)=1\;\text{for all}\;\muj\in \gwdh\}\\&=&\bigcap_{\muj\in \gwdh}\ker_{\wdh}(\muj)=\gwdh^\perp.
\end{eqnarray*}
By duality, we obtain the second equality.
\epf
\subsection{Quotients of hypergroups} 
In this subsection, we define the quotient hypergroup of a RN-hypergroup $\hb$ by a sub-hypergroup $\hs$. This is the analogue of the quotient construction from \cite{hdk}, in the case of probability groups. Recall that a \emph{probability group} is a \invisible{symmetric}normalized RN-hypergroup. 

We can define the following equivalence relation on $\mtc B$. 
For two elements $a,b \in \mtc B$, say $a\sim_\mtcs b$ if and only if there are $s_1,s_2\in \mtcs, x\in \mtcb$ such that $m(x, as_1)>0$ and $m(x,s_2b)>0$.

Let $[a]_\mtcs$ denote the equivalence class of any $a\in \mtc B$, with respect to $\sim_\mtcs$. Define $H//\mtcs$ as the set of all these equivalence classes of $\sim_\mtcs$.  Recall that $\lam_\mtcs$ is the central idempotent of $\mtcs$ corresponding to $\mu_1=\fp$.

For two elements $a,b \in \mtc B$, we can see that $[a]_\mtcs=[b]_\mtcs$ if and only if $\lam_\mtcs a\lam_\mtcs=\lam_\mtcs b\lam_\mtcs$. Therefore there is a set bijection
$$\begin{array}{cccc}
             \phi : & \comp[H//\mtcs]  & \to & \lam_\mtcs H \lam_\mtcs \\ 
                            & [a]_\mtcs & \mapsto & \lam_\mtcs a\lam_\mtcs 
\end{array}$$
Then, it is easy to verify that $H//\mtcs$ becomes a RN-hypergroup, with the multiplication inherited from $\lam_\mtcs H\lam_\mtcs$, via the above isomorphism. We denote by $\overline{m}({[c]},\;[a][b])$ the multiplicity  structure of $H//\mtcs$. Therefore,
$$[a][b]=\sum_{[c]\in H//\mtcs}\overline{m}({[c]},\;[a][b])[c].$$
We write shortly $[a]$ instead of $[a]_\mtcs$ when no confusion is possible.
Following \cite{hdk}, if $\hb$ is abelian, we can show that 
\beq\label{coset}\overline{m}({[c]},\;[a][b])=\sum_{w \in [c]}m(w, ab).
\eeq

It was proven in \cite[Proposition 2.11]{hdk} that if $\hb$ is an abelian  dualizable probability group then,
$$\begin{array}{cccc}
             \al : & \mtcs^\perp  & \to & \widehat{H//\mtcs} \\ 
                            & \psi & \mapsto & \al(\psi)
\end{array},$$
with $\al(\psi)([a]_\mtcs):=\psi(a)$, is an isomorphism of probability groups. 
It is clear that the same results remains true for dualizable ARN-hypergroups.

\subsection{Applying Harrison's results}
In this subsection, we also use the short notation $H$ (or $\mtc B$) for a hypergroup $(H, \mtc B)$, if $\mtc B$ (or $H$) are implicitly understood. For any sub-hypergroup $\ls$ of a dualizable ANR-hypergroup $\hb$, note that \cite[Proposition 2.11]{hdk}  implies that there is an isomorphism of hypergroups:
$$
\wdh//\mtc S^\perp\simeq \widehat{\mtc S}, \;\;[\mu]_{\mtc S^\perp}\sent \mu\big|_{L}.
$$
In particular, for $\mtc S=H_{ad}$,
\beq\label{w:qout:gh}
\wdh//G(\wdh)\simeq \widehat{H_{ad}}
\eeq
Dually, for $\mtcs=\wdhad$, we obtain that
\beq\label{quot:gh}
H//G(H)\simeq \widehat{\wdhad}
\eeq
For $S=G(H)$, then
$\wdh//\wdhad\simeq \widehat{G(H)}$,
and dually, for $\mtcs=\gwdh$, then
$H//H_{ad}\simeq \widehat{\gwdh}.$
\section{On lower and upper central series of hypergroups}\label{lucs}
Let $H$ be a RN-hypergroup.
Following \cite{NG}, let $H^{(0)}=H$, $H^{(1)}=\had$, and  $H^{(n)}=H^{(n-1)}_{ad}$, for all $n\geq 1$.

\bn{defn}
The non-increasing sequence 
$$
H=H^{(0)}\supseteq  H^{(1)}\supseteq \dots \supseteq H^{(n)} \supseteq \dots
$$
will be called the \emph{upper central series}.
\end{defn}
\bn{defn} 
Let $\hs$ be a sub-hypergroup  of an ARN-hypergroup $\hb$. Let $\mtcs^\coo$ be the set of standard elements $x\in \mtc B$ such that $xx^*\in \mtc \mtcs$. The \emph{commutator} $\hs^{\coo}$ of $\hs$ in $\hb$ as the sub-hypergroup generated by $\mtcs^\coo$.
\end{defn}

As in \cite[Remark 4.9]{NG}, it follows that the linear span of $ \mtcs^\coo$ is already a sub-hypergroup of $H$. Moreover, \cite[Lemma 4.15]{NG} works as well in the settings of ARN-hypergroups. Thus, 
\beq\label{gn}
(\mtc S^{\coo})_{ad}\subseteq \mtc S\subseteq (\mtc S_{ad})^{\coo}
\eeq
for any sub-hypergroup $\mtc S$ of $\mtc B$. By applying $\perp$ to \eqref{gn}, we obtain
\beq\label{gn:perp}
\big((\mtc S^{\coo})_{ad}\big)^\perp\supseteq \mtc S^\perp\supseteq \big((\mtc S_{ad})^{\coo}\big)^\perp.
\eeq
For an ARN-hypergroup $\hb$, define $H_{(0)} = \comp$ and $H_{(n)}=(H_{(n-1)})^{\coo}$, for all $n\geq 1$. Then, $$H_{(1)} = \comp^{\coo} = G\hb = H_{pt}.$$
\bn{defn}
Let $H$ be an ARN-hypergroup. The non-decreasing sequence
$$
\comp=H_{(0)}\subseteq  H_{(1)}\subseteq \dots \subseteq H_{(n)} \subseteq \dots
$$
 will be called the \emph{lower central series} of $H$.
 \end{defn}
Similarly to \cite[Theorem 4.16]{NG}, we can prove that $H^{(n)}=\comp$ if and only if $H_{(n)}=H$, where $H$ is an ARN-hypergroup.
\bn{defn}
An ARN-hypergroup is termed \emph{nilpotent} if there exists some $n \in \mathbb{N}$ such that $H^{(n)} = \comp$. The smallest such $n$ for which this condition is satisfied is referred to as the \emph{nilpotency class} of $H$.
\end{defn}
\bl\label{sperpad} 
Let $\hb$ be a dualizable ARN-hypergroup.
For any sub-hypergroup $\hs$,
\beq\label{main:perp:ad}
( \mtc S^\perp)_{ad}\subseteq ( \mtc S^\coo)^\perp.
\eeq
\el
\bpf
It reduces to show that $\muj \star \mu_{j^\#}\in ({\mtc S}^\coo)^\perp$, for all $\muj\in {\mtc S}^\perp$.

Let $x\in {\mtc S}^\coo$ be a standard element. It follows that $xx^*\in {\mtc S}$, and therefore, $\muj(x)\mujstar(x)=\muj(xx^*)=\fp(x)^2$, for all $\muj\in {\mtc S}^\perp$. This implies 
\begin{eqnarray*}
[\muj\star\mu_{j^\#}](\nox)&=&\muj(\nox)\mu_{j^\#}(\nox)\\&=&\frac{1}{\fp(x)^2}\muj(x)\mu_{j^\#}(x)=1,
\end{eqnarray*}
which shows that $\muj\star\mujstar\in ({\mtc S}^\coo)^\perp$. Thus $({\mtc S}^\perp)_{ad}\subseteq (S^\coo)^\perp$.
\epf
\bp\label{fst:incl}
Let $\hb$ be a dualizable ARN-hypergroup. With the above notations, for all $n\geq 0$,
\beq\label{fiq}
\wdh^{(n)}\subseteq (H_{(n)})^\perp
\eeq
and 
\beq\label{dual:fiq}
H^{(n)}\subseteq (\wdh_{(n)})^\perp.
\eeq 
\ep
\bpf
We will prove the first inclusion by induction on $n$. If $n=0$, then $$\wdh^{(0)} = \wdh = \comp^{\perp} =  H_{(0)}^\perp.$$
Now, suppose that $\wdhn\subseteq \hn^\perp$. Then
$$
\wdh^{(n+1)}=(\wdhn)_{ad}\subseteq (\hn^\perp)_{ad}\subseteq (\hn^\coo)^\perp=H_{(n+1)}^\perp.
$$
By duality, we obtain the second inclusion.
\epf
Observe that the equality holds for $n=1$ also, because $$\wdh^{(1)}=\wdhad=G(H)^\perp=H_{(1)}^\perp.$$
\subsection{Proof of Theorem \ref{iff:nilpotent}}
\bpf
If $\hb$ is nilpotent then  $H_{(n)}=H$, for some $n$.  Then, $\wdh^{(n)}=\comp$ by Equation \eqref{fiq}, and therefore $\whb$ is nilpotent. Idem for the converse with Equation \eqref{dual:fiq}.

\epf
\subsection{Burnside property for nilpotent ARN-hypergroups  }
Let $\hb$ be an ARN-hypergroup. Recall that we call $x_i\in \mtc B$ a  \emph{vanishing  element} if there is $\muj\in \wdb$ such that $\muj(x_i)=0$.
\bl\label{v:elements}
Let $\hb$ be an ARN-hypergroup and $\hs$ be a sub-hypergroup of $\hb$. Consider a standard element $x_i\in {\mtc B}$ such that  $[x_i]$ is a vanishing element in $H//\mtcs$. Then, $x_i$ is also a vanishing element in $\hb$.
\bpf
Recall that $H//\mtcs\simeq \lam_\mtcs H\lam_\mtcs$ as hypergroups, with $[x]\sent \lam_\mtcs x\lam_\mtcs$. 
If $\lam_\mtcs=\sum_{j\in \mtc I_\mtcs}F_j$ (Definition \ref{def:supp:j}), then $\lam_\mtcs H\lam_\mtcs=\oplus_{j\in \mtc I_\mtcs} \comp F_j$, as a subalgebra of $H$. Thus, the characters $\muj$, with $j\in \mtc I_\mtcs$, have distinct restriction to $\lam_\mtcs H\lam_\mtcs$ (since $\mu_i(F_j) = \delta_{i,j}$), and so cover all the characters of $\lam_\mtcs H\lam_\mtcs$ (since $\dim_{\comp}(\lam_\mtcs H\lam_\mtcs) = |I_\mtcs|$). If $[x_i]$ is a vanishing element in $H//\mtcs$, then $\lam_\mtcs x_i\lam_\mtcs$ is a vanishing element in $\lam_\mtcs H\lam_\mtcs$, thus by above, $x_i$ is a vanishing element in $H$.
\epf
\el
\bl \label{lem:subgroup}
Let $\hb$ be an ARN-hypergroup. A grouplike element $x_i \in \ghb$ is a constituent of  $x_jx_{j^*}$ if and only if $x_ix_j=\fp(x_i)x_j$. The set of all grouplike elements that are constituents of $x_jx_{j^*}$ form a subgroup of $\ghb$.
\el
\bpf
By Frobenius reciprocity (\ref{eq:FR}), $N_{j,j^*}^i = \frac{h_j}{h_i} N_{i,j}^j$, thus $x_i$ is a constituent of $x_j x_{j^*}$ if and only if $x_j$ is a constitutent of $x_ix_j$, if and only if (by Lemma \ref{perm}) $\frac{x_i}{d_i} \frac{x_j}{d_j} = \frac{x_j}{d_j}$, meaning that $x_ix_j=\fp(x_i)x_j$. The second assertion is a straightforward consequence of this. 
\epf
{\bf Proof of Theorem \ref{nilpotent:burnside}}
\vskip 0.2cm
By Theorem \ref{iff:nilpotent}, it is enough to prove the following:
\bp
A nilpotent dualizable ARN-hypergroup $\hb$ is Burnside.
\ep
\bpf
We will proceed by induction on the nilpotency class of $\hb$.

If  $\hb$ is of nilpotency class $1$, then $\hb$ is pointed and we are done.

Now, suppose that $\hb$ is nilpotent of class $n$ and that the result holds for the nilpotency classes less than $n$. We will show that any standard element $x_i\in \mtc B$ is either a vanishing or a grouplike element. Assume that $x_i\notin G(H)$. We will show that $x_i$ is a vanishing element.

Consider $[x_i]\in H//G(H)$. Since $H//G(H)\simeq \widehat{\wdhad}$ by (\ref{quot:gh}), it is nilpotent, with nilpotency class $n-1$, by Theorem \ref{iff:nilpotent}. By induction hypothesis, we know that $[x_i]$ is either a grouplike element or a vanishing element in $H//G(H)$. 

If $[x_i]$ is a vanishing element, then Lemma \ref{v:elements} shows that $x_i$ is a vanishing element.

Now, if $[x_i]$ is a grouplike element of $H//G(H)$ and by abelian assumption, $$\lam_{G(H)} = h_{[x_i]}\lam_{G(H)}x_i\lam_{G(H)}x_i^*\lam_{G(H)} = \lam_{G(H)}h_{[x_i]}x_ix_i^*\lam_{G(H)},$$ which implies that $h_{[x_i]}x_ix_i^*\in G(H)$.
Thus, $h_{[x_i]}x_ix_i^*=\sum_{g \in G_1}g$, for a subgroup $G_1 \subset G$ (by Lemma \ref{lem:subgroup}). Since $x_i$ is not a grouplike element in $H$, $G_1\neq \{1\}$. Then, there is $\muj$ such that $\muj(\blam_{G_1})=0$ where $\blam_{G_1}=\frac{1}{|G_1|}\big(\sum_{x\in G_1}x\big)$ is the integral of $G_1$. This implies that $\muj(x_ix_i^*)=0$. Thus $\muj(x_i)=0$, which means that $x_i$ is a vanishing element.
\epf
Now, since $H$ and $\wdh$ are simultaneously nilpotent, Theorem \ref{nilpotent:burnside} follows.
\section{Applications to fusion categories}\label{afc}
Let $\cc$ be a pivotal fusion category with a commutative Grothendieck ring.  As in \S \ref{introd}, we denote the set of isomorphism class representatives of simple objects of $\cc$ by $\irr(\cc):=\{X_1,\dots,X_m\}$  and let $\mtc I: =\{1,\dots,m\}$. Let also $d_i:=\dim(X_i)$  be the categorical dimension of $X_{i}$ for all $i\in \mtc I$. For any simple object $X_i\in \irr(\cc)$ denote by $x_i:=[X_i]$ the class of $X_i$ in the Grothendieck ring $K(\cc)$ of $\cc$.  

Let  $\czcc$ be the Drinfeld center of the fusion category  $\cc$. The forgetful functor $F:\czcc\ra\cc$  admits a right adjoint functor $R:\cc \ra \czcc$.  It is well-known that $A:=FR(\unu)$ has the structure of a central commutative algebra in $\cc$ (meaning that $R(\unu)$ is a commutative algebra in $\czcc$), where $\unu$ is the unit object of $\cc$. The vector spaces $$\cecc:= \hm_{\C}(\unu, A) \text{ and } \cfcc:=\hm_\cc(A, \unu)$$ are respectively called \emph{the space of central elements} and the \emph{space of class functions} of $\cc$. Recall \cite{scalg} that $\cfcc\simeq K(\cc)_\comp$, the complexification of the Grothendieck ring  $K(\cc)$ of $\cc$.

For any simple object $X_i$ of $\cc$ we denote by $\ch_i:=\mtr{ch}(X_i)\in \cfcc$ its associated character \cite[Definition 3.9]{scalg}. The central element space $\cecc$ has a basis of primitive orthogonal idempotents $E_i$ such that   $<\ch_i, E_j>=\delta_{i,j}d_i$ for all $1\leq i,j\leq m$.

As explained in \cite[Theorem 3.8]{scalg} the adjunction between $F$ and $R$  gives a canonical isomorphism of algebras
\beq\label{adjisom}
\cfcc \xra{\cong} \mtr{End}_{\czcc}(R(\unu)),
\eeq
where $\cfcc$ is the algebra of class functions on $\cc$.  Since $\czcc$ is also a fusion category we can write $R(\unu)=\bigoplus_{j \in \mathcal I}\mathcal C^j$ as a direct sum of simple objects in $\czcc$. Recall that $\cc^j$ are called \emph{conjugacy classes } for $\cc$. The above isomorphism also gives a canonical bijection between the set of primitive central idempotents  $\{F_j\}_{j \in \mathcal I}$ and the set of conjugacy classes $\{\cc^j\}$. 

For any fusion category $\cc$, by abuse of notation,  the \invisible{symmetric}RN-hypergroup $\hbz:=(K(\cc), \irrcc, \fp)$ is also denoted by $K(\cc)$. We denote by $\wkc$ its dual hypergroup.

Recall also from \cite{NG}  the notion of a universal group grading  of a fusion category $\cc$ and its  adjoint subcategory $\ccad$.
\br\label{hj}
 By  \cite[Equation (4.8)]{ccc-march}, it follows that for any spherical fusion category, $\dim(\cc^j)=\frac{\dim(\cc)}{n_j}$. See also \cite[Theorem 2.13]{O3} for a related statement. Then Equation \eqref{whj:eq} implies that 
\beq\label{whj}
\whj=\dim(\cc^j)
\eeq 
in the dual hypergroup $\wkc$.  In particular, for a weakly-integral fusion category $\cc$, by \cite[Proposition 8.27]{eno-annals}, $\dim(\cc^j)\in \mathbb Z$, and therefore $\wkc$ is $h$-integral.
\er
\bc \label{cor:BurCoProd}
Let $\cc$ be a a fusion category with a Burnside commutative Grothendieck ring and let $\wkc$ be its dual. Then 
\beq\label{fus:categ}
\prod_{j=1}^m\mu_j=\sum_{X_i\in \ccpt}\sgn(x_i)\wtf_i
\eeq
where $\sgn(x_i)$ is the determinant of the permutation matrix $L_{\nxi}$ on $K(\cc)$.
\ec
\bpf
Since  $K(\cc)$ is a Burnside hypergroup  the result follows from Corollary \ref{bsd:iff:pt:id}.
\epf
Note that the dual hypergroup $\wkc$ is denoted by $\wcfcc$ in \cite{b-blms}. Theorem 3.4 from the same paper implies that for any pivotal fusion category there is a canonical isomorphism of $\comp$-algebra 
\beq\label{al:isom}
\al:\wcfcc \ra \cecc,\;\muj\sent \frac{C_j}{\dim(\cc^j)},
\eeq
where $C_j \in \cecc$ is the \emph{conjugacy class sum} corresponding to the  conjugacy class $\mathcal C^j$, defined for a pivotal fusion category by Shimizu as ${\mtf}^{-1}(F_j)$. Here  $\lam\in \cfcc$ is  a cointegral of $\cc$ such that $\langle \lam, u\rangle =1$, see \cite[\S 5]{scalg}.  Recall also the \emph{Fourier transform} of $\cc$  associated to $\lambda$ is the linear map
\beq
\mtc F_{\lambda}:\cecc\ra \cfcc\;\;\text{given by}\;\;a \mapsto \lambda \lh \mtc S(a)
\eeq

By \cite[Lemma 4.1 and (4.7)]{ccc-march}, for any spherical fusion category $\cc$ over $\mathbb{C}$ with a commutative Grothendieck ring, we have 
$$
\langle F_i, C_j \rangle = \delta_{i,j} \dim(\cc^j).
$$

This relation implies that the set $(F_i)$ forms a dual basis to $(C_j / \dim(\cc^j))$ with respect to the bilinear form $\langle , \rangle$. Consequently, this provides an easier way for defining the conjugacy class sums $(C_j)$ in this case (see \cite[(4.11)]{ccc-march} for a formula).

\bc
For any weakly-integral fusion category $\cc$ with a commutative Grothendieck ring the following identity holds in $\cecc$:
\beq\label{cls:sums:prod}
\prod_{j=1}^m \frac{C_j}{\dim(\cc^j)}=\sum_{X_i\in \ccpt}\sgn(x_i)E_i
\eeq
where $E_i\in \cecc$ is the primitive central idempotent of $X_i\in \ccpt$.
\ec
\bpf
By Corollary \ref{hhj}, we can apply Corollary \ref{cor:BurCoProd}. The result follows by applying the canonical isomorphism $\alpha$.
\epf
For a spherical fusion category note the following:
\bp \label{grp:ccj}
Let $\cc$ be a spherical fusion category with a commutative Grothendieck ring $\kc$ and RN dual. Let $\muj\in \wkc$ be a linear character. Then $\muj$ is a grouplike element of $\wkc$ if and only if $\dim(\cc^j)=1$.
\ep
\bpf
By Lemma \ref{grouplike:nn} and the fact that $\wkc$ is normalized, $\muj$ is a grouplike element if and only if $\wh_j=1$. This is also equivalent to $\dim(\mathcal C^j)=1$, by Equation \eqref{whj}.
\epf
Now, suppose that $\cc$ is a pivotal fusion category and $\cd \subseteq \cc$ is a fusion subcategory. Following the notion of support defined in Definition \ref{def:supp:j}, there exists a subset $\mtc I_\cd \subseteq \mtc I$ such that 

$$
\lam_{\cd} = \sum_{j \in \mtc I_\cd} F_j,
$$

Please note that this set is called $\mtc L_\cd$ in \cite[\S 4.2]{ccc-march}. Additionally, \cite[Lemma 4.6]{ccc-march} implies \beq\label{jccad}
{\mtc I}_{ _{\ccad}}=\{j\in \mtc I\;|\; \dim(
\cc^j)=1\},
\eeq
for any fusion category with a commutative Grothendieck ring $\kc$.
\bp\label{wkc:v}
Let $\cc$ be a fusion category with a commutative Grothendieck ring. Then  $\wkc$ is Burnside if and only if:
\beq\label{dual:muj:gen}
\prod_{i=1}^m\nxi=\sum_{j \in \jccad}\sgn(\muj)F_{j}
\eeq
\ep
\bpf
It follows from Corollary \ref{dl:brsd:id} and Theorem \eqref{grpk:eq:hbad}. 
\epf
{\bf Proof of Theorem \ref{dual:v:prop}}
\bpf
It follows from Corollary \ref{d-b:sq:idm} and Theorem \eqref{grpk:eq:hbad}, together with 
\beqn
\sum_{j \in \jccad}F_{j}=\lam_{\ccad}=\frac{1}{\dim(\ccad)}\big(\sum_{x_i\in \ccad}d_ix_i\big),
\eeqn
where the last equality comes from Equation (\ref{eq:lamls}).
\epf

%
 
Here is a reformulation of Theorem \ref{dual:v:prop} for the case of a nilpotent finite group:
 
 \bt
For any finite nilpotent group $G$, the dual $\widehta{\text{ch}(G)}$  is Burnside and therefore 
\beq
\left(\prod_{x_{i}\in \irr(G)}  \nxi \right)^2 =\frac{|Z(G)|}{|G|}\left(\sum_{x_i\in \irr(G/Z(G))}d_ix_i \right).
\eeq
\et
\bpf
Recall that $ K(\text{Rep}(G)) $ is the character ring $ \text{ch}(G) $. The result that $ \widehat{\text{ch}(G)} $ is a RN hypergroup is established in \cite[Equation (3.12)]{b-blms}. Furthermore, \cite[Theorem B]{inw} states that if $ G $ is a nilpotent group, then $ \text{ch}(G) $ is dual-Burnside.
\epf
Theorem \ref{nilpotent:burnside} leads to the following conclusion:

\bc \label{cor:NilBur}
If $\cc$ is a nilpotent fusion category with a commutative Grothendieck ring and RN dual, then both $\kc$ and $\wkc$ are Burnside.
\ec

\br
A nilpotent fusion category is weakly-integral \cite{NG}, which means it is also pseudo-unitary and spherical \cite{EGNO15}. It remains an open question whether every pseudo-unitary fusion category can be given a unitary structure \cite[Remark 9.4.7]{EGNO15}. However, a commutative fusion ring with a unitary categorification does have a RN dual \cite{lpw}. Therefore, we believe that the assumption of a RN dual in Corollary \ref{cor:NilBur} could be omitted. Additionally, in the braided case, its Grothendieck ring possesses a RN dual according to \cite[Theorem 1.2]{b-blms}.
\er
\section{Premodular categories}\label{pmc}
Recall that a premodular category is defined as a braided spherical fusion category. In the rest of this section, we will assume that $\mathcal{C}$ is a premodular category.
 %
By \cite[Example 6.14]{scalg}  there is $\comp$-algebra map $\fq: {\cfcc}\ra \cecc$ given by the following formula:
\beq\label{sh}
\fq(\ch_i)=\sum_{j \in \mtc I}\frac{s_{ij}}{d_{j}}E_{j}.
\eeq
where $S=(s_{ij})$ is the $S$-matrix of $\cc$ and $(E_{j})$ are the primitive central idempotents of $\cecc$ as defined in the previous section.
 
As described in \cite[\S 4]{b-ant}, there is a map $ M: \mathcal{I} \to \mathcal{I} $ (referred to as the \emph{braided partition function}) such that if 
\beqn
\fq(F_{j})=\sum_{i \in \mtc B_{j}}E_{i},
\eeqn
then $ M(i) = j $. Let $ \mathcal{I}_2 \subseteq \mathcal{I} $ denote the set of all indices $ j $ for which $ \fq(F_j) \neq 0 $, meaning that $ \mathcal{B}_j $ is not an empty set. Given that $ \fq(1) = 1 $, we can thus partition the set of all isomorphism classes of simple objects $ \text{Irr}(\mathcal{C}) = \bigsqcup_{j \in \mathcal{I}_2} \mtc R_j $, where $ \mtc R_j = \{ [X_i] \;|\; i \in \mathcal{B}_j \} $. In other words, $\mtc R_j=\{[X_i]\;|\; M(i)=j\}$.
Thus, we obtain a unique function $ M: \mathcal{I} \to \mathcal{I}_2 $ with the property that $ E_i \fq(F_{M(i)}) \neq 0 $ for all $ i \in \mathcal{I} $.
 
The paper \cite{ntens} introduces the notion of cosets of a fusion category $\mathcal{C}$ with respect to a fusion subcategory $\mathcal{D}$: two simple objects $X,Y$ in $\mathcal{C}$ are in the same (right) coset if and only if there is a simple object $S$ in $\mathcal{D}$ such that $X$ is a constituent of $Y \otimes S$. By \cite[Theorem 4.10]{b-ant} two simple objects $X_i, X_{i'}$ of a pseudo-unitary  premodular fusion category $\cc$ are in the same coset with respect to $\cc':=\cz_2(\cc)$ if and only if $M(i)=M(i')$. In other words, $(\mtc R_j)_{j\in \mtc I_2}$ correspond exactly to these cosets.  Denote also ${\mtr R}_j:=\sum_{[X_i] \in \mtc R_j}d_i\ch_i\in \cfcc$ the regular part of their characters. By \cite[Equation (4.20)]{b-ant},
\beq\label{newr}
\dim({\mtr R}_j)=\dim({\cc'})\dim(\cc^j),\;\text{for all}\; j \in \mtc I_2.
\eeq
\bp \label{prop:dualbwrt}
Suppose that $\cc$ is a weakly-integral premodular  category such that $\cc'\subseteq \ccpt$ and acting freely on the set $\irrcc$.
Then for any $j\in \mtc I_2$ such that $\dim(\cc^{j})>1$ there is an $x_{i}$ such that $\muj(x_{i})\neq 0$. 
\ep
\bpf 

First of all, $\kc$ is $h$-integral because as a fusion ring, $h_i=1$ for all $i$. If $\cc$ is weakly integral, then the dual $\wkc$ is a rational RN hypergroup, as shown by \cite[Equation (4.3)]{b-blms}. Therefore, we can attempt to apply Theorem \ref{dual:burnside} for $\hbz = \kc$. To do this, we need to ensure that the additional condition \eqref{duality:alg} holds for all $i \in \mtc I$, specifically that $\dim(\cc^j)\frac{|\alij|^2}{d_i^2} \in \mathbb A$ for all $i \in \mtc I$, where $\alij:= \mu_j(x_i)$. It was demonstrated in \cite[Lemma 4.2]{b-ant} that

\beq\label{sym:fus}
\frac{\al_{ _{iM(i')}}}{d_i}=\frac{s_{ _{ii'}}}{d_id_{i'}}=\frac{\al_{ _{i'M(i)}}}{d_{i'}}.
\eeq
for all $i,i' \in \mtc I$. Fix an index $i_j\in \mtc I$ such that $M(i_j)=j \in \mtc I_2$, and take $i' = i_j$. It follows that
\beq\label{ali:s3}
\alij=\frac{d_i}{d_{i_j}}\al_{i_jM(i)}.
\eeq
Therefore 
\beq
\dim(\cc^j)\frac{|\alij|^2}{d_i^2}\numeq{\ref{ali:s3}}\dim(\cc^j)\frac{|d_i|^2}{|d_{i_j}|^2}\frac{\mid\al_{i_jM(i)}\mid^2}{d_i^2}
\eeq
But $\cc$ is spherical, so $d_i=\overline{d_i}$ and therefore $|d_i|^2=d_i^2$. Thus 
\beq
\dim(\cc^j)\frac{|\alij|^2}{d_i^2}=\frac{\dim(\cc^j)}{d_{i_j}^2}\mid\al_{i_jM(i)}\mid^2
\eeq

But $\cc'\subseteq \ccpt$, acting freely on $\irrcc$, thus $|\mtc R_j| = \dim(\mtc C')$, and so by Equation \eqref{newr}, $\dim(\cc^j)=\frac{\dim(R_j)}{\dim(\cc')}=d_{i_j}^2$, and the proof follows since $\al_{i_jM(i)}\in \mathbb A$.
\epf

Roughly speaking, Proposition \ref{prop:dualbwrt} can be intuitively understood as follows: a weakly integral premodular fusion category $\mathcal{C}$ is dual-Burnside with respect to $\mathcal{C}'$.

\subsection{Modular category case}\label{mtc}
In this section, let $\cc$ denote a modular fusion category. We will utilize the results from the previous section to $K(\mathcal{C})$ and prove Theorems \ref{mtc:even:order} and \ref{first:div}. Recall the definition of $ \mtc I_{ _\cd} $ as provided in Definition \eqref{def:supp:j}.

By composing Drinfeld's map $\fq:\cfcc\ra \cecc$ with the inverse of the isomorphism $\al:\cecc\ra \wcfcc$ from Equation \eqref{al:isom}, we obtain an algebra isomorphism $\tfq:\cfcc\ra\wcfcc$. This can also be expressed as an algebra isomorphism $\tfq:\kc\ra\wkc$. It was shown in \cite[Remark 4.2]{b-blms} that this is, in fact, an isomorphism of normalized hypergroups: $\tfq:\overline{\kc}\ra \wkc$, where $\overline{\kc}$ denotes the normalized version of $\kc$, so that $\tfq(\frac{x_i}{d_i}) = \mu_i$.

\br\label{bij:ij}
Drinfeld's map $\fq$ is an algebra isomorphism since the S-matrix is invertible. It establishes a canonical bijection between the index sets of $\{F_i \}$ and $\{E_i \}$ which is given by $\fq(F_i)= E_i$, for all $i\in \mtc I$. 
\er
\bc\label{mtc:case}
Let $\cc$ be a modular fusion category. Then $K(\cc)$ is Burnside if and only if it is dual-Burnside.
\ec
\bpf
As mentioned above  $\tfq:\overline{\kc}\ra \wkc$  is an isomorphism of normalized hypergroup. Therefore $ \overline{\kc}$ (and thus $\kc$) is Burnside if and only if $\wkc$ is also Burnside.
\epf
\bl\label{j:ccad}
Let $\cc$ be a modular fusion category. Then with the above notations,
\beq
X_i\in \irr(\ccpt) \iff i\in {\mtc I}_\ccad
\eeq
\el
\bpf
Apply Theorem \ref{grpk:eq:hbad} with the isomorphism of normalized hypergroups $\tfq$.
\epf
\bt 
In any modular fusion category $\cc$ with $\kc$ Burnside the following identity holds:
\beq\label{gen:case}
\prod_{i \in \mtc I} \frac{x_i}{d_i}
=\sum_{j\in {\mtc I}_\ccad}\sgn(x_j)F_j
\eeq
\et
\bpf
Since $K(\cc)$ is Burnside, we can apply $\tfq^{-1}$ to  Equation \eqref{fus:categ},
but $\tfq(\nxi)=\mui$, so we get
$$
\prod_{i \in \mtc I} \nxi=\sum_{X_i\in \ccpt}\sgn(x_i)F_i
$$
and Lemma \ref{j:ccad} finishes the proof.
\epf
{\bf Proof of Theorem \ref{mtc:even:order}}
\bpf
Squaring Equation \eqref{gen:case} we obtain that 
$\big(\prod_{i \in \mtc I} \frac{x_i}{d_i}\big)^{2}=\sum_{j\in {\mtc I}_\ccad}F_j = \lam_{\ccad}$ by the definition of $\jccad$. The result follows by Equation (\ref{eq:lamls}).
\epf
Recall that the above theorem holds for weakly-integral modular categories since their Grothendieck rings are Burnside by \cite[Appendix]{gnn} or \cite[Theorem 2]{b-galois}.
\vskip 0.7cm
{\bf Proof of Corollary \ref{mtc:odd:order}}

\bpf
In this case every invertible element $x_j$ has odd order since this order divides $|\gcc|$. Thus every $\sgn(x_j)=1$ for any invertible object $X_j\in \ccpt$.
\epf
\bp\label{ccad:d:burn}
Let $\cc$ be a fusion category such that $\kc$ is commutative and the dual $\wkc$ is Burnside. Then 
$$
\frac{(\prod_{i \in \mtc I} d_i)^2}{\fp(\ccad)}\in \mathbb Z.
$$
\ep
\bpf 
Let $P_\cc$ be $\prod_{i \in \mtc I} d_i$.  Suppose that $\big(\prod_{i \in \mtc I} x_i\big)^2=\sumitom N_ix_i$ with $N_i\in \mathbb Z_{\geq 0}$. Then equalizing coefficients in Equation \eqref{dual:muj:gen:squared}, we obtain the following $\frac{N_i}{P_\cc^2}=\frac{d_i}{\fp(\ccad)}$. Take $i=1$ then  $N_1\fp(\ccad)=P_\cc^2$ which proves the divisibility.
\epf
\vskip 0.3cm
{\bf Proof of Theorem \ref{dual:burns:div}.}
\bpf
By \cite[Theorem 3.10]{NG} $d_i^2\in \mathbb Z$ for all $i$ since $\cc$ is weakly-integral. Note that since $\wkc$ is Burnside the first divisibility result follows from Proposition \ref{ccad:d:burn}. Thus  $\fpccad$ is an integer (as a rational algebraic integer). Moreover, if $\cc$ is nilpotent  then \cite[Corollary 5.3]{NG} implies $d_i^2\mid \fpccad$. Thus  $\mtc V(\ccad)=\bigcup_{i \in \mtc I} \mtc V(d_i^{2})$.
\epf
\br
In the case of an integral fusion category Equation \eqref{nilp:ccad} can be written as $$\mtc V(\ccad)=\bigcup_{i \in \mtc I} \mtc V(d_i).
$$
\er
Recall that a weakly-integral fusion category with a commutative Grothendieck ring is Burnside by \cite[Theorem 2]{b-galois}. Moreover, a Burnside modular fusion category is dual-Burnside by Corollary \ref{mtc:case}. So a weakly-integral modular fusion category is dual-Burnside.
\vskip 0.3cm
{\bf Proof of Theorem \ref{first:div} }
\bpf 
 The first divisibility  follows from Proposition \eqref{ccad:d:burn}. Now, let $\cc$ be a weakly-integral modular fusion category. By \cite[Proposition 8.27]{eno-annals} we have that $\fp(\ccad)$ is an integer. By \cite[Theorem 3.10]{NG} $d_i^2\in \mathbb Z$ for all $i$. The second item follows since for a modular category we have $\fp(\cc)=\fp(\ccad)\fp(\ccpt)$, because the group associated to $\ccpt$ is the universal grading group by \cite[Lemma 8.22.5]{EGNO15}. Thus $\mtc V(\cc)=\mtc V(\ccad)\cup \mtc V(\ccpt)$. On the other hand, from the first item, $\mtc V(\ccad)\subseteq \bigcup_{i \in \mtc I} \mtc V(d_i^{{2}})$ which proves the inclusion $\mtc V(\cc)\subseteq \mtc V(\ccpt)\cup \big(\bigcup_{i \in \mtc I} \mtc V(d_i^{{2}})\big)
$.

Conversely, $\mtc V(\ccpt)\cup \big(\bigcup_{i \in \mtc I} \mtc V(d_i^2)\big)\subseteq \mtc V(\cc)$ since $\fp(\ccpt)\mid \fpcc$, and by \cite[Proposition 8.14.6]{EGNO15}, $d_i^{{2}}\mid \fp(\cc)$ for all $i \in \mtc I$.
\epf
{{\bf Proof of Theorem \ref{gen:o-yu}}
\bpf
Since $\gcd(d_i^2, d) = 1$ and $d$ is square-free, it follows from Equation \eqref{prime:set} that $d$ divides $\fp(\ccpt)$. Now, consider a pointed fusion subcategory $\cd \subseteq \ccpt$ such that $\fp(\cd) = d$. The existence of such a subcategory is straightforward to demonstrate because the associated groups are abelian, given that the category is braided. According to \cite[Theorem 8.21.5]{EGNO15}, with $\cc' = \vect$, we have $\fp(\cd') = m$, but $\gcd(d, m) = 1$. Consequently, by Lagrange's theorem \cite[Theorem 7.17.6]{EGNO15}, $\cd \cap \cd' = \vect$, which implies that $\cd$ is non-degenerate, as stated in \cite[Corollary 8.20.10]{EGNO15}. Thus, by \cite[Theorem 8.21.4]{EGNO15}, $\cc \simeq \cd \boxtimes \cd'$, and $\cd'$ is also non-degenerate.
\epf
\br
Note that Theorem \ref{gen:o-yu} improves \cite[Theorem 4.5]{o-yu} since the authors also assume the existence of a Tannakian subcategory $\ce=\rep(G)\subseteq \cc$ such that $\cc^0_G\simeq \cc(\mathbb Z_d, q)\boxtimes \cca$ for some non-degenerate braided fusion category $\cca$.
\er 
{\bf Proof of Corollary \ref{ccpt}:}
\bpf
Immediate from Theorem \ref{gen:o-yu}, we just need to check that $gcd(d_i^2,d)=1$, but it is clear since $d_i^2$ divides $\fp(\cc)$ and $d$ is a factor of the square-free part.
\epf
}


\bc\label{cor:maxd}
An integral modular fusion category $\cc$ can always be decomposed into $\cd \boxtimes \cd'$, where both $\cd$ and $\cd'$ are modular, $\cd$ is pointed, and $\fp(\cd)$ is the maximal square-free part of $\fp(\cc)$.
\ec

This means that the classification of integral modular fusion categories $\cc$ simplifies to those where $\fp(\cc)$ lacks a square-free part, up to a Deligne tensor product with a pointed modular fusion category.

\bc\label{p:square}
For an integral perfect modular fusion category $\cc$, $\fp(\cc)$ has no square-free part. In other words, if a prime $p$ divides $\fp(\cc)$, then $p^2$ also divides $\fp(\cc)$.
\ec

Corollary \ref{p:square} obviously extends to the weakly-integral case due to the following:

\bp\label{prop:wint}
A weakly-integral perfect modular fusion category $\cc$ is integral.
\ep

\bpf
Since $\ccpt$ is trivial then $\ccad=\cc$, because $(\ccpt)' = \ccad$ by \cite[Corollary 8.22.8]{EGNO15}, see also Remark \ref{rem:modperf}. Moreover, \cite[Exercise 9.6.12]{EGNO15} asserts that the adjoint subcategory of a weakly-integral fusion category is integral.
\epf

\bc
If $\cc$ is a perfect, integral, even-dimensional modular fusion category, then $4$ divides $\fp(\cc)$.
\ec
\subsection{On Statements \eqref{5.4} and \eqref{5.5} and Conjecture \eqref{odd-j81}}
%
%
In \cite[Theorem 5.6]{lpr2} it was shown that Statements \eqref{5.4} and \eqref{5.5} are equivalent. 
We will show that any of these statements  implies Conjecture \ref{odd-j81} formulated in  \cite{j-plav-odd}.
\bp\label{main:link}
 The existence of a perfect odd dimensional  modular fusion category  implies the existence of a non-pointed simple {integral} modular fusion category.
 \ep
\bpf
Let $\cc$ be a perfect modular fusion category with an odd dimension. By Proposition \ref{prop:wint}, $\cc$ must be integral.

Firstly, since $\cc$ is perfect, it cannot have any (non-trivial) pointed subcategories. Additionally, it cannot have (non-trivial) symmetric subcategories. By Deligne's theorem, symmetric subcategories would be equivalent to $\rep(G,z)$ for a perfect group $G$, implying that $|G|$ is even by Feit–Thompson theorem. This contradicts the odd dimension of $\cc$.

However, due to the finite number of isomorphism classes of simple objects, $\cc$ must contain at least one minimal subcategory. Let us call this minimal subcategory $\cd$. If $\cd$ is itself a modular fusion category, then $\cd$ serves as an example of a non-pointed, simple, integral modular fusion category.

If $\cd$ is not modular, then $\cd' \cap \cd$ is non-trivial according to \cite[Corollary 8.20.10]{EGNO15}, yet it is symmetric by design. But as discussed earlier, such symmetric subcategories cannot exist, leading to a contradiction.
\epf
Recall that in \cite{j-plav-odd} the authors have also shown that Conjecture \ref{odd-j81} is equivalent to either of the following:
\begin{conj}\label{odd-j82}
Odd-dimensional fusion categories are solvable.
\end{conj}
\begin{conj}\label{odd-j83}
Odd-dimensional modular fusion categories are solvable.
\end{conj}  

\section{Applications and extra results} \label{sec:AppExtra}
\subsection{Perfect Drinfeld center} \label{sub:perfmod}
We will characterize the perfect fusion categories with a perfect Drinfeld center. That will provide a large class of perfect modular fusion categories.

%
%
%

\bl \label{lem:Forg1} Let $\mathcal{C}$ be a fusion category. Let $F: \mathcal{Z}(\mathcal{C}) \to \mathcal{C}$ be the forgetful functor. Let $F_1:\czcc_{\mtr{pt}}\ra \ccpt$ be the induced group homomorphism. Let $G$ be the  universal grading group of $\mathcal{C}$. Then, 
\begin{enumerate}
\item $\ker F_1$ is isomorphic to $\widehat{G_{ab}}$, with $G_{ab}:= G/G'$ the abelianization of $G$,
\item $G$ is perfect if and only if $\ker F_1$ is trivial (i.e. $F(Z)=1$ implies $Z=1$).
\end{enumerate}
\el
\bpf
By definition of the center $\mathcal{Z}(\mathcal{C})$, see \cite[Definition 7.13.1]{EGNO15}, $F(Z)=1$ if and only if $Z=(1,\gamma)$ with $\gamma \in \Aut_{\otimes}(\id_{\mathcal{C}})$. So $\Aut_{\otimes}(\id_{\mathcal{C}})$ is isomorphic to $\ker F_1$. Now, \cite[Proposition 3.9]{NG} states that $\Aut_{\otimes}(\id_{\mathcal{C}})$ is isomorphic to $\widehat{G_{ab}}$. In particular, $\ker F_1$ is trivial if and only if $G'=G$ (i.e. $G$ is perfect).
\epf

Note that Lemma \ref{lem:Forg1} covers \cite[Lemma 2.1]{DNV15}.

\bp \label{prop:perfZ}
Let $\mathcal{C}$ be a perfect fusion category. Then $\mathcal{Z}(\mathcal{C})$ is perfect if and only if the universal grading group of $\mathcal{C}$ is perfect.
\ep 
\bpf
Let $Z$ be a simple object of $\mathcal{Z}(\mathcal{C})$ with $\fp(Z)=1$. The forgetful functor $F: \mathcal{Z}(\mathcal{C}) \to \mathcal{C}$ is a tensor functor, so $\fp(F(Z))=1$, but $\mathcal{C}$ is perfect, so $F(Z)=1$. The result follows by Lemma \ref{lem:Forg1}.
\epf

\bc \label{cor:perfZbr}
Let $\mathcal{C}$ be a braided fusion category. Then $\mathcal{Z}(\mathcal{C})$ is perfect if and only if $\mathcal{C}$ is perfect with a trivial universal grading group.
\ec 
\bpf
If $\mathcal{Z}(\mathcal{C})$ is perfect, then $\mathcal{C}$ is also perfect because it can be embedded into $\mathcal{Z}(\mathcal{C})$ as a braided fusion category. According to Proposition \ref{prop:perfZ}, the universal grading group of $\mathcal{C}$ is perfect. However, since $\mathcal{C}$ is braided, its Grothendieck ring is commutative, which means its universal grading group must be abelian due to the faithful grading. But a perfect abelian group is trivial. Conversely, the result follows directly from Proposition \ref{prop:perfZ}, as the trivial group is perfect.
\epf

Proposition \ref{prop:perfZ} and Corollary \ref{cor:perfZbr} provide a complete characterization of when the Drinfeld center $\mathcal{Z}(\mathcal{C})$ is perfect assuming the fusion category $\mathcal{C}$ is perfect or braided. However, it is important to note that there are fusion categories that are neither perfect nor braided, yet their Drinfeld centers are perfect. An example of this is the fusion category $\mathcal{C} = \mathrm{Vec}(G)$, where $G$ is a non-abelian finite simple group. This fusion category is Morita equivalent to $\rep(G)$, as illustrated in \cite[Example 7.12.19]{EGNO15}. As a result, their Drinfeld centers are braided equivalent, as shown in \cite[Theorem 8.12.3]{EGNO15}. Nonetheless, Corollary \ref{cor:simpleDG} below confirms that $\mathcal{Z}(\rep(G))$ is perfect in this case.

\bl \label{lem:centergrp}
The universal grading group of $\rep(G)$ is isomorphic to the center $Z(G)$ of the finite group $G$.
\el
\bpf
Immediate by Theorem \ref{univ:grd}, because a grouplike element $\mu_j$ in the dual hypergroup of the Grothendieck ring of $\rep(G)$ corresponds to a column of the character table of $G$ with squared norm $n_j=|G|$, by Lemma \ref{max:formal:codeg}, so to a conjugacy class of size $|G|/n_j = 1$, thus to a central element. 
\epf

\bc \label{cor:perfDG}
Let $G$ be a finite group. The Drinfeld center $\mathcal{Z}(\rep(G))$ is perfect if and only if $G$ is perfect with a trivial center $Z(G)$.
\ec 
\bpf
Immediate by Corollary \ref{cor:perfZbr} and  Lemma \ref{lem:centergrp}, because $\rep(G)$ is perfect if and only if $G$ is perfect.
\epf
\br \label{rk:CV} Corollary \ref{cor:perfDG} can also be proven using the fact (as referenced in \cite[\S 8.5]{EGNO15}, even when $G$ is infinite) that the simple objects of $\mathcal{Z}(\mathrm{Vec}(G))$ correspond to pairs $(C, V)$. Here, $C$ is a finite conjugacy class in $G$, and $V$ is an irreducible finite-dimensional representation of the centralizer of an element $g$ in $C$. The Frobenius-Perron dimension of the object associated with a pair $(C, V)$ is $|C| \cdot \dim_{\mathbb{C}}(V)$. But, as noted earlier, $\mathcal{Z}(\mathrm{Vec}(G))$ and $\mathcal{Z}(\rep(G))$ are braided equivalent if $G$ is finite.
\er

\bc \label{cor:simpleDG}
Let $G$ be a non-abelian finite simple group. Then $\mathcal{Z}(\rep(G))$ is perfect.
\ec
\bpf Immediate by Corollary \ref{cor:perfDG}, because a non-abelian finite simple group is perfect with a trivial center.
\epf

By Corollary \ref{cor:simpleDG}, $\mathcal{Z}(\rep(A_5))$ is a perfect integral modular fusion category of $\fp$ $60^2=3600$. Using Remark \ref{rk:CV} and GAP, it is of rank $22$ and type $$[ [ 1, 1 ], [ 3, 2 ], [ 4, 1 ], [ 5, 1 ], [ 12, 10 ], [ 15, 4 ], [ 20, 3 ] ].$$

\begin{question} \label{q:22}
Is there a non-trivial perfect integral modular fusion category of rank less than $22$, or $\fp$ of less than $3600$?
\end{question}

Recall that a Hopf algebra $A$ is called  \emph{perfect} if $\rep(A)$ is perfect, i.e. $G(A^*)$ is a trivial group; that a finite dimensional semisimple Hopf algebra $A$ (over $\mathbb{C}$) is \emph{factorizable} if and only if $\rep(A)$ is modular; that $\mathcal{Z}(\rep(A)) = \rep(D(A))$ where the Hopf algebra $D(A)$ is the \emph{Drinfeld double} of $A$. So, for any non-abelian finite simple group $G$ then $D(G)$ is a finite dimensional semisimple factorizable perfect Hopf algebra (over $\mathbb{C}$) of dimension $|G|^2$. Let us clarify Question \ref{q:22} for the case of Hopf algebras:

\begin{question} \label{q:22:hopf}
Is there a non-trivial finite-dimensional semisimple factorizable perfect Hopf algebra over $\mathbb{C}$ that has less than $22$ irreducible representations (up to isomorphism), or a dimension of less than $3600$?
\end{question}

\bl \label{lem:A*toA}
Let $A$ be a semisimple factorizable Hopf algebra. If $A^*$ is perfect then so is $A$.
\el
\bpf
By \cite[Theorem 2.3(b)]{schfact} there is an isomorphism $G(A^*)\simeq G(A)\cap Z(A)$. 
If $G(A)=\{1\}$ then clearly $G(A^*)=\{1\}$.  
\epf
Please note that the converse of Lemma \ref{lem:A*toA} does not hold. Specifically, for any centerless finite perfect group $ G $, the Drinfeld double $ D(G) $ is both factorizable and perfect, as established by Corollary \ref{cor:perfDG}. However, according to \cite[\S 7.12-7.16]{EGNO15}, the fusion category $\mathrm{Rep}(D(G)^*)$ is equivalent to $\mathrm{Vec}(G) \boxtimes \mathrm{Rep}(G)$. Consequently, $ D(G)^* $ is not perfect. Hopf algebras $A$ such that both $A$ and $A^*$ are perfect are called \emph{biperfect} Hopf algebras. The smallest known example of biperfect Hopf algebra (let us call it $H$) was described in \cite{egg-biperf} as a bicrossed product. This comes from an exact factorization $M_{24}=G_1G_2$ of the Mathieu group  of degree $24$. Here $G_1=\mtr{PSL}(2,23)$ and $G_2=(C_2)^4\rtimes A_7$. Its dimension is $\dimk(H)=|M_{24}|=2^{10}\times 3^3\times 5\times7\times11\times 23=244823040$. By \cite{rad-min}, $A$ is biperfect if and only if $D(A)^*$ is perfect (so biperfect, by Lemma \ref{lem:A*toA}). Thus $D(H)$ is also a biperfect Hopf algebra. 
Note that the authors of \cite{egg-biperf} suspect $H$ to be the smallest example of a biperfect Hopf algebra, and $M_{24}$ may also be the only finite simple group with a factorization that produces a biperfect Hopf algebra.

\subsection{Burnside and integrality properties} \label{sub:BurnInt}
\br \label{rem:uniper} According to Theorem \ref{univ:grd}, the dual hypergroup of a commutative fusion ring is perfect if and only if the universal grading group of the fusion ring is trivial.
\er
\br \label{rem:biper} The Grothendieck ring of a modular fusion category is perfect if and only if its dual is perfect (since it is isomorphic to its dual as normalized hypergroup, see \S \ref{mtc}).
\er
\br \label{rem:modperf} By Remarks \ref{rem:uniper} and \ref{rem:biper}, a modular fusion category $\cc$ is perfect if and only if its universal grading group is trivial (i.e. $\cc_{ad} = \cc$). Alternatively, recall that $(\ccpt)' = \ccad$ by \cite[Corollary 8.22.8]{EGNO15}.
\er
{\bf Proof of Theorem \ref{db:integral:hypgs}}
\bpf
We utilize the Galois action on $\widehat{\mathcal B}$ as described in Lemma \ref{def:sghstar}. Let $\mu_1$ be the $\fp$ character. According to Proposition \ref{sg:nk:eq:ntauk}, any Galois conjugate of $\mu_1$, denoted by $\mu_j = \sigma \cdot \mu_1$, has a formal codegree given by $n_j = \sigma(\fp(H, \mathcal{B}))$. However, because $\mu_j = \sigma \cdot \mu_1$ is also a non-vanishing character, it must be a grouplike element due to the definition of a dual-Burnside hypergroup (Definition \ref{def:d-B}). Therefore, its formal codegree must also be $\fp(H, \mathcal{B})$, by Corollary \ref{frmlcdg:grplike}. Thus, $\sigma(\fp(H, \mathcal{B})) = \fp(H, \mathcal{B})$ for every $\sigma$ in the Galois group. It follows that $\fp(H, \mathcal{B})$ is a rational number, completing the proof.
\epf
As a consequence of Theorem \ref{db:integral:hypgs}, every commutative dual-Burnside fusion ring is weakly integral. However, the converse is not true, even for unitary integral fusion categories such as $\rep(A_7)$; refer to \S \ref{sub:dualburn} for more examples. The following result is partially attributed to Andrew Schopieray; see \cite{AndrewMO}.
\bc \label{cor:andrew}
Let $\hbz$ be rational abelian normalizable hypergroup. If the dual  $\whbz$ is perfect and Burnside then $\hbz$ is integral, i.e $d_i\in \mathbb Q$.
\ec
\bpf
The dual $\whbz$ is perfect, meaning that $\mu_1$ is the only grouplike linear character. As in the proof of Theorem \ref{db:integral:hypgs}, for any $\sg \in \gal( \overline{\mathbb Q}/\mathbb Q)$, $\sg.\mu_1$ is also non-vanishing and thus grouplike by the dual-Burnside assumption. Consequently, $\sg.\mu_1 = \mu_1$, and $\sg(d_i) = d_i$ for all $\sg$, which implies that $d_i \in \mathbb{Q}$.
\epf
\bc \label{cor:andrewb}
A commutative dual-Burnside fusion ring with a trivial universal grading group is integral.
\ec
\bpf
It follows from Corollary \ref{cor:andrew} and Remark \ref{rem:uniper}.
\epf
Recall Remark \ref{rem:uni} about a trivial universal grading group.  
\bc \label{cor:andrew2}
A fusion category with a dual-Burnside commutative Grothendieck ring and a trivial universal grading group is integral.
\ec
\bpf
Immediate from Corollary \ref{cor:andrewb}.
\epf

\textbf{Proof of Theorem \ref{thm:perfmodburn}}
\begin{proof}
Recall that a modular fusion category is Burnside if and only if it is dual-Burnside (Corollary \ref{mtc:case}). Now, by Remark \ref{rem:modperf}, the universal grading group is trivial, so by Corollary \ref{cor:andrew2}, it must be integral. The converse follows from \cite[Theorem 6.1]{gnn}.
\end{proof} 
\br
According to Theorem \ref{thm:perfmodburn}, a simple non-integral modular fusion category cannot be (dual-)Burnside. Many modular fusion categories of Lie type, known as Verlinde categories, belong to this class, as do the one described in \cite[Theorem 1]{schop-np} and the Drinfeld center of the Extended-Haagerup fusion categories. Consequently, none of these are (dual-)Burnside. Therefore, for the Grothendieck rings in this class, we can still infer the existence of a non-group-like and non-vanishing linear character (and basic element).
\er

\subsection{Near-group modular fusion categories} \label{sub:nearmod}
Let $G$ be a finite abelian group and $m$ be a non-negative integer. Let $K(G,m)$ be the fusion ring  with basis $\mathcal{B} = G\cup \{\ro\}$ and fusion rules: $$\ro^2=\sum_{g \in G}g+m\ro \text{ and } g\ro=\ro g=\ro,$$ for all $g\in G$. A fusion category whose Grothendieck ring is $K(G,m)$ is called \emph{near-group} (or also \emph{Tambara-Yamagami} when $m=0$), see for example \cite{Izu17}.
\bp \label{prop:modnear}
There is no modular fusion category of Grothendieck ring $K(G,m)$ if $G$ is non-trivial and $m > 0$.
\ep
\bpf
Let $\mathcal{C}$ be a modular fusion category of Grothendieck ring $K(G,m)$. Given that $\ro^2 = m\ro + \sum_{g \in G} g$, it follows that $\fp(\ro)^2 = m\fp(\ro) + |G|$. Let $x_{+} > 0$ and $x_{-} < 0$ be the solutions to the equation $x^2 - mx - |G| = 0$. Thus, $\fp(\ro) = x_{+}$.
Now, consider $\eta$, a linear character of $R$. For all $g \in G$, we have $\eta(\ro) = \eta(\ro g) = \eta(\ro) \eta(g)$. Therefore, $\eta(\ro)$ is nonzero if and only if $\eta(g) = 1$ for all $g \in G$. This condition is equivalent to solving $\eta(\ro)^2 - m\eta(\ro) - |G| = 0$, which implies $\eta(\ro) = x_{\pm}$. We denote such characters as $\psi_{\pm}$.
It follows that $\psi_+ = \fp$ is the unit of $\wdb$. Let $\eps$ be the trivial character of $G$. We identify any $\eta$ in $\widehat{G} \setminus \{\eps\}$ (non-empty since $G$ is non-trivial) with the corresponding element in $\wdb$ that vanishes on $\ro$.
Let $q = -\frac{x_{-}}{x_{+}}$, which is positive. The multiplication on $\widehat{R}$ is defined by:
\begin{align*}
\eta\star \psi_{-}&=\eta, \\
\eta\star \eta'&=(1-\delta_{\eta', \eta^*})\eta\eta' + \delta_{\eta', \eta^*}(\frac{q}{1+q} \psi_{+} + \frac{1}{1+q}\psi_{-}), \\
\psi_{-}  \star \psi_{-} &=q\psi_++(1-q)\psi_-.
\end{align*}
To verify this, consider the evaluations at $\frac{g}{\fp(g)} = g$ and $\frac{\rho}{\fp(\rho)} = \frac{\rho}{x_{+}}$. Note that $\eta$ is not grouplike, as defined in Definition \ref{grplike:def}. Additionally, $\psi_{-}$ is also not grouplike because $1 - q > 0$ (implying that $R$ has an RN-dual). This follows from the inequality $x_{+}^2 = mx_{+} + |G| > mx_{-} + |G| = x_{-}^2$, given $m > 0$. Consequently, $G(\widehat{R}) = \{\fp\}$. Now, as discussed in \S \ref{mtc}, the normalization $\overline{R}$ of $R$ forms a selfdual hypergroup, assuming $\mathcal{C}$ is modular. Drinfeld's map $\tfq: \overline{R} \to \widehat{R}$ establishes this isomorphism of normalized hypergroups. It induces an isomorphism from $G(\overline{R})$ to $G(\widehat{R})$, contradicting $|G(\widehat{R})| = 1 < |G| = |G(\overline{R})|$, given $G$ non-trivial.
\epf
\bp \label{prop:modnear2}
Let $\mathcal{C}$ be a modular fusion category of Grothendieck ring $K(G,0)$. Then $G\simeq C_1$ or  $C_2$. 
\ep
\bpf
Assume that $G$ is non-trivial. Following the proof of Proposition \ref{prop:modnear}, $G(\widehat{R})$ is isomorphic to $G(\overline{R})=G$, by Drinfeld's map. And, if $m=0$ then $q=1$ and $\psi_-^2=\fp$. Therefore $G(\widehat{R}) = \{ \psi_+, \psi_-\}$, and $|G|=2$. 
\epf
\br \label{rk:ising}
Ising categories are the only modular fusion categories whose Grothendieck ring are $K(C_2, 0)$, see \cite[Appendix B]{dgno2}. 
\er
\bc \label{cor:nearmod}
A modular fusion category is near-group if and only if its Grothendieck ring is $K(G,m)$ with $(G,m) = (C_1,0), (C_1,1), (C_2,0)$.
\ec
\begin{proof}
According to Proposition \ref{prop:modnear}, either $G$ is trivial or $m = 0$. If $G$ is trivial, the category has rank 2, and the conclusion follows from \cite{OstR2}. On the other hand, if $G$ is non-trivial and $m = 0$, the result follows from Proposition \ref{prop:modnear2} and Remark \ref{rk:ising}.
\end{proof}
\br
Note that modular (generalized) Tambara-Yamagami categories were classified before, see \cite[Lemma 5.3 and Theorem 5.4]{NataleFaitful},  \cite{Thornton} and \cite{Siehler}.
\er

\section{Concrete examples} \label{sec:exacrit}
This section focuses on providing concrete examples of fusion rings and categories to which the main results of the paper apply as categorification criteria, along with counter-examples demonstrating the limits of these results.

\subsection{Burnside property} \label{sub:burn}
Recall that a fusion ring is called \emph{Burnside} if, for every basis element, its fusion matrix has a norm of $1$ if and only if its determinant is nonzero. In other words, being grouplike is equivalent to being invertible (see Proposition \ref{grouplike:set}). William Burnside originally proved that the Grothendieck ring of $\Rep(G)$ is Burnside for every finite group $G$. This result was expanded in \cite{b-galois} to include all weakly integral fusion categories with a commutative Grothendieck ring, but first within the context of modular categories in \cite[Appendix]{gnn}. In Theorem \ref{burnside-fr}, it extends to every commutative fusion ring with an $h$-integral dual, and further to hypergroups in Theorem \ref{burnside}.

\br
By Theorems \ref{hbz:hyp} and \ref{dual:univ:grd}, along with the fact that $\rep(G)_{pt} \simeq \rep(G/G')$, Burnside's result can be restated as follows: for any finite group $G$,
$$
\left(\prod_{j \in \mathcal{I}} \frac{C_j}{|\mathcal C^j|}\right)^2=\frac{1}{|G'|}\sum_{\mathcal C^j\subseteq G'} C_j,
$$

where $ G' $ is the commutator subgroup of $ G $, $(\mathcal C^j)$ are the conjugacy classes of $ G $, and $ C_j := \sum_{g \in \mathcal{C}^j} g $ is the class sum associated with $\mathcal{C}^j$. It was already observed in \cite{har08}.
\er

A fusion ring $ R $ is termed \emph{$\alpha$-Frobenius} if, for every basis element $ x $, the expression $\frac{\mathrm{FPdim}(R)^{\alpha}}{\mathrm{FPdim}(x)}$ is an algebraic integer \cite{eno-nec}. Kaplansky's 6th conjecture posits that the Grothendieck ring of a complex fusion category is $1$-Frobenius \cite[Question 1]{eno-weakly}. An exhaustive classification of all $1$-Frobenius simple integral fusion rings, within certain specified limits, was provided in \cite{lpw}. These limits have been updated in \cite{BP24} as follows:

$$\begin{array}{c|cccccccc}
\text{Rank} & \le 5 & 6 & 7 & 8 & 9 & 10 & 11 & 12 \\ \hline
\FPdim \le  &  10^7 & 10^6 & 10^5 & 20000 & 10000 & 5000 & 3000 & 1000
\end{array}$$  
There are exactly $505$ non-pointed examples (including $8$ that are character rings of groups). Among them, we found only $4$ non-Burnside ones. They have the type of $\Rep(A_7)$, but different fusion data.
\begin{itemize}
\item Rank: $9$,
\item FPdim: $2520 = 2^3 \times 3^2 \times 5 \times 7$,
\item Type: $[[1, 1], [6, 1], [10, 2], [14, 2], [15, 1], [21, 1], [35, 1]]$,
\item Fusion data 1: 
$$\begin{smallmatrix}1&0&0&0&0&0&0&0&0 \\ 0&1&0&0&0&0&0&0&0 \\ 0&0&1&0&0&0&0&0&0 \\ 0&0&0&1&0&0&0&0&0 \\ 0&0&0&0&1&0&0&0&0 \\ 0&0&0&0&0&1&0&0&0 \\ 0&0&0&0&0&0&1&0&0 \\ 0&0&0&0&0&0&0&1&0 \\ 0&0&0&0&0&0&0&0&1\end{smallmatrix}, \ \begin{smallmatrix}0&1&0&0&0&0&0&0&0 \\ 1&0&0&0&0&1&0&1&0 \\ 0&0&0&1&0&0&1&0&1 \\ 0&0&1&0&0&0&1&0&1 \\ 0&0&0&0&1&2&0&2&0 \\ 0&1&0&0&2&0&1&0&1 \\ 0&0&1&1&0&1&0&1&1 \\ 0&1&0&0&2&0&1&2&1 \\ 0&0&1&1&0&1&1&1&4\end{smallmatrix}, \ \begin{smallmatrix}0&0&1&0&0&0&0&0&0 \\ 0&0&0&1&0&0&1&0&1 \\ 0&1&0&1&1&1&0&1&1 \\ 1&0&0&0&1&1&1&1&1 \\ 0&0&1&1&0&0&1&0&3 \\ 0&0&1&1&0&1&1&1&2 \\ 0&1&1&0&1&1&1&1&2 \\ 0&0&1&1&0&1&1&1&4 \\ 0&1&1&1&3&2&2&4&4\end{smallmatrix}, \ \begin{smallmatrix}0&0&0&1&0&0&0&0&0 \\ 0&0&1&0&0&0&1&0&1 \\ 1&0&0&0&1&1&1&1&1 \\ 0&1&1&0&1&1&0&1&1 \\ 0&0&1&1&0&0&1&0&3 \\ 0&0&1&1&0&1&1&1&2 \\ 0&1&0&1&1&1&1&1&2 \\ 0&0&1&1&0&1&1&1&4 \\ 0&1&1&1&3&2&2&4&4\end{smallmatrix}, \ \begin{smallmatrix}0&0&0&0&1&0&0&0&0 \\ 0&0&0&0&1&2&0&2&0 \\ 0&0&1&1&0&0&1&0&3 \\ 0&0&1&1&0&0&1&0&3 \\ 1&1&0&0&3&3&0&5&0 \\ 0&2&0&0&3&0&2&2&2 \\ 0&0&1&1&0&2&1&2&3 \\ 0&2&0&0&5&2&2&4&2 \\ 0&0&3&3&0&2&3&2&9\end{smallmatrix},$$
$$
\begin{smallmatrix}0&0&0&0&0&1&0&0&0 \\ 0&1&0&0&2&0&1&0&1 \\ 0&0&1&1&0&1&1&1&2 \\ 0&0&1&1&0&1&1&1&2 \\ 0&2&0&0&3&0&2&2&2 \\ 1&0&1&1&0&3&0&3&2 \\ 0&1&1&1&2&0&2&1&3 \\ 0&0&1&1&2&3&1&4&3 \\ 0&1&2&2&2&2&3&3&8\end{smallmatrix}, \ \begin{smallmatrix}0&0&0&0&0&0&1&0&0 \\ 0&0&1&1&0&1&0&1&1 \\ 0&1&1&0&1&1&1&1&2 \\ 0&1&0&1&1&1&1&1&2 \\ 0&0&1&1&0&2&1&2&3 \\ 0&1&1&1&2&0&2&1&3 \\ 1&0&1&1&1&2&1&2&3 \\ 0&1&1&1&2&1&2&2&5 \\ 0&1&2&2&3&3&3&5&7\end{smallmatrix}, \ \begin{smallmatrix}0&0&0&0&0&0&0&1&0 \\ 0&1&0&0&2&0&1&2&1 \\ 0&0&1&1&0&1&1&1&4 \\ 0&0&1&1&0&1&1&1&4 \\ 0&2&0&0&5&2&2&4&2 \\ 0&0&1&1&2&3&1&4&3 \\ 0&1&1&1&2&1&2&2&5 \\ 1&2&1&1&4&4&2&6&4 \\ 0&1&4&4&2&3&5&4&12\end{smallmatrix}, \ \begin{smallmatrix}0&0&0&0&0&0&0&0&1 \\ 0&0&1&1&0&1&1&1&4 \\ 0&1&1&1&3&2&2&4&4 \\ 0&1&1&1&3&2&2&4&4 \\ 0&0&3&3&0&2&3&2&9 \\ 0&1&2&2&2&2&3&3&8 \\ 0&1&2&2&3&3&3&5&7 \\ 0&1&4&4&2&3&5&4&12 \\ 1&4&4&4&9&8&7&12&15\end{smallmatrix}
$$
\item Fusion data 2: 
$$\begin{smallmatrix}1&0&0&0&0&0&0&0&0 \\ 0&1&0&0&0&0&0&0&0 \\ 0&0&1&0&0&0&0&0&0 \\ 0&0&0&1&0&0&0&0&0 \\ 0&0&0&0&1&0&0&0&0 \\ 0&0&0&0&0&1&0&0&0 \\ 0&0&0&0&0&0&1&0&0 \\ 0&0&0&0&0&0&0&1&0 \\ 0&0&0&0&0&0&0&0&1\end{smallmatrix}, \ \begin{smallmatrix}0&1&0&0&0&0&0&0&0 \\ 1&0&0&0&0&1&0&1&0 \\ 0&0&0&1&0&0&1&0&1 \\ 0&0&1&0&0&0&1&0&1 \\ 0&0&0&0&1&2&0&2&0 \\ 0&1&0&0&2&0&1&0&1 \\ 0&0&1&1&0&1&0&1&1 \\ 0&1&0&0&2&0&1&2&1 \\ 0&0&1&1&0&1&1&1&4\end{smallmatrix}, \ \begin{smallmatrix}0&0&1&0&0&0&0&0&0 \\ 0&0&0&1&0&0&1&0&1 \\ 0&1&1&2&1&1&1&1&0 \\ 1&0&1&1&1&1&2&1&0 \\ 0&0&1&1&0&0&1&0&3 \\ 0&0&1&1&0&1&1&1&2 \\ 0&1&2&1&1&1&2&1&1 \\ 0&0&1&1&0&1&1&1&4 \\ 0&1&0&0&3&2&1&4&5\end{smallmatrix}, \ \begin{smallmatrix}0&0&0&1&0&0&0&0&0 \\ 0&0&1&0&0&0&1&0&1 \\ 1&0&1&1&1&1&2&1&0 \\ 0&1&2&1&1&1&1&1&0 \\ 0&0&1&1&0&0&1&0&3 \\ 0&0&1&1&0&1&1&1&2 \\ 0&1&1&2&1&1&2&1&1 \\ 0&0&1&1&0&1&1&1&4 \\ 0&1&0&0&3&2&1&4&5\end{smallmatrix}, \ \begin{smallmatrix}0&0&0&0&1&0&0&0&0 \\ 0&0&0&0&1&2&0&2&0 \\ 0&0&1&1&0&0&1&0&3 \\ 0&0&1&1&0&0&1&0&3 \\ 1&1&0&0&3&3&0&5&0 \\ 0&2&0&0&3&0&2&2&2 \\ 0&0&1&1&0&2&1&2&3 \\ 0&2&0&0&5&2&2&4&2 \\ 0&0&3&3&0&2&3&2&9\end{smallmatrix},$$
$$\begin{smallmatrix}0&0&0&0&0&1&0&0&0 \\ 0&1&0&0&2&0&1&0&1 \\ 0&0&1&1&0&1&1&1&2 \\ 0&0&1&1&0&1&1&1&2 \\ 0&2&0&0&3&0&2&2&2 \\ 1&0&1&1&0&3&0&3&2 \\ 0&1&1&1&2&0&2&1&3 \\ 0&0&1&1&2&3&1&4&3 \\ 0&1&2&2&2&2&3&3&8\end{smallmatrix}, \ \begin{smallmatrix}0&0&0&0&0&0&1&0&0 \\ 0&0&1&1&0&1&0&1&1 \\ 0&1&2&1&1&1&2&1&1 \\ 0&1&1&2&1&1&2&1&1 \\ 0&0&1&1&0&2&1&2&3 \\ 0&1&1&1&2&0&2&1&3 \\ 1&0&2&2&1&2&2&2&2 \\ 0&1&1&1&2&1&2&2&5 \\ 0&1&1&1&3&3&2&5&8\end{smallmatrix}, \ \begin{smallmatrix}0&0&0&0&0&0&0&1&0 \\ 0&1&0&0&2&0&1&2&1 \\ 0&0&1&1&0&1&1&1&4 \\ 0&0&1&1&0&1&1&1&4 \\ 0&2&0&0&5&2&2&4&2 \\ 0&0&1&1&2&3&1&4&3 \\ 0&1&1&1&2&1&2&2&5 \\ 1&2&1&1&4&4&2&6&4 \\ 0&1&4&4&2&3&5&4&12\end{smallmatrix}, \ \begin{smallmatrix}0&0&0&0&0&0&0&0&1 \\ 0&0&1&1&0&1&1&1&4 \\ 0&1&0&0&3&2&1&4&5 \\ 0&1&0&0&3&2&1&4&5 \\ 0&0&3&3&0&2&3&2&9 \\ 0&1&2&2&2&2&3&3&8 \\ 0&1&1&1&3&3&2&5&8 \\ 0&1&4&4&2&3&5&4&12 \\ 1&4&5&5&9&8&8&12&14\end{smallmatrix}$$
\item Fusion data 3: 
$$\begin{smallmatrix}1&0&0&0&0&0&0&0&0 \\ 0&1&0&0&0&0&0&0&0 \\ 0&0&1&0&0&0&0&0&0 \\ 0&0&0&1&0&0&0&0&0 \\ 0&0&0&0&1&0&0&0&0 \\ 0&0&0&0&0&1&0&0&0 \\ 0&0&0&0&0&0&1&0&0 \\ 0&0&0&0&0&0&0&1&0 \\ 0&0&0&0&0&0&0&0&1\end{smallmatrix}, \ \begin{smallmatrix}0&1&0&0&0&0&0&0&0 \\ 1&1&0&0&0&1&1&0&0 \\ 0&0&0&1&0&0&1&0&1 \\ 0&0&1&0&0&0&1&0&1 \\ 0&0&0&0&2&0&0&1&1 \\ 0&1&0&0&0&3&1&1&0 \\ 0&1&1&1&0&1&1&0&1 \\ 0&0&0&0&1&1&0&3&1 \\ 0&0&1&1&1&0&1&1&4\end{smallmatrix}, \ \begin{smallmatrix}0&0&1&0&0&0&0&0&0 \\ 0&0&0&1&0&0&1&0&1 \\ 0&1&0&1&0&1&0&0&2 \\ 1&0&0&0&0&1&1&0&2 \\ 0&0&0&0&1&1&0&2&2 \\ 0&0&1&1&1&0&1&1&2 \\ 0&1&1&0&0&1&1&0&3 \\ 0&0&0&0&2&1&0&3&3 \\ 0&1&2&2&2&2&3&3&4\end{smallmatrix}, \ \begin{smallmatrix}0&0&0&1&0&0&0&0&0 \\ 0&0&1&0&0&0&1&0&1 \\ 1&0&0&0&0&1&1&0&2 \\ 0&1&1&0&0&1&0&0&2 \\ 0&0&0&0&1&1&0&2&2 \\ 0&0&1&1&1&0&1&1&2 \\ 0&1&0&1&0&1&1&0&3 \\ 0&0&0&0&2&1&0&3&3 \\ 0&1&2&2&2&2&3&3&4\end{smallmatrix}, \ \begin{smallmatrix}0&0&0&0&1&0&0&0&0 \\ 0&0&0&0&2&0&0&1&1 \\ 0&0&0&0&1&1&0&2&2 \\ 0&0&0&0&1&1&0&2&2 \\ 1&2&1&1&1&1&2&0&3 \\ 0&0&1&1&1&0&1&2&3 \\ 0&0&0&0&2&1&0&3&3 \\ 0&1&2&2&0&2&3&0&5 \\ 0&1&2&2&3&3&3&5&6\end{smallmatrix},$$
$$\begin{smallmatrix}0&0&0&0&0&1&0&0&0 \\ 0&1&0&0&0&3&1&1&0 \\ 0&0&1&1&1&0&1&1&2 \\ 0&0&1&1&1&0&1&1&2 \\ 0&0&1&1&1&0&1&2&3 \\ 1&3&0&0&0&6&2&3&0 \\ 0&1&1&1&1&2&2&2&2 \\ 0&1&1&1&2&3&2&3&3 \\ 0&0&2&2&3&0&2&3&9\end{smallmatrix}, \ \begin{smallmatrix}0&0&0&0&0&0&1&0&0 \\ 0&1&1&1&0&1&1&0&1 \\ 0&1&1&0&0&1&1&0&3 \\ 0&1&0&1&0&1&1&0&3 \\ 0&0&0&0&2&1&0&3&3 \\ 0&1&1&1&1&2&2&2&2 \\ 1&1&1&1&0&2&2&0&4 \\ 0&0&0&0&3&2&0&5&4 \\ 0&1&3&3&3&2&4&4&7\end{smallmatrix}, \ \begin{smallmatrix}0&0&0&0&0&0&0&1&0 \\ 0&0&0&0&1&1&0&3&1 \\ 0&0&0&0&2&1&0&3&3 \\ 0&0&0&0&2&1&0&3&3 \\ 0&1&2&2&0&2&3&0&5 \\ 0&1&1&1&2&3&2&3&3 \\ 0&0&0&0&3&2&0&5&4 \\ 1&3&3&3&0&3&5&0&7 \\ 0&1&3&3&5&3&4&7&10\end{smallmatrix}, \ \begin{smallmatrix}0&0&0&0&0&0&0&0&1 \\ 0&0&1&1&1&0&1&1&4 \\ 0&1&2&2&2&2&3&3&4 \\ 0&1&2&2&2&2&3&3&4 \\ 0&1&2&2&3&3&3&5&6 \\ 0&0&2&2&3&0&2&3&9 \\ 0&1&3&3&3&2&4&4&7 \\ 0&1&3&3&5&3&4&7&10 \\ 1&4&4&4&6&9&7&10&17\end{smallmatrix}$$
\item Fusion data 4: 
$$\begin{smallmatrix}1&0&0&0&0&0&0&0&0 \\ 0&1&0&0&0&0&0&0&0 \\ 0&0&1&0&0&0&0&0&0 \\ 0&0&0&1&0&0&0&0&0 \\ 0&0&0&0&1&0&0&0&0 \\ 0&0&0&0&0&1&0&0&0 \\ 0&0&0&0&0&0&1&0&0 \\ 0&0&0&0&0&0&0&1&0 \\ 0&0&0&0&0&0&0&0&1\end{smallmatrix}, \ \begin{smallmatrix}0&1&0&0&0&0&0&0&0 \\ 1&0&0&0&0&1&0&1&0 \\ 0&0&0&1&0&0&1&0&1 \\ 0&0&1&0&0&0&1&0&1 \\ 0&0&0&0&1&2&0&2&0 \\ 0&1&0&0&2&0&1&0&1 \\ 0&0&1&1&0&1&0&1&1 \\ 0&1&0&0&2&0&1&2&1 \\ 0&0&1&1&0&1&1&1&4\end{smallmatrix}, \ \begin{smallmatrix}0&0&1&0&0&0&0&0&0 \\ 0&0&0&1&0&0&1&0&1 \\ 1&0&0&2&1&1&2&1&0 \\ 0&1&2&1&1&1&1&1&0 \\ 0&0&1&1&0&0&1&0&3 \\ 0&0&1&1&0&1&1&1&2 \\ 0&1&2&1&1&1&2&1&1 \\ 0&0&1&1&0&1&1&1&4 \\ 0&1&0&0&3&2&1&4&5\end{smallmatrix}, \ \begin{smallmatrix}0&0&0&1&0&0&0&0&0 \\ 0&0&1&0&0&0&1&0&1 \\ 0&1&2&1&1&1&1&1&0 \\ 1&0&1&1&1&1&2&1&0 \\ 0&0&1&1&0&0&1&0&3 \\ 0&0&1&1&0&1&1&1&2 \\ 0&1&1&2&1&1&2&1&1 \\ 0&0&1&1&0&1&1&1&4 \\ 0&1&0&0&3&2&1&4&5\end{smallmatrix}, \ \begin{smallmatrix}0&0&0&0&1&0&0&0&0 \\ 0&0&0&0&1&2&0&2&0 \\ 0&0&1&1&0&0&1&0&3 \\ 0&0&1&1&0&0&1&0&3 \\ 1&1&0&0&3&3&0&5&0 \\ 0&2&0&0&3&0&2&2&2 \\ 0&0&1&1&0&2&1&2&3 \\ 0&2&0&0&5&2&2&4&2 \\ 0&0&3&3&0&2&3&2&9\end{smallmatrix},$$
$$\begin{smallmatrix}0&0&0&0&0&1&0&0&0 \\ 0&1&0&0&2&0&1&0&1 \\ 0&0&1&1&0&1&1&1&2 \\ 0&0&1&1&0&1&1&1&2 \\ 0&2&0&0&3&0&2&2&2 \\ 1&0&1&1&0&3&0&3&2 \\ 0&1&1&1&2&0&2&1&3 \\ 0&0&1&1&2&3&1&4&3 \\ 0&1&2&2&2&2&3&3&8\end{smallmatrix}, \ \begin{smallmatrix}0&0&0&0&0&0&1&0&0 \\ 0&0&1&1&0&1&0&1&1 \\ 0&1&2&1&1&1&2&1&1 \\ 0&1&1&2&1&1&2&1&1 \\ 0&0&1&1&0&2&1&2&3 \\ 0&1&1&1&2&0&2&1&3 \\ 1&0&2&2&1&2&2&2&2 \\ 0&1&1&1&2&1&2&2&5 \\ 0&1&1&1&3&3&2&5&8\end{smallmatrix}, \ \begin{smallmatrix}0&0&0&0&0&0&0&1&0 \\ 0&1&0&0&2&0&1&2&1 \\ 0&0&1&1&0&1&1&1&4 \\ 0&0&1&1&0&1&1&1&4 \\ 0&2&0&0&5&2&2&4&2 \\ 0&0&1&1&2&3&1&4&3 \\ 0&1&1&1&2&1&2&2&5 \\ 1&2&1&1&4&4&2&6&4 \\ 0&1&4&4&2&3&5&4&12\end{smallmatrix}, \ \begin{smallmatrix}0&0&0&0&0&0&0&0&1 \\ 0&0&1&1&0&1&1&1&4 \\ 0&1&0&0&3&2&1&4&5 \\ 0&1&0&0&3&2&1&4&5 \\ 0&0&3&3&0&2&3&2&9 \\ 0&1&2&2&2&2&3&3&8 \\ 0&1&1&1&3&3&2&5&8 \\ 0&1&4&4&2&3&5&4&12 \\ 1&4&5&5&9&8&8&12&14\end{smallmatrix}$$
\end{itemize}
The determinant of the second fusion matrix for each fusion data set mentioned above is $\pm 36$, while its norm is $6$. This indicates that these integral fusion rings are non-Burnside and, consequently, cannot be categorified.

Similar classifications have been conducted for the non-1-Frobenius case. Interestingly, most—but not all—of the fusion rings identified are non-Burnside, and thus, not amenable to categorification (see the example below). This serves as supporting evidence for Kaplansky's 6th conjecture.

\begin{itemize}
\item Rank: $7$,
\item FPdim: $798 = 2 \times 3 \times 7 \times 19$,
\item Type: $[[1, 1], [7, 1], [8,1], [9,3], [21, 1]]$,
\item Fusion data: $$ 
\begin{smallmatrix}1&0&0&0&0&0&0 \\ 0&1&0&0&0&0&0 \\ 0&0&1&0&0&0&0 \\ 0&0&0&1&0&0&0 \\ 0&0&0&0&1&0&0 \\ 0&0&0&0&0&1&0 \\ 0&0&0&0&0&0&1\end{smallmatrix}, \
\begin{smallmatrix}0&1&0&0&0&0&0 \\ 1&0&0&1&1&1&1 \\ 0&0&1&1&1&1&1 \\ 0&1&1&1&1&1&1 \\ 0&1&1&1&1&1&1 \\ 0&1&1&1&1&1&1 \\ 0&1&1&1&1&1&5\end{smallmatrix}, \
\begin{smallmatrix}0&0&1&0&0&0&0 \\ 0&0&1&1&1&1&1 \\ 1&1&1&1&1&1&1 \\ 0&1&1&2&1&1&1 \\ 0&1&1&1&2&1&1 \\ 0&1&1&1&1&2&1 \\ 0&1&1&1&1&1&6\end{smallmatrix}, \
\begin{smallmatrix}0&0&0&1&0&0&0 \\ 0&1&1&1&1&1&1 \\ 0&1&1&2&1&1&1 \\ 1&1&2&1&1&2&1 \\ 0&1&1&1&2&2&1 \\ 0&1&1&2&2&1&1 \\ 0&1&1&1&1&1&7\end{smallmatrix}, \
\begin{smallmatrix}0&0&0&0&1&0&0 \\ 0&1&1&1&1&1&1 \\ 0&1&1&1&2&1&1 \\ 0&1&1&1&2&2&1 \\ 1&1&2&2&1&1&1 \\ 0&1&1&2&1&2&1 \\ 0&1&1&1&1&1&7\end{smallmatrix}, \
\begin{smallmatrix}0&0&0&0&0&1&0 \\ 0&1&1&1&1&1&1 \\ 0&1&1&1&1&2&1 \\ 0&1&1&2&2&1&1 \\ 0&1&1&2&1&2&1 \\ 1&1&2&1&2&1&1 \\ 0&1&1&1&1&1&7\end{smallmatrix}, \
\begin{smallmatrix}0&0&0&0&0&0&1 \\ 0&1&1&1&1&1&5 \\ 0&1&1&1&1&1&6 \\ 0&1&1&1&1&1&7 \\ 0&1&1&1&1&1&7 \\ 0&1&1&1&1&1&7 \\ 1&5&6&7&7&7&8 \end{smallmatrix}
$$ 
\end{itemize}
The determinant of the third fusion matrix is $16$, while its norm is $8$. This shows that the integral fusion ring is non-Burnside and therefore cannot be categorified.
\subsection{Dual-Burnside property} \label{sub:dualburn}
We abbreviate ``(A)RN" for ``(abelian) real non-negative." According to \cite[Corollary 8.5]{lpw}, the Grothendieck ring of a unitary fusion category, if commutative, has an RN dual (hypergroup).

\begin{question} Is the dual of the Grothendieck ring of a complex fusion category always RN?
\end{question}

Recall that a commutative fusion ring $ F $ is called \emph{dual-Burnside} if a column of its character table has a zero entry if and only if the squared norm of the column (the corresponding formal codegree) is strictly less than $\mathrm{FPdim}(F)$. When $ F $ has an RN dual, it is equivalent to $ P^2 $ being idempotent, where
$$
P = \prod_{i \in \mathcal{I}} \frac{x_i}{\mathrm{FPdim}(x_i)}
$$
and $\{x_i\}_{i \in I}$ is the set of basis elements (see Corollary \ref{corP2DualB}). Note that for a non-abelian group $ G $, the product of all its elements depends on the ordering; in fact, the set of all possible products forms a $ G' $-coset (see \cite{DeHe}).

According to \cite[Theorem B]{inw}, for any finite nilpotent group $ G $, the Grothendieck ring of $\Rep(G)$ is dual-Burnside. This extends to every nilpotent dualizable ARN-hypergroup by Theorem \ref{nilpotent:burnside}. Thus, every commutative nilpotent fusion ring with an RN dual is both Burnside and dual-Burnside. All commutative nilpotent fusion rings we examined have RN duals.

\begin{question} Is there a commutative nilpotent fusion ring without an RN dual?
\end{question}

We define a finite group $ G $ as \emph{dual-Burnside} if $ K(\Rep(G)) $ is dual-Burnside. Therefore, a nilpotent finite group is dual-Burnside. There are exactly 144 finite groups of order less than 32, among which 30 are non-nilpotent. Of these, only $\mathrm{SL}(2,3)$ is dual-Burnside, while the other 29, such as $ S_3 $ and $ D_5 $, are not.

A finite group $ G $ is centerless if and only if $ K(\Rep(G)) $ has a perfect dual (see Definition \ref{perfect}, Theorem \ref{univ:grd}, and Lemma \ref{lem:centergrp}). Thus, a centerless finite group is dual-Burnside if every non-$\mathrm{FPdim}$ column of its character table has a zero entry. We verified that every centerless and dual-Burnside finite group $ G $ of order $|G| \le 1000$ is \emph{almost simple} (i.e., $ S \subseteq G \subseteq \Aut(S) $ with $ S $ non-abelian simple); however, $ G = A_5 \times A_5 $, of order 3600, is not almost simple. A finite non-abelian simple group is centerless, so per \cite{PalMO} and its answers, all non-alternating finite simple groups are dual-Burnside, except the Mathieu groups $ M_{22} $ and $ M_{24} $. The alternating group $ A_n $, for $ 5 \le n \le 19 $, is dual-Burnside if and only if $ n \not \in \{7, 11, 13, 15, 16, 18, 19\} $. Recall that a finite group $ G $ is simple if and only if $ K(\Rep(G)) $ is simple.

\br
By Theorem \ref{hbz:hyp}, a finite group $ G $ is dual-Burnside if and only if

$$
\left(\prod_{\chi \in \mathrm{Irr}(G)} \frac{\chi}{\chi(1)} \right)^2 = \frac{|Z(G)|}{|G|} \left(\sum_{\chi \in \mathrm{Irr}(G/Z(G))}\chi(1) \chi \right).
$$

\er

\subsection{Modular categories} \label{mtc:app}
Let us present some applications of Theorem \ref{first:div}. There are exactly 71 distinct types of half-Frobenius integral fusion rings with ranks up to 12 (see \cite{abpp}). Of these, 15 types (listed below) are excluded from modular categorification by Theorems \ref{first:div} or \ref{ccpt}.
{\scriptsize
$$[1, 1, 1, 1, 2, 2],
 [1, 1, 1, 1, 2, 2, 2, 2],
 [1, 1, 1, 1, 2, 2, 2, 2, 2],
 [1, 1, 1, 1, 2, 2, 2, 4, 4],
 [1, 1, 1, 1, 1, 1, 1, 1, 1, 3],$$ $$
 [1, 1, 1, 1, 2, 2, 2, 2, 2, 2],
 [1, 1, 1, 1, 2, 2, 2, 4, 4, 4, 4], 
 [1, 1, 1, 1, 1, 1, 1, 1, 1, 3, 3, 3],
 [1, 1, 1, 1, 1, 1, 1, 1, 2, 2, 2, 2], $$ $$
 [1, 1, 1, 1, 1, 1, 1, 1, 2, 2, 4, 4], 
 [1, 1, 1, 1, 2, 2, 2, 2, 2, 2, 2, 2],
 [1, 1, 1, 1, 2, 2, 2, 2, 2, 2, 2, 4],
 [1, 1, 1, 1, 2, 2, 2, 4, 4, 4, 4, 4],  $$ $$
 [1, 1, 1, 1, 2, 2, 2, 4, 4, 4, 8, 8],
 [1, 1, 2, 2, 2, 2, 3, 3, 6, 6, 6, 6].$$}

For example, let us provide all the fusion data for the first type above:  
\begin{itemize}
\item Rank: $6$,
\item FPdim: $12 = 2^2 \times 3$,
\item Type: $[1, 1, 1, 1, 2, 2]$,
\item Fusion data 1: 
$$ \begin{smallmatrix}1&0&0&0&0&0 \\ 0&1&0&0&0&0 \\ 0&0&1&0&0&0 \\ 0&0&0&1&0&0 \\ 0&0&0&0&1&0 \\ 0&0&0&0&0&1\end{smallmatrix}, \
\begin{smallmatrix}0&1&0&0&0&0 \\ 1&0&0&0&0&0 \\ 0&0&0&1&0&0 \\ 0&0&1&0&0&0 \\ 0&0&0&0&0&1 \\ 0&0&0&0&1&0\end{smallmatrix}, \
\begin{smallmatrix}0&0&1&0&0&0 \\ 0&0&0&1&0&0 \\ 1&0&0&0&0&0 \\ 0&1&0&0&0&0 \\ 0&0&0&0&0&1 \\ 0&0&0&0&1&0\end{smallmatrix}, \
\begin{smallmatrix}0&0&0&1&0&0 \\ 0&0&1&0&0&0 \\ 0&1&0&0&0&0 \\ 1&0&0&0&0&0 \\ 0&0&0&0&1&0 \\ 0&0&0&0&0&1\end{smallmatrix}, \
\begin{smallmatrix}0&0&0&0&1&0 \\ 0&0&0&0&0&1 \\ 0&0&0&0&0&1 \\ 0&0&0&0&1&0 \\ 1&0&0&1&0&1 \\ 0&1&1&0&1&0\end{smallmatrix}, \
\begin{smallmatrix}0&0&0&0&0&1 \\ 0&0&0&0&1&0 \\ 0&0&0&0&1&0 \\ 0&0&0&0&0&1 \\ 0&1&1&0&1&0 \\ 1&0&0&1&0&1\end{smallmatrix}$$
\item Fusion data 2:
$$ \begin{smallmatrix}1&0&0&0&0&0 \\ 0&1&0&0&0&0 \\ 0&0&1&0&0&0 \\ 0&0&0&1&0&0 \\ 0&0&0&0&1&0 \\ 0&0&0&0&0&1\end{smallmatrix}, \
\begin{smallmatrix}0&1&0&0&0&0 \\ 1&0&0&0&0&0 \\ 0&0&0&1&0&0 \\ 0&0&1&0&0&0 \\ 0&0&0&0&1&0 \\ 0&0&0&0&0&1\end{smallmatrix}, \
\begin{smallmatrix}0&0&1&0&0&0 \\ 0&0&0&1&0&0 \\ 0&1&0&0&0&0 \\ 1&0&0&0&0&0 \\ 0&0&0&0&0&1 \\ 0&0&0&0&1&0\end{smallmatrix}, \
\begin{smallmatrix}0&0&0&1&0&0 \\ 0&0&1&0&0&0 \\ 1&0&0&0&0&0 \\ 0&1&0&0&0&0 \\ 0&0&0&0&0&1 \\ 0&0&0&0&1&0\end{smallmatrix}, \
\begin{smallmatrix}0&0&0&0&1&0 \\ 0&0&0&0&1&0 \\ 0&0&0&0&0&1 \\ 0&0&0&0&0&1 \\ 1&1&0&0&0&1 \\ 0&0&1&1&1&0\end{smallmatrix}, \
\begin{smallmatrix}0&0&0&0&0&1 \\ 0&0&0&0&0&1 \\ 0&0&0&0&1&0 \\ 0&0&0&0&1&0 \\ 0&0&1&1&1&0 \\ 1&1&0&0&0&1\end{smallmatrix}$$
\item Fusion data 3:
$$ \begin{smallmatrix}1&0&0&0&0&0 \\ 0&1&0&0&0&0 \\ 0&0&1&0&0&0 \\ 0&0&0&1&0&0 \\ 0&0&0&0&1&0 \\ 0&0&0&0&0&1\end{smallmatrix}, \
\begin{smallmatrix}0&1&0&0&0&0 \\ 1&0&0&0&0&0 \\ 0&0&0&1&0&0 \\ 0&0&1&0&0&0 \\ 0&0&0&0&1&0 \\ 0&0&0&0&0&1\end{smallmatrix}, \
\begin{smallmatrix}0&0&1&0&0&0 \\ 0&0&0&1&0&0 \\ 1&0&0&0&0&0 \\ 0&1&0&0&0&0 \\ 0&0&0&0&1&0 \\ 0&0&0&0&0&1\end{smallmatrix}, \
\begin{smallmatrix}0&0&0&1&0&0 \\ 0&0&1&0&0&0 \\ 0&1&0&0&0&0 \\ 1&0&0&0&0&0 \\ 0&0&0&0&1&0 \\ 0&0&0&0&0&1\end{smallmatrix}, \
\begin{smallmatrix}0&0&0&0&1&0 \\ 0&0&0&0&1&0 \\ 0&0&0&0&1&0 \\ 0&0&0&0&1&0 \\ 0&0&0&0&0&2 \\ 1&1&1&1&0&0\end{smallmatrix}, \
\begin{smallmatrix}0&0&0&0&0&1 \\ 0&0&0&0&0&1 \\ 0&0&0&0&0&1 \\ 0&0&0&0&0&1 \\ 1&1&1&1&0&0 \\ 0&0&0&0&2&0\end{smallmatrix}$$
\item Fusion data 4:
$$ \begin{smallmatrix}1&0&0&0&0&0 \\ 0&1&0&0&0&0 \\ 0&0&1&0&0&0 \\ 0&0&0&1&0&0 \\ 0&0&0&0&1&0 \\ 0&0&0&0&0&1\end{smallmatrix}, \
\begin{smallmatrix}0&1&0&0&0&0 \\ 1&0&0&0&0&0 \\ 0&0&0&1&0&0 \\ 0&0&1&0&0&0 \\ 0&0&0&0&1&0 \\ 0&0&0&0&0&1\end{smallmatrix}, \
\begin{smallmatrix}0&0&1&0&0&0 \\ 0&0&0&1&0&0 \\ 0&1&0&0&0&0 \\ 1&0&0&0&0&0 \\ 0&0&0&0&1&0 \\ 0&0&0&0&0&1\end{smallmatrix}, \
\begin{smallmatrix}0&0&0&1&0&0 \\ 0&0&1&0&0&0 \\ 1&0&0&0&0&0 \\ 0&1&0&0&0&0 \\ 0&0&0&0&1&0 \\ 0&0&0&0&0&1\end{smallmatrix}, \
\begin{smallmatrix}0&0&0&0&1&0 \\ 0&0&0&0&1&0 \\ 0&0&0&0&1&0 \\ 0&0&0&0&1&0 \\ 0&0&0&0&0&2 \\ 1&1&1&1&0&0\end{smallmatrix}, \
\begin{smallmatrix}0&0&0&0&0&1 \\ 0&0&0&0&0&1 \\ 0&0&0&0&0&1 \\ 0&0&0&0&0&1 \\ 1&1&1&1&0&0 \\ 0&0&0&0&2&0\end{smallmatrix}$$
\end{itemize}
They are excluded from modular categorification by Theorem \ref{first:div} (\ref{prime:set:int}) because $\mtc{V}(\mathcal{C}) = \mtc{V}(12) = \{2,3\}$, while $3 \not\in \mtc{V}(d_i)$ for any $i$, and $\mtc{V}(\mathcal{C}_{pt}) = \{2\}$. Additionally, Theorem \ref{ccpt} rules them out since 3 is a (powerless) factor of $\fp(\mathcal{C})$, but does not divide $\fp(\mathcal{C}_{pt}) = 4$. 

Note that the fusion data for cases 3 and 4 can be generalized. Consider a group $G$ with order $n^2$ and a group $K$ with order $m+1$. Then, analyze the following half-Frobenius integral fusion ring:
\begin{itemize}
\item Rank: $n^2 + m$,
\item FPdim: $n^2 (m+1)$,
\item Type: $[[1,n^2],[n,m]]$,
\item Basis: $\mathcal{B} = \{x_g\}_{g \in G} \cup \{\rho_k\}_{k \in K \setminus \{ e \}}$,
\item Fusion data:
\begin{itemize}
\item $\rho_e := \frac{1}{n}\sum_{g \in G} x_g$,
\item $x_g x_h = x_{gh}$, for all $g,h \in G$,
\item $\rho_k \rho_l = n\rho_{kl}$, for all $k,l \in K$,
\item $x_g \rho_k = \rho_k x_g = \rho_k$, for all $(g,k) \in G \times K$,
\end{itemize}
\end{itemize}
According to Theorem \ref{first:div} (\ref{prime:set:int}), it is excluded from modular categorification if $n > 1$ and $m+1$ has a prime factor that does not divide $n$.

Note that the types $[1, 1, 1, 1, 1, 1, 1, 1, 1, 3, 3, 3]$ and $[1, 1, 1, 1, 2, 2, 2, 2, 2, 2, 2, 2]$ are excluded from modular categorification by Theorem \ref{first:div}, but not by Theorem \ref{ccpt}, since their $\fp = 2^2 3^2$ have no powerless prime factor.

\vspace*{.5cm}

\noindent \textbf{Availability of data and materials.} Data for the computations in this paper are available on reasonable request from the authors.

\vspace*{.15cm}

\noindent \textbf{Conflict of interest statement.} On behalf of all authors, the corresponding author declares that there are no conflicts of interest.

\bibliographystyle{alpha}
\bibliography{24nov}
\ed